\documentclass[english,10pt]{article}
\usepackage[T1]{fontenc}
\usepackage[latin9]{inputenc}
\usepackage{geometry}
\usepackage{textcomp}
\usepackage{amsthm}
\usepackage{amsmath}
\usepackage{amssymb}
\usepackage[all]{xy}

\theoremstyle{plain}

\usepackage{mathrsfs}
\usepackage{hyperref}
\usepackage{xy}
\hypersetup{colorlinks=true,citecolor=blue}

\newcommand{\scx}{\mathscr{X}}
\newcommand{\sx}{\mathscr{X}}

\newcommand{\se}{\mathscr{E}}
\newcommand{\scz}{\mathscr{Z}}
\newcommand{\sz}{\mathscr{Z}}
\newcommand{\scw}{\mathscr{W}}
\newcommand{\sw}{\mathscr{W}}
\newcommand{\scm}{\mathscr{M}}
\newcommand{\sm}{\mathscr{M}}

\newcommand{\sj}{\mathscr{J}}
\newcommand{\st}{\mathscr{T}}
\newcommand{\scc}{\mathscr{C}}
\newcommand{\sco}{\mathscr{O}}
\newcommand{\scd}{\mathscr{D}}
\newcommand{\sca}{\mathscr{A}}
\newcommand{\scl}{\mathscr{L}}
\newcommand{\srr}{\mathscr{R}}
\newcommand{\ppar}{\mathscr{P} \mathrm{ar}}

\usepackage{babel}

\begin{document}

\title{Universal eigenvarieties, trianguline Galois representations, and
\emph{p-}adic Langlands functoriality}

\author{by David Hansen%
\thanks{Department of Mathematics, Columbia University, 2990 Broadway, New
York NY 10027; hansen@math.columbia.edu%
}\\
with an appendix by James Newton%
\thanks{Department of Mathematics, Imperial College London, London, SW7 2AZ;
j.newton@imperial.ac.uk%
}}
\maketitle
\begin{abstract}
Using the overconvergent cohomology modules introduced by Ash and
Stevens, we construct eigenvarieties associated with reductive groups
and establish some basic geometric properties of these spaces, building
on work of Ash-Stevens, Urban, and others. We also formulate a precise
modularity conjecture linking trianguline Galois representations with
overconvergent cohomology classes. In the course of giving evidence
for this conjecture, we establish several new instances of \emph{p-}adic
Langlands functoriality. Our main technical innovations are a family
of universal coefficients spectral sequences for overconvergent cohomology
and a generalization of Chenevier's interpolation theorem.

\tableofcontents{}
\end{abstract}

\section{Introduction}

\subsection{Eigenvarieties and overconvergent cohomology}

Since the pioneering works of Serre \cite{SeAntwerp}, Katz \cite{Katzmodschemes},
and especially Hida \cite{Hidordinv,Hidtotreal,HidordGL2} and Coleman \cite{Colemanclassicaloc,Colemanbanachfamilies},
\emph{p-}adic families of modular forms have become a major topic
in modern number theory. Aside from their intrinsic beauty, these
families have found applications towards Iwasawa theory, the Bloch-Kato
conjecture, modularity lifting theorems, and the local and global
Langlands correspondences \cite{BCnontemp,BCast,BCsign,BLGGT,Sknote,Wilesiwasawatotreal}.
One of the guiding examples in the field is Coleman and Mazur's\emph{
eigencurve} \cite{CMeigencurve,BuEigen}, a universal object parametrizing
all overconvergent \emph{p-}adic modular forms of fixed tame level
and finite slope. Concurrently with their work, Stevens introduced
his beautifully simple idea of \emph{overconvergent cohomology} \cite{Strigidsymbs},
a group-cohomological avatar of overconvergent \emph{p}-adic modular
forms\emph{.} 

Ash and Stevens developed this idea much further in \cite{AS}: as
conceived there, overconvergent cohomology works for any connected
reductive $\mathbf{Q}$-group $\mathbf{G}$ split at $p$, and leads
to natural candidates for quite general eigenvarieties. When the group
$\mathbf{G}^{\mathrm{der}}(\mathbf{R})$ possesses discrete series
representations, Urban \cite{UrEigen} used overconvergent cohomology
to construct eigenvarieties interpolating classical forms occuring
with nonzero Euler-Poincaré multiplicities, and showed that his construction
yields spaces which are equidimensional of the same dimension as weight
space. In this article, we develop the theory of eigenvarieties for
a connected reductive group over a number field $F$, building on
the ideas introduced in \cite{AS} and \cite{UrEigen}, and we formulate
precise conjectures relating these spaces with representations of
the absolute Galois group $\mathrm{Gal}(\overline{F}/F)$.

In order to state our results, we first establish some notation. Fix
a prime $p$ and a number field $ $$F$ with ring of integers $\mathcal{O}_{F}$.
For any place $v$ of $F$ we write $\mathcal{O}_{v}$ for the $v$-adic
completion of $\mathcal{O}_{F}$ and $F_{v}$ for its fraction field,
and we set $F_{\infty}=F\otimes_{\mathbf{Q}}\mathbf{R}$. Fix once
and for all an algebraic closure $\overline{\mathbf{Q}_{p}}$ and
an isomorphism $\iota:\mathbf{C}\overset{\sim}{\to}\overline{\mathbf{Q}_{p}}$.
Let $\mathbf{G}$ be a connected reductive group over $\mathbf{Q}$
of the form $\mathbf{G}=\mathrm{Res}_{F/\mathbf{Q}}\mathbf{H}$, with
$\mathbf{H}$ a connected reductive group over $F$ split over $F_{v}$
for each $v|p$. For any open compact subgroup $K_{f}\subset\mathbf{G}(\mathbf{A}_{f})=\mathbf{H}(\mathbf{A}_{F,f})$,
we have the associated locally symmetric space\[
Y(K_{f})=\mathbf{G}(\mathbf{Q})\backslash\mathbf{G}(\mathbf{A})/K_{\infty}^{\circ}K_{f},\]
where $K_{\infty}^{\circ}$ denotes the identity component of a maximal
compact-modulo-center subgroup $K_{\infty}\subset\mathbf{G}(\mathbf{R})=\mathbf{H}(F_{\infty})$.

For each place $v|p$, we set $H_{v}=\mathbf{H}\times_{\mathrm{Spec}F}\mathrm{Spec}F_{v}$.
Let $(B_{v},T_{v})$ be a split Borel pair in $H_{v}$. The group
$H_{v}$ spreads out to a group scheme over $\mathrm{Spec}\mathcal{O}_{v}$
(also denoted $H_{v}$) with reductive special fiber, and we assume
that $B_{v}$ and $T_{v}$ are defined compatibly over $\mathcal{O}_{v}$
as well. Set $G=\prod_{v|p}\mathrm{Res}_{\mathcal{O}_{v}/\mathbf{Z}_{p}}H_{v}$,
and define $B$ and $T$ analogously; we regard $G$, $B$ and $T$
as group schemes over $\mathrm{Spec}\mathbf{Z}_{p}$, so e.g. $G(R)=\prod_{v|p}H_{v}(\mathcal{O}_{v}\otimes_{\mathbf{Z}_{p}}R)$
for $R$ any $\mathbf{Z}_{p}$-algebra and likewise for $B$ and $T$.
Let $I$ be the Iwahori subgroup of $G(\mathbf{Z}_{p})$ associated
with $B$.

Let $\sw$ be the rigid analytic space whose $\overline{\mathbf{Q}_{p}}$-points
are given by\[
\scw(\overline{\mathbf{Q}_{p}})=\mathrm{Hom}_{\mathrm{cts}}(T(\mathbf{Z}_{p}),\overline{\mathbf{Q}_{p}}^{\times}).\]
We denote by $\lambda$ both a $\overline{\mathbf{Q}_{p}}$-point
of $\scw$ and the corresponding character of $T(\mathbf{Z}_{p})$.
Given any such character $\lambda:T(\mathbf{Z}_{p})\to\overline{\mathbf{Q}_{p}}^{\times}$,
the image of $\lambda$ generates a subfield $k_{\lambda}\subset\overline{\mathbf{Q}_{p}}$
finite over $\mathbf{Q}_{p}$. Following ideas of Stevens and Ash-Stevens \cite{Strigidsymbs,AS},
we define a Fréchet $k_{\lambda}$-module $\scd_{\lambda}$ of locally
analytic distributions equipped with a continuous $k_{\lambda}$-linear
left action of $I$. For any open compact subgroup $K^{p}\subset\mathbf{G}(\mathbf{A}_{f}^{p})$,
the quotient\[
\left(\mathbf{G}(\mathbf{Q})\backslash\mathbf{G}(\mathbf{A})/K_{\infty}^{\circ}K^{p}\times\scd_{\lambda}\right)/I\to Y(K^{p}I)\]
defines a local system on $Y(K^{p}I)$ which we also denote by $\scd_{\lambda}$.
This local system is nontrivial if and only if $\lambda$ is trivial
on $Z_{\mathbf{G}}(\mathbf{Q})\cap K^{p}I\subset T(\mathbf{Z}_{p})$,
and this condition cuts out a closed equidimensional rigid subspace
$\scw_{K^{p}}\subset\scw$. The cohomology $H^{\ast}(Y(K^{p}I),\scd_{\lambda})$
is naturally a Hecke module, and when $\lambda$ is a $B$-dominant
algebraic weight with associated algebraic representation $\scl_{\lambda}$,
there is a surjective $I$-equivariant map $i_{\lambda}:\scd_{\lambda}\twoheadrightarrow\scl_{\lambda}$
which induces a Hecke-equivariant and degree-preserving map\[
H^{\ast}\left(Y(K^{p}I),\scd_{\lambda}\right)\to H^{\ast}\left(Y(K^{p}I),\scl_{\lambda}\right).\]
The target, by Matsushima's formula and its generalizations, is isomorphic
as a Hecke module to a finite-dimensional space of classical automorphic
forms; the source, on the other hand, is much larger. When $F=\mathbf{Q}$
and $\mathbf{G}=\mathrm{GL}_{2}$, Stevens proved that the finite-slope
systems of Hecke eigenvalues appearing in $H^{\ast}(Y(K^{p}I),\scd_{\lambda})$
are exactly those appearing in a corresponding space of overconvergent
modular forms. These cohomology modules make sense, however, even
for groups without associated Shimura varieties, and in our opinion
they are the {}``correct'' surrogate for spaces of overconvergent
\emph{p-}adic modular forms.

To explain our results, we need to be precise about the Hecke algebras
under consideration. After choosing a uniformizer $\varpi_{v}$ of
$\mathcal{O}_{v}$ for each $v|p$, the action of $I$ on $\scd_{\lambda}$
extends canonically to an action of the monoid $\Delta\subset G(\mathbf{Q}_{p})$
generated by $I$ and by the monoid \[
T^{+}=\left\{ t\in T(\mathbf{Q}_{p})\mid t^{-1}B(\mathbf{Z}_{p})t\subseteq B(\mathbf{Z}_{p})\right\} .\]
The algebra $\mathcal{A}_{p}^{+}=\scc_{c}^{\infty}(I\backslash\Delta/I,\mathbf{Q}_{p})$
is a commutative subalgebra of the Iwahori-Hecke algebra of $G$.
Letting $U_{t}=[ItI]\in\mathcal{A}_{p}^{+}$ be the Hecke operator
associated with any $t\in T^{+}$, the map $t\mapsto U_{t}$ extends
to an algebra isomorphism\[
\mathbf{Q}_{p}[T^{+}/T(\mathbf{Z}_{p})]\overset{\sim}{\to}\mathcal{A}_{p}^{+}.\]
We define $U_{t}$ to be a \emph{controlling operator }if $t\in T^{+}$
satisfies $\cap_{i\geq1}t^{-i}B(\mathbf{Z}_{p})t^{i}=\{1\}$.

Let $S(K^{p})$ denote the finite set of places of $F$ where either
$v|p$, or $\mathbf{H}/F_{v}$ is ramified, or $K_{v}^{p}=\mathbf{H}(F_{v})\cap K^{p}$
is not a hyperspecial maximal compact subgroup of $\mathbf{H}(F_{v})$.
The main Hecke algebra of interest for us is\[
\mathbf{T}(K^{p})=\mathcal{A}_{p}^{+}\otimes_{\mathbf{Q}_{p}}\bigotimes_{v\notin S(K^{p})}'\scc_{c}^{\infty}(K_{v}^{p}\backslash\mathbf{H}(F_{v})/K_{v}^{p},\mathbf{Q}_{p}).\]

\textbf{Definition 1.1.1. }\emph{A }finite-slope eigenpacket (of weight
$\lambda$ and level $K^{p}$)\emph{ is an algebra homomorphism $\phi:\mathbf{T}(K^{p})\to\overline{\mathbf{Q}_{p}}$
such that the space\[
\left\{ v\in H^{\ast}(Y(K^{p}I),\scd_{\lambda})\otimes_{k_{\lambda}}\overline{\mathbf{Q}_{p}}\mid T\cdot v=\phi(T)v\,\forall T\in\mathbf{T}(K^{p})\right\} \]
is nonzero and such that $\phi(U_{t})\neq0$ for some controlling
operator $U_{t}$.}

This definition is independent of the specific choice of controlling
operator, and has some consequences which are not entirely obvious:
in particular, the set of finite-slope eigenpackets of given weight
and level which satisfy $v_{p}(\phi(U_{t}))\leq h$ is \emph{finite
}for any fixed $h$. The image of any finite-slope eigenpacket $\phi$
generates a subfield of $\overline{\mathbf{Q}_{p}}$ finite over $\mathbf{Q}_{p}$,
and we denote by $k_{\phi}$ the compositum of this field with $k_{\lambda}$.

Our first main result is the existence of an eigenvariety parametrizing
the finite-slope eigenpackets appearing in $H^{\ast}(Y(K^{p}I),\scd_{\lambda})$
as $\lambda$ varies over $\scw_{K^{p}}$.

\textbf{Theorem 1.1.2. }\emph{Notation and assumptions as above, there
is a canonical separated rigid analytic space $\scx=\scx_{\mathbf{G},K^{p}}$
equipped with a morphism $w:\scx\to\scw_{K^{p}}$ and an algebra homomorphism
$\phi_{\scx}:\mathbf{T}(K^{p})\to\sco(\scx)$ such that:}
\begin{description}
\item [{i.}] \emph{The morphism $w$ has discrete fibers and is locally
finite in the domain.}
\item [{ii.}] \emph{For any point $\lambda\in\scw_{K^{p}}(\overline{\mathbf{Q}_{p}})$,
there is a canonical bijection between points in the fiber $w^{-1}(\lambda)\subset\scx(\overline{\mathbf{Q}_{p}})$
and finite-slope eigenpackets of weight $\lambda$ and level $K^{p}$,
realized by the map sending $x\in w^{-1}(\lambda)$ to the algebra
homomorphism\[
\phi_{\scx,x}:\mathbf{T}(K^{p})\overset{\phi_{\scx}}{\to}\sco(\scx)\to k_{x}.\]
Writing $\phi\mapsto x(\phi)$ for the inverse of this map, we have
$k_{\phi}=k_{x(\phi)}$.}
\item [{iii.}] \emph{There are canonically defined sheaves of automorphic
forms on $\sx$ interpolating the modules $H^{\ast}(Y(K^{p}I),\scd_{\lambda})$.}
\item [{iv.}] \emph{For any controlling operator $U_{t}$, there is a canonical
closed immersion $\sz_{t}\hookrightarrow\scw_{K^{p}}\times\mathbf{A}^{1}$
and finite morphism $s:\scx\to\scz_{t}$ such that $w$ factors as
the composite of $s$ with the morphism $\scz_{t}\to\scw_{K^{p}}\times\mathbf{A}^{1}\overset{\mathrm{pr}_{1}}{\to}\scw_{K^{p}}$.}
\end{description}
Note that part iv. of this theorem characterizes the Grothendieck
topology on $\scx_{\mathbf{G},K^{p}}$. We construct $\scx_{\mathbf{G},K^{p}}$
through a by-now-familiar process of gluing suitable affinoid local
pieces. These affinoids are already defined in \cite{AS}, but gluing
them turns out to be a somewhat subtle affair. The main novelty of
our analysis is not the gluing, however, but rather the construction
of some universal coefficients spectral sequences which allow us to
give a fairly soft analysis of the resulting spaces. We shall say
more about both these points below.

In order to put this theorem in context, and to partially explain
the title of this paper, we introduce a little formalism. Suppose
we are given, for each $\lambda\in\scw_{K^{p}}(\overline{\mathbf{Q}_{p}})$,
a $k_{\lambda}\otimes\mathbf{T}(K^{p})$-module $M_{\lambda}$ of
{}``overconvergent modular forms of weight $\lambda$ and tame level
$K^{p}$ on $\mathbf{G}$.'' If the assignment $\lambda\mapsto M_{\lambda}$
varies analytically with $\lambda$, then one has a chance of constructing
an eigenvariety as in Theorem 1.1.2 whose fiber over any $\lambda$
parametrizes the finite-slope systems of Hecke eigenvalues occurring
in $M_{\lambda}$. In practice there are several definitions of suitable
families $\{M_{\lambda}\}_{\lambda}$ for which this construction
can be carried out, which admit no \emph{a priori }comparison. We
shall take up the task of systematically comparing the resulting eigenvarieties
in \cite{Hcomp}. For now, let us make the following remarks:
\begin{itemize}
\item The union of irreducible components of $\sx$ where $\phi_{\scx}(U_{t})$
is a \emph{p-}adic unit admits a canonical integral model. The study
of these components with their integral structures is usually known
as {}``Hida theory'', following on the pioneering works of H. Hida \cite{Hidordinv,Hidtotreal,HidaSLn},
and is best carried out from a formal-schematic perspective, as opposed
to the rigid analytic setup of the present paper.
\item When $\mathbf{G}^{\mathrm{der}}(\mathbf{R})$ is compact, the spaces
$Y(K_{f})$ are finite sets of points, and $H^{\ast}(Y(K^{p}I),\scd_{\lambda})=H^{0}(Y(K^{p}I),\scd_{\lambda})$
is a space of {}``algebraic overconvergent modular forms.'' In this
case, a number of authors have worked out special cases of the construction
underlying Theorem 1.1.2 \cite{BCast,BuEigen,BuFam,ChFamilies,Loeffler,Taibi}.
\item If $\mathbf{G}$ gives rise to Shimura varieties of PEL type, one
might try to construct suitable $M_{\lambda}$'s using coherent cohomology.
In the case $\mathbf{G}=\mathrm{GL}_{2}/\mathbf{Q}$, Coleman and
Mazur initiated the entire theory of eigenvarieties in this way with
the construction of their famous eigencurve \cite{CMeigencurve},
building on the results in \cite{Colemanbanachfamilies,Colemanclassicaloc}.
Until recently, further coherent-cohomological constructions of $M_{\lambda}$
for general weights $\lambda$, e.g. \cite{KL,MokTan}, have relied
(as did Coleman and Mazur) on tricks involving Eisenstein series or
lifts of Hasse invariants, and have only given rise to one-dimensional
families. However, canonical constructions of overconvergent modular
forms of arbitrary \emph{p-}adic weight on Shimura varieties have
recently been discovered \cite{AIS,AIP,Br,CHJ} which do give rise
to universal coherent-cohomological eigenvarieties.
\item When $\mathbf{G}^{\mathrm{der}}(\mathbf{R})$ has a discrete series,
there is a natural closed immersion from Urban's eigenvariety into
$\scx_{\mathbf{G},K^{p}}$, and the image of this map is exactly the
union of the $\mathrm{dim}\scw_{K^{p}}$-dimensional irreducible components
of $\scx_{\mathbf{G},K^{p}}$.
\item For general $\mathbf{G}$, Emerton \cite{Eminterpolate} gave a construction
of eigenvarieties in which $M_{\lambda}$ is (essentially) the weight-$\lambda$
subspace of the locally analytic Jacquet module of completed cohomology.
We'll carefully compare Emerton's construction with the construction
of the present paper in \cite{Hcomp}, generalizing the comparison
given in \cite{Loeffler} in the case when $\mathbf{G}^{\mathrm{der}}(\mathbf{R})$
is compact.
\end{itemize}
One of the main results we shall prove in \cite{Hcomp} is that for
a given group $\mathbf{G}$ and tame level $K^{p}$, every relevant
eigenvariety on the above list (and in particular, every eigenvariety
known to the author) admits a natural closed immersion into $\scx_{\mathbf{G},K^{p}}$.
In this light, it seems reasonable to regard the spaces $\scx_{\mathbf{G},K^{p}}$
as universal eigenvarieties, hence the title.

Let us outline the proof of Theorem 1.1.2. Again following \cite{AS},
we define for any affinoid subdomain $\Omega\subset\scw_{K^{p}}$
a Fréchet $\sco(\Omega)$-module $\scd_{\Omega}$ with continuous
$I$-action, such that $\scd_{\Omega}\otimes_{\sco(\Omega)}k_{\lambda}\cong\scd_{\lambda}$
for all $\lambda\in\Omega(\overline{\mathbf{Q}_{p}})$. The cohomology
modules $H^{\ast}(Y(K^{p}I,\scd_{\Omega})$ don't \emph{a priori }carry
any natural topology, but after making some noncanonical choices,
we obtain a chain complex $C_{\Omega}^{\bullet}$ of Fréchet modules
such that $H^{\ast}(C_{\Omega}^{\bullet})\cong H^{\ast}(Y(K^{p}I),\scd_{\Omega})$
together with an action of an operator $\tilde{U}_{t}$ on $C_{\Omega}^{\bullet}$
lifting the action of any given controlling operator $U_{t}$ on cohomology.
Using a fundamental result of Buzzard and the coherence of the assignment
$\Omega\rightsquigarrow C_{\Omega}^{\bullet}$, we are then able to
glue the finite-slope part of $C_{\Omega}^{\bullet}$ over varying
$\Omega$ into a complex of coherent sheaves on a certain noncanonical
Fredholm hypersurface $\scz\subset\scw_{K^{p}}\times\mathbf{A}^{1}$.
The cohomology sheaves $\sm^{\ast}$ of this complex are completely
canonical, and Buzzard's result gives an admissible covering of $\scz$
by affinoids $\scz_{\Omega,h}$, for suitably varying $\Omega\subset\scw_{K^{p}}$
and $h\in\mathbf{Q}_{\geq0}$, such that $\scm^{\ast}(\scz_{\Omega,h})$
is canonically isomorphic to the {}``slope-$\leq h$ part'' $H^{\ast}(Y(K^{p}I),\scd_{\Omega})_{h}$
of $H^{\ast}(Y(K^{p}I),\scd_{\Omega})$. Recall that $H^{\ast}(Y(K^{p}I),\scd_{\Omega})_{h}$
- when it exists\emph{ - }is an $\sco(\Omega)$-module-finite Hecke-stable
direct summand of $H^{\ast}(Y(K^{p}I),\scd_{\Omega})$ characterized
as the maximal subspace of the latter on which (roughly) every eigenvalue
of $U_{t}$ has valuation $\leq h$. Using this identification, we
are able to glue the Hecke actions on $H^{\ast}(Y(K^{p}I),\scd_{\Omega})_{h}$
into an action on these cohomology sheaves, from which point we easily
obtain $\scx$ by a simple {}``relative Spec''-type construction.
We formalize this latter process in the definition of an \emph{eigenvariety
datum }(Definition 4.2.1).

Our main tool in analyzing this construction is family of spectral
sequences which allows us to recover $H^{\ast}(Y(K^{p}I),\scd_{\lambda})_{h}$
from $H^{\ast}(Y(K^{p}I),\scd_{\Omega})_{h}$. For example, we prove
the following result (our most general result in this direction is
Theorem 3.3.1).

\textbf{Theorem 1.1.3. }\emph{There is a Hecke-equivariant second-quadrant
spectral sequence\[
E_{2}^{i,j}=\mathrm{Tor}_{-i}^{\sco(\Omega)}\left(H^{j}(Y(K^{p}I),\scd_{\Omega})_{h},k_{\lambda}\right)\Rightarrow H^{i+j}(Y(K^{p}I),\scd_{\lambda})_{h}\]
for any $\lambda\in\Omega(\overline{\mathbf{Q}_{p}})$.}

This sequence and its relatives turn out to be powerful tools for
studying the geometry of eigenvarieties. To explain our results in
this direction, recall that when $\mathbf{G}^{\mathrm{der}}(\mathbf{R})$
has a discrete series, standard limit multiplicity results yield an
abundance of classical automorphic forms of essentially every arithmetic
weight. One expects that correspondingly e$ $very irreducible component
$\mathscr{X}_{i}$ of the eigenvariety $\mathscr{X}$ which contains
a {}``suitably general'' classical point has maximal dimension,
namely $\mathrm{dim}\mathscr{X}_{i}=\mathrm{dim}\scw_{K^{p}}$. This
numerical coincidence is\emph{ characteristic }of the groups for which
$\mathbf{G}^{\mathrm{der}}(\mathbf{R})$ has a discrete series. More
precisely, define the \emph{defect }and the \emph{amplitude }of $\mathbf{G}$,
respectively, as the integers $l(\mathbf{G})=\mathrm{rank}\mathbf{G}-\mathrm{rank}K_{\infty}$
and $q(\mathbf{G})=\frac{1}{2}(\mathrm{dim}(\mathbf{G}(\mathbf{R})/K_{\infty})-l(\mathbf{G}))$.%
\footnote{Here {}``rank'' denotes the absolute rank, i.e. the dimension of
any maximal torus, split or otherwise.%
} Note that $l(\mathbf{G})$ is zero if and only if $\mathbf{G}^{\mathrm{der}}(\mathbf{R})$
has a discrete series, and that algebraic representations with regular
highest weight contribute to $(\mathfrak{g},K_{\infty})$-cohomology
exactly in the unbroken range of degrees $[q(\mathbf{G}),q(\mathbf{G})+l(\mathbf{G})]$.
We say a point $x\in\mathscr{X}(\overline{\mathbf{Q}_{p}})$ is \emph{classical
}if the weight $\lambda_{x}$ factors as $\lambda_{x}=\lambda^{\mathrm{alg}}\varepsilon$
with $\varepsilon$ finite-order and $\lambda^{\mathrm{alg}}$ $B$-dominant
algebraic and the associated eigenpacket $\phi_{x}$ matches the Hecke
data of an algebraic automorphic representation $\pi$ of $\mathbf{G}(\mathbf{A_{Q}})$
such that $\pi_{\infty}$ contributes to $(\mathfrak{g},K_{\infty})$-cohomology
with coefficients in an irreducible algebraic representation of highest
weight $\lambda^{\mathrm{alg}}$. A classical point is \emph{regular
}if $\lambda^{\mathrm{alg}}$ is regular. The definitions of\emph{
non-critical}, \emph{interior}, and \emph{strongly interior} points
are slightly more subtle and we defer them until §3.2. The following
conjecture is a special case of a conjecture of Urban (Conjecture
5.7.3 of \cite{UrEigen}).

\textbf{Conjecture 1.1.4. }\emph{Every irreducible component $\mathscr{X}_{i}$
of $\mathscr{X}_{\mathbf{G},K^{p}}$ containing a strongly interior,
non-critical, regular classical point has dimension $\mathrm{dim}\scw_{K^{p}}-l(\mathbf{G})$.}

Using the spectral sequences, we verify Conjecture 1.1.4 in many cases.
Our techniques yield new results even for groups with discrete series.

\textbf{Theorem 1.1.5.}\emph{ If $l(\mathbf{G})\leq1$, then Conjecture
1.1.4 is true, and if $l(\mathbf{G})\geq1$, every irreducible component
of $\mathscr{X}_{\mathbf{G},K^{p}}$ containing a strongly interior,
non-critical, regular classical point has dimension at most $\mathrm{dim}\scw_{K^{p}}-1$.}

In fact we prove slightly more, cf. Theorem 4.5.1. Some basic examples
of groups with $l(\mathbf{G})=1$ include $\mathrm{GL}_{3}/\mathbf{Q}$
and $\mathrm{Res}_{F/\mathbf{Q}}\mathbf{H}$ where $F$ is a number
field with exactly one complex embedding and $\mathbf{H}$ is an $F$-inner
form (possibly split) of $\mathrm{GL}_{2}$.

It turns out that our techniques imply much more. Given any $x\in\scx_{\mathbf{G},K^{p}}(\overline{\mathbf{Q}_{p}})$,
set\[
l(x)=\mathrm{sup}\left\{ i\mid H^{i}(Y(K^{p}I),\scd_{\lambda_{x}})_{(\ker\phi_{x})}\neq0\right\} -\mathrm{inf}\left\{ i\mid H^{i}(Y(K^{p}I),\scd_{\lambda_{x}})_{(\ker\phi_{x})}\neq0\right\} .\]
We shall see that interior, noncritical, regular classical points
satisfy $l(x)=l(\mathbf{G})$. Upon reading an earlier version of this
paper, James Newton discovered a proof of the following result, which
is given here in Appendix B (cf. also the remark at the end of §4.5).

\textbf{Theorem 1.1.6 (Newton). }\emph{The dimension of any irreducible
component $\scx_{i}$ of $\scx_{\mathbf{G},K^{p}}$ containing a given
point $x$ satisfies $\mathrm{dim}\scx_{i}\geq\mathrm{dim}\scw_{K^{p}}-l(x)$.
In particular, the lower bound of Conjecture 1.1.4 is true.}

\subsection{The conjectural connections with Galois representations}

$ $In this section we restrict our attention to the case $\mathbf{G}=\mathrm{Res}_{F/\mathbf{Q}}\mathrm{GL}_{n}$,
choosing $B$ upper triangular and $T$ diagonal, so $T(\mathbf{Z}_{p})\cong\prod_{v|p}(\mathcal{O}_{F_{v}}^{\times})^{n}$
and $T(\mathbf{Q}_{p})\cong\prod_{v|p}(F_{v}^{\times})^{n}$ with
the obvious diagonal coordinates. For a given tame level $K^{p}$
we abbreviate $\scx_{K^{p}}=\scx_{\mathbf{G},K^{p}}$, and we write
$T_{v,i}$ and $U_{v,i}$ for the usual Hecke operators at places
$v\notin S(K^{p})$ and $v|p$, respectively (cf. §4.6).

Let $G_{F}$ be the absolute Galois group of $F$, and let $\rho:G_{F}\to\mathrm{GL}_{n}(\overline{\mathbf{Q}_{p}})$
be a continuous semisimple representation which is unramified almost
everywhere. We say a tame level $K^{p}$ is admissible for $\rho$
if the set $S(K^{p})$ contains the set of places where $\rho$ is
ramified. 

\textbf{Definition 1.2.1. }\emph{If $K^{p}$ is admissible for $\rho$
and $\phi:\mathbf{T}(K^{p})\to\overline{\mathbf{Q}_{p}}$ is an algebra
homomorphism, $\rho$ and $\phi$ are }\textbf{associated }\emph{if
the equality\[
\det(I_{n}-X\cdot\rho(\mathrm{Frob}_{v}))=\sum_{i=0}^{n}(-1)^{i}\mathbf{N}v^{\frac{i(i-1)}{2}}\phi(T_{v,i})X^{i}\,\mathrm{in}\,\overline{\mathbf{Q}_{p}}[X]\]
holds for all $v\notin S(K^{p})$.}

This is a standard incarnation of the usual reciprocity between Galois
representations and automorphic forms. Note that for any given $\phi$,
there is at most one isomorphism class of continuous semisimple Galois
representations associated with $\phi$, by the Brauer-Nesbitt theorem.
Note also that if $\phi$ is associated with a Galois representation,
the image of $\phi$ is necessarily contained in a finite extension
of $\mathbf{Q}_{p}$.

\textbf{Conjecture 1.2.2. }\emph{Given any point $x\in\sx_{K^{p}}(\overline{\mathbf{Q}_{p}})$
with weight $\lambda_{x}$ and corresponding eigenpacket $\phi_{x}$,
there is a continuous n-dimensional semisimple representation $\rho_{x}:G_{F}\to\mathrm{GL}_{n}(\overline{\mathbf{Q}_{p}})$
with the following properties:}
\begin{description}
\item [{i.}] \emph{The tame level $K^{p}$ is admissible for $\rho_{x}$,
and $\rho_{x}$ and $\phi_{x}$ are associated.}
\item [{ii.}] \emph{The representation $\rho_{x}$ is odd: for any real
infinite place $v$ with complex conjugation $c_{v}$, we have\[
\mathrm{tr}\rho_{x}(c_{v})=\begin{cases}
\pm1 & \mathrm{if}\, n\,\mathrm{is\, odd}\\
0 & \mathrm{if}\, n\,\mathrm{is\, even.}\end{cases}\]
}
\item [{iii.}] \emph{For each place $v|p$, $\rho_{x}|G_{F_{v}}$ is trianguline,
and the space of crystalline periods \[
\mathbf{D}_{\mathrm{crys}}^{+}\left((\wedge^{i}\rho_{x}|G_{F_{v}})\otimes(\lambda_{x,n}\cdots\lambda_{x,n+1-i}\cdot\mathbf{N}^{\frac{i(i-1)}{2}})\circ\chi_{v}\right)^{\varphi_{v}^{f_{v}}=\phi_{x}(U_{v,i})}\]
is nonzero for each $1\leq i\leq n$, where $\chi_{v}:G_{F_{v}}\to\mathcal{O}_{v}^{\times}$
is the Lubin-Tate character associated with our chosen uniformizer
$\varpi_{v}$.}
\end{description}
Part iii. of this conjecture is naturally inspired by a famous result
of Kisin \cite{KisinOCFM} and its generalizations due to Bellaïche-Chenevier,
Hellmann and others \cite{BCast,Hellmann}. When $F$ is totally real
or CM, the existence of $\rho_{x}$ satisfying parts i. and ii. of
this conjecture can be deduced from the recent work of Scholze, but
part iii. seems quite difficult even for $\mathrm{GL}_{3}/\mathbf{Q}$. 

We would like to formulate a converse to this conjecture. More precisely,
suppose we are given a continuous, absolutely irreducible, almost
everywhere unramified representation $\rho:G_{F}\to\mathrm{GL}_{n}(\overline{\mathbf{Q}_{p}})$
which is odd at all real infinite places. Let $N(\rho)\subset\mathcal{O}_{F}$
be the prime-to-\emph{p }Artin conductor of $\rho$.%
\footnote{When $F=\mathbf{Q}$ we conflate $N(\rho)$ with its positive generator
in the obvious way.%
} We define \[
\scx[\rho]=\left\{ x\in\scx_{K_{1}(N(\rho))}(\overline{\mathbf{Q}_{p}})\mid\rho\,\mathrm{and}\,\phi_{x}\,\mathrm{are}\,\mathrm{associated}\right\} ,\]
where $K_{1}(N)$ is the usual level subgroup appearing in the theory
of new vectors for $\mathrm{GL}_{n}$.%
\footnote{We sometimes write {}``$\rho\simeq\rho_{x}$'' for {}``$\rho$
and $\phi_{x}$ are associated'', since the latter certainly implies
that $\rho_{x}$ exists!%
} What can we say about this set of points?

\textbf{Conjecture 1.2.3. }\emph{The set $\scx[\rho]$ is nonempty
if and only if $\rho$ is trianguline at all places dividing $p$.}

To formulate a more quantitative statement, let $\st=\st_{n,F}$ denote
the rigid space with $\st(L)=\mathrm{Hom}_{\mathrm{cts}}(T(\mathbf{Q}_{p}),L)$.
We denote by $\delta$ both a point of $\st(\overline{\mathbf{Q}_{p}})$
and the associated character $\delta:T(\mathbf{Q}_{p})\to\overline{\mathbf{Q}_{p}}^{\times}$,
and we identify any such character $\delta$ with an ordered $n$-tuple
of continuous characters $\delta_{i}:\prod_{v|p}F_{v}^{\times}\to\overline{\mathbf{Q}_{p}}^{\times}$
in the natural way. Given any point $x\in\scx_{K^{p}}(\overline{\mathbf{Q}_{p}})$,
with associated weight $\lambda_{x}$, define $\delta_{x}\in\st(\overline{\mathbf{Q}_{p}})$
as follows:

i. $\delta_{x}(\underbrace{\varpi_{v},\dots,\varpi_{v}}_{i},1,\dots,1)=\phi_{\scx,x}(U_{v,i})$,
and

ii. $\delta_{x}(t)=\prod_{i=1}^{n}\lambda_{x,n+1-i}(t_{i}^{-1})\mathbf{N}t_{i}^{1-i}$
for $ $$t=\mathrm{diag}(t_{1},\dots,t_{n})\in T(\mathbf{Z}_{p})$.\\
Clearly there is a unique global character $\delta_{\scx}:T(\mathbf{Q}_{p})\to\sco(\scx_{K^{p}})^{\times}$
specializing to $\delta_{x}$ at every point $x$. Let $\st[\rho]\subset\st(\overline{\mathbf{Q}_{p}})$
be the image of the map\begin{eqnarray*}
\scx[\rho] & \to & \st\\
x & \mapsto & \delta_{x}.\end{eqnarray*}
The map $\scx[\rho]\to\st[\rho]$ is a bijection, since the data of
the pair $(\rho,\delta_{x})$ is enough to recover the point $x$
uniquely by Theorem 1.1.2. Describing the set $\st[\rho]$ amounts
roughly to a characteristic zero analogue of the {}``weight part''
of Serre's modularity conjecture and its generalizations, and just
as with Serre's conjecture, some elements of $\st[\rho]$ are easier
to predict than others.

Define $\ppar(\rho)$, the set of \emph{parameters of $\rho$}, as
the set of characters $\delta\in\st(\overline{\mathbf{Q}_{p}})$ such
that $\delta|_{T_{v}(\mathbf{Q}_{p})}$ is the parameter of some triangulation
of $\mathbf{D}_{\mathrm{rig}}^{\dagger}(\rho|G_{F_{v}})$ for all
$v|p$. This set is clearly nonempty if and only if $\rho|G_{F_{v}}$
is trianguline for each $v|p$, and depends only on the possible triangulations
of these local representations (we refer the reader to §6.1 for some
background on trianguline representations).

\textbf{Conjecture 1.2.4. }\emph{If $\delta$ is a parameter of $\rho$,
there is a unique point $x=x(\rho,\delta)\in\scx_{K_{1}(N(\rho))}(\overline{\mathbf{Q}_{p}})$
such that $\rho\simeq\rho_{x}$ and $\delta=\delta_{x}$. Equivalently,
$\ppar(\rho)\subset\st[\rho]$.}

Our next conjecture gives a complete description of the set $\st[\rho]$.
To keep this introduction at a reasonable length, we refer the reader
to §6.2 for two key definitions: briefly, given any character $\delta\in\st_{n,F}(\overline{\mathbf{Q}_{p}})$
we define a group $W(\delta)<S_{n}^{\mathrm{Hom}(F,\overline{\mathbf{Q}_{p}})}$
and a finite set of characters $W(\delta)\cdot\delta\subset\st_{n,F}(\overline{\mathbf{Q}_{p}})$,
together with a certain partial ordering {}``$\preceq$'' on $W(\delta)\cdot\delta$.
(The set $W(\delta)\cdot\delta$ is the $W(\delta)$-orbit of $\delta$
under a certain action, and in particular contains $\delta$.)

\textbf{Conjecture 1.2.5. }\emph{The set $\st[\rho]$ consists exactly
of those characters $\eta$ such that for some $\delta\in\ppar(\rho)$
we have $\eta\in W(\delta)\cdot\delta$ and $\delta\preceq\eta$.}

This gives a purely local and Galois-theoretic description of the
automorphically defined set $\scx[\rho]$. When $\rho$ is crystabelline
this conjecture (or rather, its natural analogue for unitary groups)
recovers a recent conjecture of Breuil (Conjecture 6.5 of \cite{Breuil}).
We refer the reader to §6.2 for some further discussion of Conjecture
1.2.5 and its local analogues.

By way of evidence, we have the following results. For the remainder
of this section we assume $F=\mathbf{Q}$, in which case an element
$\delta\in\ppar(\rho)$ is simply the parameter of some triangulation
of $\mathbf{D}_{\mathrm{rig}}^{\dagger}(\rho|G_{\mathbf{Q}_{p}})$.

\textbf{Theorem 1.2.6. }\emph{Notation and assumptions as above, Conjecture
1.2.5 is true when $n=2$, $\overline{\rho}:G_{\mathbf{Q}}\to\mathrm{GL}_{2}(\overline{\mathbf{F}_{p}})$
is absolutely irreducible, and $\overline{\rho}|G_{\mathbf{Q}_{p}}$
is not isomorphic to a twist of $\left(\begin{array}{cc}
1 & \ast\\
 & \overline{\chi}_{\mathrm{cyc}}\end{array}\right)$.}

This result is due almost entirely to others, and the proof is simply
a matter of assembling their results. More precisely, under the hypotheses
of Theorem 1.2.6, Emerton \cite{EmertonLG} proved that $\rho$ is
a twist of the Galois representation $\rho_{f}$ associated with a
finite-slope overconvergent cuspidal eigenform $f$. The result then
follows from work of Stevens and Bellaïche \cite{Stfamilies,BeCrit}
showing that the Hecke data associated with overconvergent eigenforms
appears in overconvergent cohomology. The situation is most interesting
when the weight of $f$ is an integer $k\geq2$ and $U_{p}f=\alpha f$
with $v_{p}(\alpha)>k-1$, in which case Conjecture 1.2.4 actually
predicts the existence of the \emph{companion form }of\emph{ $f$}:
a form $g$ of weight $2-k$ such that $\rho_{f}\simeq\rho_{g}\otimes\chi_{\mathrm{cyc}}^{1-k}$ \cite{Colemanclassicaloc}. We should remind the reader that Emerton's
magisterial work relies on the full force of the \emph{p-}adic local
Langlands correspondence for $\mathrm{GL}_{2}/\mathbf{Q}_{p}$ \cite{BergerBreuil,ColPGlong,KisDef},
not to mention Khare and Wintenberger's proof of Serre's conjecture \cite{KW1,KW2}.
We should also note that in general we need a certain level-lowering
result to deduce that $f$ really does occur at the minimal level
$N(\rho_{f})$ (cf. §6.3).

Surprisingly, we are able to offer some evidence for these conjectures
beyond the cases where $n=2$ or $\rho$ is geometric. To state our
results in this direction, let $\rho_{f}$ be the Galois representation
associated with a finite-slope cuspidal overconvergent eigenform $f$,
and let $\delta_{f}=(\delta_{f,1},\delta_{f,2})$ be the canonical
parameter of $\rho_{f}$ (cf. §6.1). Suppose the residual representation
$\overline{\rho}_{f}$ satisfies the hypotheses of Theorem 1.2.6.

\textbf{Theorem 1.2.7. }\emph{Conjecture 1.2.4 is true for the representation
$\mathrm{sym}^{2}\rho_{f}$ and the parameter $ $\[
\mathrm{sym}^{2}\delta_{f}=(\delta_{f,1}^{2},\delta_{f,1}\delta_{f,2},\delta_{f,2}^{2})\in\ppar(\mathrm{sym}^{2}\rho_{f}).\]
}Now choose a second eigenform $g$ with Galois representation $\rho_{g}$
and canonical parameter $\delta_{g}$. Define a map\begin{eqnarray*}
\st_{2,\mathbf{Q}}\times\st_{2,\mathbf{Q}} & \to & \st_{4,\mathbf{Q}}\\
(\delta,\delta') & \mapsto & \delta\boxtimes\delta'\end{eqnarray*}
by $\delta\boxtimes\delta'=(\delta_{1}\delta_{1}',\delta_{1}\delta_{2}',\delta_{2}\delta_{1}',\delta_{2}\delta_{2}')$.

\textbf{Theorem 1.2.8. }\emph{Supposing $f$ and $g$ have tame level
one, Conjecture 1.2.4 is true for the representation $\rho=\rho_{f}\otimes\rho_{g}$
and the parameters $\delta_{f}\boxtimes\delta_{g}$ and $\delta_{g}\boxtimes\delta_{f}$.}

The assumption on the level is only for simplicity and can easily
be removed. Note that even when $f$ and $g$ are refinements of classical
level one eigenforms of distinct weights (and therefore $\rho$ is
crystalline with distinct Hodge-Tate weights at $p$), at least one
of these two parameters is \emph{critical }in the sense of Bellaïche-Chenevier,
and the weight of the associated point is algebraic but not $B$-dominant!
In this case, Theorem 1.2.8 asserts the existence of certain overconvergent
companion forms on the split form of $\mathrm{GL}_{4}/\mathbf{Q}$.

\subsection{\emph{p}-adic Langlands functoriality}

Given an $\mathbf{R}$-nonsplit quaternion algebra $D/\mathbf{Q}$
of discriminant $d$, Buzzard \cite{BuFam,BuEigen} constructed a
certain eigencurve $\scc_{D}$ using overconvergent algebraic modular
forms on $D$, and raised the question of whether there exists a closed
immersion $\iota_{\mathrm{JL}}:\scc_{D}\hookrightarrow\scc(d)$ into
the tame level $d$ eigencurve interpolating the Jacquet-Langlands
correspondence on classical points. Chenevier affirmatively answered
this question in a beautiful paper \cite{CHjacquetlanglands} as a
consequence of an abstract interpolation theorem.

In §5.1, inspired by Chenevier's results, we establish rather flexible
interpolation theorems (cf. Theorems 5.1.2 and 5.1.6). As sample applications
of these tools, we prove the following results.

\textbf{Theorem 1.3.1. }\emph{Let $N$ be a squarefree integer and
set $\scx=\scx_{\mathrm{GL}_{3}/\mathbf{Q},K_{1}(N^{2})}$. Let $\scc_{0}(N)$
denote the cuspidal locus in the tame level $N$ eigencurve $\scc(N)$.
Then there is a finite morphism $\mathbf{s}:\scc_{0}(N)\to\scx$ such
that $\rho_{\mathbf{s}(x)}\simeq\mathrm{sym}^{2}\rho_{x}$ and $\delta_{\mathbf{s}(x)}=\mathrm{sym}^{2}\delta_{x}$
for all $x\in\scc_{0}(N)(\overline{\mathbf{Q}_{p}})$.}

In fact we prove a much more precise result for arbitrary levels,
taking into account the inertial behavior of $\rho_{x}$ at all primes
$\ell\nmid p$; this immediately implies Theorem 1.2.7.

\textbf{Theorem 1.3.2. }\emph{Set $\scx=\scx_{\mathrm{GL}_{4}/\mathbf{Q},K_{1}(1)}$.
Then there is a finite morphism $\mathbf{t}:\scc_{0}(1)\times\scc_{0}(1)\to\scx$
such that $\rho_{\mathbf{t}(x,y)}\simeq\rho_{x}\otimes\rho_{y}$ and
$\delta_{\mathbf{t}(x,y)}=\delta_{x}\boxtimes\delta_{y}$ for all
$(x,y)\in\scc_{0}(1)(\overline{\mathbf{Q}_{p}})\times\scc_{0}(1)(\overline{\mathbf{Q}_{p}})$.}

This quickly implies Theorem 1.2.8. It's tempting to proliferate \emph{p-}adic
functorialities (and, simultaneously, evidence towards Conjectures
1.2.4 and 1.2.5) by combining Theorem 5.1.6 with known classical functorialities.%
\footnote{Although in general there will be subtle issues involving $L$-packets;
see \cite{Ludwig}.%
} For example, we invite the reader to construct a symmetric eighth
power map\[
\scx_{\mathrm{GL}_{2}/\mathbf{Q},K^{p}}\to\scx_{\mathrm{GL}_{9}/\mathbf{Q},K'^{p}}\]
(for compatible tame levels $K^{p},K'^{p}$) by applying Clozel and
Thorne's recent work \cite{CT}.

\subsection{Notation and terminology}

Our notation and terminology is mostly standard. For $p$ the prime
with respect to which things are -adic, we fix once and for all an
algebraic closure $\overline{\mathbf{Q}_{p}}$ and an isomorphism
$\iota:\overline{\mathbf{Q}_{p}}\overset{\sim}{\to}\mathbf{C}$. We
denote by $F$ (resp. $K$) a finite extension of $\mathbf{Q}$ (resp.
$\mathbf{Q}_{p}$). Unless otherwise noted, $L$ denotes a sufficiently
large subfield of $\overline{\mathbf{Q}_{p}}$ finite over $\mathbf{Q}_{p}$,
where the meaning of {}``sufficiently large'' may change from one
line to the next. If $F$ is a number field and $\rho:G_{F}\to\mathrm{GL}_{n}(L)$
is a Galois representation, and {}``blah'' is an adjective from
\emph{p-}adic Hodge theory (crystalline, semistable, de Rham, Hodge-Tate,
trianguline, etc.), we say {}``$\rho$ is blah'' as shorthand for
{}``$\rho|G_{F_{v}}$ is blah for all places $v|p$''.

We normalize the reciprocity maps of local class field theory so uniformizers
map to geometric Frobenii. If $\pi$ is an irreducible admissible
representation of $\mathrm{GL}_{n}(K)$, we write $\mathrm{rec}(\pi)$
for the Frobenius-semisimple Weil-Deligne representation associated
with $\pi$ via the local Langlands correspondence, normalized as
in Harris and Taylor's book. If $f=\sum_{n=1}^{\infty}a_{f}(n)q^{n}\in S_{k}(\Gamma_{1}(N))$
is a classical newform, we write $\rho_{f,\iota}$ (or just $\rho_{f}$)
for the two-dimensional semisimple $\overline{\mathbf{Q}_{p}}$-linear
representation of $G_{\mathbf{Q}}$ characterized by the equality
$\iota\mathrm{tr}\rho_{f}(\mathrm{Frob}_{\ell})=a_{f}(\ell)$ for
all $\ell\nmid Np$.

In nonarchimedean functional analysis and rigid analytic geometry
we follow \cite{BGR}. If $M$ and $N$ are topological $\mathbf{Q}_{p}$-vector
spaces, we write $\mathcal{L}(M,N)$ for the space of continuous $\mathbf{Q}_{p}$-linear
maps between $M$ and $N$; if $M$ and $N$ are $\mathbf{Q}_{p}$-Banach
spaces, the operator norm\[
|f|=\sup_{m\in M,\,|m|_{M}\leq1}|f(m)|_{N}\]
makes $\mathcal{L}(M,N)$ into a Banach space. If $(A,|\bullet|_{A})$
is a Banach space which furthermore is a commutative $\mathbf{Q}_{p}$-algebra
whose multiplication map is (jointly) continuous, we say $A$ is a
\textbf{$\mathbf{Q}_{p}$}-\emph{Banach algebra}. An $A$-module $M$
which is also a Banach space is a \emph{Banach $A$-module }if the
structure map $A\times M\to M$ extends to a continuous map $A\widehat{\otimes}_{\mathbf{Q}_{p}}M\to M$,
or equivalently if the norm on $M$ satisfies $\left|am\right|_{M}\leq C|a|_{A}|m|_{M}$
for all $a\in A$ and $m\in M$ with some fixed constant $C$. For
a topological ring $R$ and topological $R$-modules $M,N$, we write
$\mathcal{L}_{R}(M,N)$ for the $R$-module of continuous $R$-linear
maps $f:M\to N$. When $A$ is a Banach algebra and $M,N$ are Banach
$A$-modules, we topologize $\mathcal{L}_{A}(M,N)$ via its natural
Banach $A$-module structure. We write $\mathrm{Ban}_{A}$ for the
category whose objects are Banach $A$-modules and whose morphisms
are elements of $\mathcal{L}_{A}(-,-)$. If $I$ is any set and $A$
is a Banach algebra, we write $c_{I}(A)$ for the module of sequences
$\mathbf{a}=(a_{i})_{i\in I}$ with $|a_{i}|_{A}\to0$; the norm $|\mathbf{a}|=\sup_{i\in I}|a_{i}|_{A}$
gives $c_{I}(A)$ the structure of a Banach $A$-module. If $M$ is
any Banach $A$-module, we say $M$ is \emph{orthonormalizable }if
$M$ is \emph{isomorphic} to $c_{I}(A)$ for some $I$ (such modules
are called {}``potentially orthonormalizable'' in \cite{BuEigen}).

If $A$ is an affinoid algebra, then $\mathrm{Sp}A$, the \emph{affinoid
space }associated with $A$, denotes the locally G-ringed space $(\mathrm{Max}A,\mathcal{O}_{A})$
where $\mathrm{Max}A$ is the set of maximal ideals of $A$ endowed
with the Tate topology and $\mathcal{O}_{A}$ is the extension of
the assignment $U\mapsto A_{U}$, for affinoid subdomains $U\subset\mathrm{Max}A$
with representing algebras $A_{U}$, to a structure sheaf on $\mathrm{Max}A$.
If $X$ is an affinoid space, we write $\sco(X)$ for the coordinate
ring of $X$, so $A\simeq\sco(\mathrm{Sp}A)$. If $A$ is reduced
we equip $A$ with the canonical supremum norm. If $X$ is a rigid
analytic space, we write $\sco_{X}$ for the structure sheaf and $\sco(X)$
for the ring of global sections of $\sco_{X}$. Given a point $x\in X$,
we write $\mathfrak{m}_{x}$ for the corresponding maximal ideal in
$\sco_{X}(U)$ for any admissible affinoid open $U\subset X$ containing
$x$, and $k(x)$ for the residue field $\sco_{X}(U)/\mathfrak{m}_{x}$.
Given a point $x\in X(\overline{\mathbf{Q}_{p}})$, we write $k_{x}\subset\overline{\mathbf{Q}_{p}}$
for the image of $k(x)$ under the associated embedding $k(x)\hookrightarrow\overline{\mathbf{Q}_{p}}$.
$\sco_{X,x}$ denotes the local ring of $\mathcal{O}_{X}$ at $x$
in the Tate topology, and $\widehat{\mathcal{O}_{X,x}}$ denotes the
$\mathfrak{m}_{x}$-adic completion of $\mathcal{O}_{X,x}$. A Zariski-dense
subset $S$ of a rigid analytic space $X$ is a \emph{very Zariski-dense
subset}, or a \emph{Zariski-dense accumulation subset},\emph{ }if
for any connected affinoid open $U\subset X$ either $U\cap S=\emptyset$
or $U\cap S$ is Zariski-dense in $U$.

In homological algebra our conventions follow \cite{Weibelhomalg}.
If $R$ is a ring, we write $\mathbf{K}^{?}(R)$, $?\in\{+,-,b,\emptyset\}$
for the homotopy category of $?$-bounded $R$-module complexes and
$\mathbf{D}^{?}(R)$ for its derived category.

\subsection*{Acknowledgments}

This article is a revised and expanded version of my 2013 Boston College
Ph.D. thesis \cite{Hthesis}. First and foremost, I heartily thank
my advisor, Avner Ash, for his invaluable suggestions, sage advice,
patient readings of preliminary drafts, and overall kindness and generosity.
I'm grateful to be his student. I'm indebted to Glenn Stevens for
kindly encouraging me to work on overconvergent cohomology, for explaining
the utility of slope decompositions to me, and for generally serving
as a second advisor. The intellectual debt this article owes to the
ideas of Ash and Stevens will be obvious to the reader. I'm also grateful
to Ben Howard for some helpful conversations and for his detailed
remarks on a preliminary draft of my thesis.

During the development of this material, I enjoyed stimulating conversations
with Joël Bellaïche, John Bergdall, Kevin Buzzard, Przemys\l{}aw Chojecki,
Giovanni Di Matteo, Michael Harris, Eugen Hellmann, Keenan Kidwell,
Judith Ludwig, Barry Mazur, James Newton, Jay Pottharst, and Jack
Thorne, and it's a pleasure to acknowledge the help and influence
of all these mathematicians. I'm also especially grateful to Jack
for his detailed comments on several preliminary drafts over the years.
It's a particular pleasure to thank Barry Mazur for a number of inspiring
discussions, for his generous and invaluable encouragement, and for
providing financial support during the summer of 2013.  Finally, I'm grateful to the referee for helpful comments and corrections.

This work was carried out at Boston College and l\textbf{'}Institut
de Mathématiques de Jussieu, and I'm pleased to acknowledge the hospitality
of these institutions. The research leading to these results has received
funding from the European Research Council under the European Community's
Seventh Framework Programme (FP7/2007-2013) / ERC Grant agreement
n\textdegree{} 290766 (AAMOT).

\section{Background}

We maintain the notation of the introduction. Set $X^{\ast}=\mathrm{Hom}(T,\mathbf{G}_{m})$
and $X_{\ast}=\mathrm{Hom}(\mathbf{G}_{m},T)$, and let $\Phi$, $\Phi^{+}$
and $\Phi^{-}$ be the sets of roots, positive roots, and negative
roots respectively, for the Borel $B$. We write $X_{+}^{\ast}$ for
the cone of $B$-dominant weights; $\rho\in X^{\ast}\otimes_{\mathbf{Z}}\frac{1}{2}\mathbf{Z}$
denotes half the sum of the positive roots.

We write $\overline{B}$ for the opposite Borel, $N$ and $\overline{N}$
for the unipotent radicals of $B$ and $\overline{B}$, and $I$ for
the Iwahori subgroup\[
I=\left\{ g\in G(\mathbf{Z}_{p})\,\mathrm{with}\, g\,\mathrm{mod}\, p\in B(\mathbf{Z}/p\mathbf{Z})\right\} .\]
For any integer $s\geq1$, set $\overline{B}^{s}=\left\{ b\in\overline{B}(\mathbf{Z}_{p}),\, b\equiv1\,\mathrm{in}\, G(\mathbf{Z}/p^{s}\mathbf{Z})\right\} $,
$\overline{N}^{s}=\overline{N}(\mathbf{Z}_{p})\cap\overline{B}^{s}$
and $T^{s}=T(\mathbf{Z}_{p})\cap\overline{B}^{s}$, so the Iwahori
decomposition reads $I=\overline{N}^{1}\cdot T(\mathbf{Z}_{p})\cdot N(\mathbf{Z}_{p})$.
We also set \[
I_{0}^{s}=\left\{ g\in I,\, g\,\mathrm{mod}\, p^{s}\in\overline{B}(\mathbf{Z}/p^{s}\mathbf{Z})\right\} \]
and\[
I_{1}^{s}=\left\{ g\in I,\, g\,\mathrm{mod}\, p^{s}\in\overline{N}(\mathbf{Z}/p^{s}\mathbf{Z})\right\} .\]
Note that $I_{1}^{s}$ is normal in $I_{0}^{s}$, with quotient $T(\mathbf{Z}/p^{s}\mathbf{Z})$.
Finally, we set $I^{s}=I\cap\ker\left\{ G(\mathbf{Z}_{p})\to G(\mathbf{Z}/p^{s}\mathbf{Z})\right\} $.

We define semigroups in $T(\mathbf{Q}_{p})$ by\[
T^{+}=\left\{ t\in T(\mathbf{Q}_{p}),\, t\overline{N}^{1}t^{-1}\subseteq\overline{N}^{1}\right\} \]
and \[
T^{++}=\left\{ t\in T(\mathbf{Q}_{p}),\,\bigcap_{i=1}^{\infty}t^{i}\overline{N}^{1}t^{-i}=\{1\}\right\} .\]
A simple calculation shows that $t\in T(\mathbf{Q}_{p})$ is contained
in $T^{+}$ (resp. $T^{++}$) if and only if $v_{p}(\alpha(t))\leq0$
(resp. $v_{p}(\alpha(t))<0$) for all $\alpha\in\Phi^{+}$. Using
these semigroups, we define a semigroup of $G(\mathbf{Q}_{p})$ by
$\Delta=IT^{+}I$. The Iwahori decompsition extends to $\Delta$:
any element $g\in\Delta$ has a unique decomposition $g=\overline{\mathrm{n}}(g)\mathrm{t}(g)\mathrm{n}(g)$
with $\overline{\mathrm{n}}\in\overline{N}^{1}$, $\mathrm{t}\in T^{+}$,
$\mathrm{n}\in N(\mathbf{Z}_{p})$. Our chosen uniformizers of $\mathcal{O}_{v},v\mid p$,
induce a canonical group homomorphism $\sigma:T(\mathbf{Q}_{p})\to T(\mathbf{Z}_{p})$
which splits the inclusion $T(\mathbf{Z}_{p})\subset T(\mathbf{Q}_{p})$,
and we set $\Lambda=T^{+}\cap\ker\sigma$ and $\Lambda^{+}=T^{++}\cap\ker\sigma$.

\subsection{Symmetric spaces and Hecke operators}

In this section we set up our conventions for the homology and cohomology
of local systems on locally symmetric spaces. Following \cite{AS},
we compute homology and cohomology using two different families of
resolutions: some extremely large {}``adelic'' resolutions
which have the advantage of making the Hecke action transparent, and
resolutions with good finiteness properties constructed from simplicial
decompositions of the Borel-Serre compactifications of locally symmetric
spaces.

\subsubsection*{Resolutions and complexes}

Let $\mathbf{G}/\mathbf{Q}$ be a connected reductive group with center
$Z_{\mathbf{G}}$. Let $\mathbf{G}(\mathbf{R})^{\circ}$ denote the
connected component of $\mathbf{G}(\mathbf{R})$ containing the identity
element, with $\mathbf{G}(\mathbf{Q})^{\circ}=\mathbf{G}(\mathbf{Q})\cap\mathbf{G}(\mathbf{R})^{\circ}$.
Fix a maximal compact-mod-center subgroup $K_{\infty}\subset\mathbf{G}(\mathbf{R})$
with $K_{\infty}^{\circ}$ the connected component containing the
identity. Given an open compact subgroup $K_{f}\subset\mathbf{G}(\mathbf{A}_{f})$,
we define the \emph{locally symmetric space of level $K_{f}$ }by\begin{eqnarray*}
Y(K_{f}) & = & \mathbf{G}(\mathbf{Q})\backslash\mathbf{G}(\mathbf{A})/K_{f}K_{\infty}^{\circ}.\end{eqnarray*}
This is a possibly disconnected Riemannian orbifold. By strong approximation
there is a finite set of elements $\gamma(K_{f})=\{x_{i},x_{i}\in\mathbf{G}(\mathbf{A}_{f})\}$
with\[
\mathbf{G}(\mathbf{A})=\coprod_{x_{i}\in\gamma(K_{f})}\mathbf{G}(\mathbf{Q})^{\circ}\mathbf{G}(\mathbf{R})^{\circ}x_{i}K_{f}.\]
Defining $Z(K_{f})=Z_{\mathbf{G}}(\mathbf{Q})\cap K_{f}$ and $\tilde{\Gamma}(x_{i})=\mathbf{G}(\mathbf{Q})^{\circ}\cap x_{i}K_{f}x_{i}^{-1}$,
we have a decomposition\[
Y(K_{f})=\mathbf{G}(\mathbf{Q})\backslash\mathbf{G}(\mathbf{A})/K_{f}K_{\infty}^{\circ}\simeq\coprod_{x_{i}\in\gamma(K_{f})}\Gamma(x_{i})\backslash D_{\infty},\]
where $D_{\infty}=\mathbf{G}(\mathbf{R})^{\circ}/K_{\infty}^{\circ}$
is the symmetric space associated with $\mathbf{G}$ and $\Gamma(x_{i})\cong\tilde{\Gamma}(x_{i})/Z(K_{f})$
denotes the image of $\Gamma(x_{i})$ in the adjoint group. If $N$
is any left $K_{f}$-module, the double quotient\[
\widetilde{N}=\mathbf{G}(\mathbf{Q})\backslash\left(D_{\infty}\times\mathbf{G}(\mathbf{A}_{f})\times N\right)/K_{f}\]
naturally gives rise to a local system on $Y(K_{f})$, which is trivial
unless $Z(K_{f})$ acts trivially on $N$. Set $D_{\mathbf{A}}=D_{\infty}\times\mathbf{G}(\mathbf{A}_{f})$,
and let $C_{\bullet}(D_{\mathbf{A}})$ denote the complex of singular
chains on $D_{\mathbf{A}}$ endowed with the natural bi-action of
$\mathbf{G}(\mathbf{Q})\times\mathbf{G}(\mathbf{A}_{f})$. If $M$
and $N$ are right and left $K_{f}$-modules, respectively, we define
the complexes of \emph{adelic chains }and \emph{adelic cochains }by\[
C_{\bullet}^{ad}(K_{f},M)=C_{\bullet}(D_{\mathbf{A}})\otimes_{\mathbf{Z}[\mathbf{G}(\mathbf{Q})\times K_{f}]}M\]
and\[
C_{ad}^{\bullet}(K_{f},N)=\mathrm{Hom}_{\mathbf{Z}[\mathbf{G}(\mathbf{Q})\times K_{f}]}(C_{\bullet}(D_{\mathbf{A}}),N),\]
and we define functors $H_{\ast}(K_{f},-)$ and $H^{\ast}(K_{f},-)$
as their cohomology.

\textbf{Proposition 2.1.1. }\emph{There is a canonical isomorphism\[
H^{\ast}(Y(K_{f}),\widetilde{N})\simeq H^{\ast}(K_{f},N)=H^{\ast}(C_{ad}^{\bullet}(K_{f},N)).\]
}

\emph{Proof. }Let $C_{\bullet}(D_{\infty})(x_{i})$ denote the complex
of singular chains on $D_{\infty}$, endowed with the natural left
action of $\Gamma(x_{i})$ induced from the left action of $G(\mathbf{Q})^{\circ}$
on $D_{\infty}$; since $D_{\infty}$ is contractible, this is a free
resolution of $\mathbf{Z}$ in the category of $\mathbf{Z}[\Gamma(x_{i})]$-modules.
Let $N(x_{i})$ denote the left $\Gamma(x_{i})$-module whose underlying
module is $N$ but with the action $\gamma\cdot_{x_{i}}n=x_{i}^{-1}\gamma x_{i}|n$.
Note that the local system $\widetilde{N}(x_{i})$ obtained by restricting
$\widetilde{N}$ to the connected component $\Gamma(x_{i})\backslash D_{\infty}$
of $Y(K_{f})$ is simply the quotient $\Gamma(x_{i})\backslash\left(D_{\infty}\times N(x_{i})\right)$.
Setting \[
C_{sing}^{\bullet}(K_{f},N)=\oplus_{i}\mathrm{Hom}_{\mathbf{Z}[\Gamma(x_{i})]}(C_{\bullet}(D_{\infty})(x_{i}),N(x_{i})),\]
the map $D_{\infty}\to(D_{\infty},x_{i})\subset D_{\mathbf{A}}$ induces
a morphism $x_{i}^{\ast}=\mathrm{Hom}(C_{\bullet}(D_{\mathbf{A}}),N)\to\mathrm{Hom}(C_{\bullet}(D_{\infty}),N)$,
which in turn induces an isomorphism\[
\oplus_{i}x_{i}^{\ast}:C_{ad}^{\bullet}(K_{f},N)\overset{\sim}{\to}\oplus_{i}\mathrm{Hom}_{\Gamma(x_{i})}(C_{\bullet}(D_{\infty})(x_{i}),N(x_{i})),\]
and passing to cohomology we have\begin{eqnarray*}
H^{\ast}(C_{ad}^{\bullet}(K_{f},N)) & \simeq & \oplus_{i}H^{\ast}(\Gamma(x_{i})\backslash D_{\infty},\widetilde{N}(x_{i}))\\
 & \simeq & H^{\ast}(Y(K_{f}),\widetilde{N})\end{eqnarray*}
as desired. $\square$

When $\Gamma(x_{i})$ is torsion-free for each $x_{i}\in\gamma(K_{f})$,
we choose a finite resolution $F_{\bullet}(x_{i})\to\mathbf{Z}\to0$
of $\mathbf{Z}$ by free left $\mathbf{Z}[\Gamma(x_{i})]$-modules
of finite rank as well as a homotopy equivalence $F_{\bullet}(x_{i})\overset{f_{i}}{\underset{g_{i}}{\rightleftarrows}}C_{\bullet}(D_{\infty})(x_{i})$.
We shall refer to the resolution $F_{\bullet}(x_{i})$ as a \emph{Borel-Serre
resolution}; the existence of such resolutions follows from taking
a finite simplicial decomposition of the Borel-Serre compactification
of $\Gamma(x_{i})\backslash D_{\infty}$ \cite{BScorners}. Setting
\[
C_{\bullet}(K_{f},N)=\oplus_{i}F_{\bullet}(x_{i})\otimes_{\mathbf{Z}[\Gamma(x_{i})]}M(x_{i})\]
and\[
C^{\bullet}(K_{f},N)=\oplus_{i}\mathrm{Hom}_{\mathbf{Z}[\Gamma(x_{i})]}(F_{\bullet}(x_{i}),N(x_{i})),\]
the maps $f_{i},g_{i}$ induce homotopy equivalences\[
C_{\bullet}(K_{f},M)\overset{f_{\ast}}{\underset{g_{\ast}}{\rightleftarrows}}C_{\bullet}^{ad}(K_{f},M)\]
and\[
C^{\bullet}(K_{f},N)\overset{g^{\ast}}{\underset{f^{\ast}}{\rightleftarrows}}C_{ad}^{\bullet}(K_{f},M).\]
We refer to the complexes $C_{\bullet}(K_{f},-)$ and $C^{\bullet}(K_{f},-)$
as \emph{Borel-Serre complexes}, and we refer to these complexes together
with a \emph{fixed }set of homotopy equivalences $\{f_{i},g_{i}\}$
as \emph{augmented Borel-Serre complexes. }When the $\Gamma(x_{i})$'s
are not torsion-free but $M$ is uniquely divisible as a $\mathbf{Z}$-module,
we may still define $C_{\bullet}(K_{f},M)$ in an ad hoc manner by
taking the $K_{f}/K_{f}'$-coinvariants of $C_{\bullet}(K_{f}',M)$
for some sufficiently small normal subgroup $K'_{f}\subset K_{f}$.

\subsubsection*{Hecke operators}

A \emph{Hecke pair }consists of a monoid $\Delta\subset\mathbf{G}(\mathbf{A}_{f})$
and a subgroup $K_{f}\subset\Delta$ such that $K_{f}$ and $\delta K_{f}\delta^{-1}$
are commensurable for all $\delta\in\Delta$. Given a Hecke pair and
a commutative ring $R$, we write $\mathbf{T}(\Delta,K_{f})_{R}$
for the $R$-algebra generated by the double coset operators $T_{\delta}=[K_{f}\delta K_{f}]$
under convolution.%
\footnote{The ring structure on $\mathbf{T}(\Delta,K_{f})_{R}$ is nicely explained
in §3.1 of \cite{ShimuraIntro}.%
}

Suppose $M$ is a right $R[\Delta]$-module. The complex $C_{\bullet}(D_{\mathbf{A}})\otimes_{\mathbf{Z}[\mathbf{G}(\mathbf{Q})]}M$
receives a right $\Delta$-action via $(\sigma\otimes m)|\delta=\sigma\delta\otimes m\delta$,
and $C_{\bullet}^{ad}(K_{f},M)$ is naturally identified with the
$K_{f}$-coinvariants of this action. Given any double coset $K_{f}\delta K_{f}=\coprod_{j}\delta_{j}K_{f}$,
the action defined on pure tensors by the formula\[
(\sigma\otimes m)\cdot[K_{f}\delta K_{f}]=\sum_{j}(\sigma\otimes m)|\delta_{j}\]
induces a well-defined algebra homomorphism\[
\xi:\mathbf{T}(\Delta,K_{f})_{R}\to\mathrm{End}_{\mathbf{Ch}(R)}(C_{\bullet}^{ad}(K_{f},M)).\]
This action induces the usual Hecke action defined by correspondences
on homology. Set $\tilde{T}=g_{\ast}\circ\xi(T)\circ f_{\ast}\in\mathrm{End}_{\mathbf{Ch}(R)}(C_{\bullet}(K_{f},M))$.
The map\[
\tilde{\xi}:\xymatrix{\mathbf{T}(\Delta,K_{f})_{R}\ar[rr]^{T\mapsto\tilde{T}\qquad\qquad\qquad\qquad} &  & \mathrm{End}_{\mathbf{Ch}(R)}(C_{\bullet}(K_{f},M))\twoheadrightarrow\mathrm{End}_{\mathbf{K}(R)}(C_{\bullet}(K_{f},M))}
\]
is a well-defined ring homomorphism, since $g_{\ast}\circ\xi(T_{1})\circ f_{\ast}\circ g_{\ast}\circ\xi(T_{2})\circ f_{\ast}$
is homotopic to $g_{\ast}\circ\xi(T_{1}T_{2})\circ f_{\ast}$. Note
that any individual lift $\tilde{T}$ is well-defined in $\mathrm{End}_{\mathbf{Ch}(R)}(C_{\bullet}(K_{f},M))$,
but if $T_{1}$ and $T_{2}$ commute in the abstract Hecke algebra,
$\tilde{T}_{1}\tilde{T}_{2}$ and $\tilde{T}_{2}\tilde{T}_{1}$ will
typically only commute up to homotopy.

Likewise, if $N$ is a left $R[\Delta]$-module, the complex $\mathrm{Hom}_{\mathbf{Z}[\mathbf{G}(\mathbf{Q})]}(C_{\bullet}(D_{\mathbf{A}}),N)$
receives a natural $\Delta$-action via the formula $\delta|\phi=\delta\cdot\phi(\sigma\delta)$,
and $C_{ad}^{\bullet}(K_{f},N)$ is naturally the $K_{f}$-invariants
of this action. The formula \[
[K_{f}\delta K_{f}]\cdot\phi=\sum_{j}\delta_{j}|\phi\]
yields an algebra homomorphism $\xi:\mathbf{T}(\Delta,K_{f})_{R}\to\mathrm{End}_{R}(C_{ad}^{\bullet}(K_{f},N))$
which induces the usual Hecke action on cohomology, and $f^{\ast}\circ\xi\circ g^{\ast}$
defines an algebra homomorphism $\mathbf{T}(\Delta,K_{f})_{R}\to\mathrm{End}_{\mathbf{K}(R)}(C^{\bullet}(K_{f},M))$.
It is extremely important for us that these Hecke actions are compatible
with the duality isomorphism\[
\mathrm{Hom}_{R}(C_{\bullet}(K_{f},M),P)\simeq C^{\bullet}(K_{f},\mathrm{Hom}_{R}(M,P)),\]
where $P$ is any $R$-module.

We shall be mostly concerned with the following Hecke algebras. For
$I$, $\Lambda$ and $\Delta$ as in the beginning of §2, set $\mathcal{A}_{p}^{+}=\mathbf{T}(\Delta,I)_{\mathbf{Q}_{p}}$.
For any $t\in T^{+}$, the double coset operator $U_{t}=[ItI]$ defines
an element of $\mathcal{A}_{p}^{+}$, and the map $\Lambda\ni t\mapsto U_{t}\in\mathcal{A}_{p}^{+}$
extends to a commutative ring isomorphism\begin{eqnarray*}
\mathbf{Q}_{p}[\Lambda] & \overset{\sim}{\to} & \mathcal{A}_{p}^{+}\\
\sum c_{t}t & \mapsto & \sum c_{t}U_{t}.\end{eqnarray*}
The operators $U_{t}=[ItI]$ are invertible in the full Iwahori-Hecke
algebra $\mathbf{T}(\mathbf{G}(\mathbf{Q}_{p}),I)_{\mathbf{Q}_{p}}$ \cite{IM},
and we define the \emph{Atkin-Lehner algebra} $\mathcal{A}_{p}$ as
the commutative subalgebra of $\mathbf{T}(\mathbf{G}(\mathbf{Q}_{p}),I)_{\mathbf{Q}_{p}}$
generated by elements of the form $U_{t}$ and $U_{t}^{-1}$ for $t\in\Lambda$.
There is a natural ring isomorphism $\mathbf{\mathcal{A}}_{p}\simeq\mathbf{Q}_{p}[T(\mathbf{Q}_{p})/T(\mathbf{Z}_{p})]$,
though note that $t\cdot T(\mathbf{Z}_{p})$ typically corresponds
to the operator $U_{t_{1}}U_{t_{2}}^{-1}$ where $t_{1},t_{2}\in\Lambda$
are any elements with with $t_{1}t_{2}^{-1}\in t\cdot T(\mathbf{Z}_{p})$.
A \emph{controlling operator} is an element of $\mathcal{A}_{p}$
of the form $U_{t}$ for $t\in\Lambda^{+}$.

Fix an open compact subgroup $K^{p}\subset\mathbf{G}(\mathbf{A}_{f}^{p})$.
We say $K^{p}$ is unramified at a place $v\nmid p$ if $\mathbf{H}/F_{v}$
is unramified and $K_{v}^{p}=K^{p}\cap\mathbf{H}(F_{v})$ is a hyperspecial
maximal compact subgroup of $\mathbf{H}(F_{v})$, and we say $K^{p}$
is ramified otherwise. Let $S=S(K^{p})$ denote the finite set of
places where $K^{p}$ is ramified or $v|p$, and set $K_{S}^{p}=K^{p}\cap\prod_{v\in S}\mathbf{G}(F_{v})$,
so $K^{p}$ admits a product decomposition $K^{p}=K_{S}^{p}\prod_{v\notin S}K_{v}^{p}$.
We mainly work with the (commutative) Hecke algebras \begin{eqnarray*}
\mathbf{T}^{p}(K^{p}) & = & \bigotimes_{v\notin S(K^{p})}'\mathbf{T}(\mathbf{H}(F_{v}),K_{v}^{p})_{\mathbf{Q}_{p}},\\
\mathbf{T}(K^{p}) & = & \mathcal{A}_{p}^{+}\otimes\mathbf{T}^{p}(K^{p}),\end{eqnarray*}
In words, $\mathbf{T}(K^{p})$ takes into account the prime-to-\emph{p
}spherical Hecke operators together with certain Atkin-Lehner operators
at $p$; we write $\mathbf{T}_{\mathbf{G}}(K^{p})$ if we need to
emphasize $\mathbf{G}$. We also set $\mathbf{T}_{\mathrm{ram}}(K^{p})=\mathbf{T}\left(\prod_{v\in S,v\nmid p}\mathbf{H}(F_{v}),K_{S}^{p}\right)$.

\subsection{Locally analytic modules}

For each $s\geq1$ fix an analytic isomorphism $\psi^{s}:\mathbf{Z}_{p}^{d}\simeq\overline{N}^{s}$,
$d=\mathrm{dim}N$.%
\footnote{Use the homomorphisms $x_{\alpha}:U_{\alpha}\to\mathbf{G}_{a}$ together
with the product decomposition $\overline{N}\simeq\prod_{\alpha\in\Phi^{-}}U_{\alpha}$
for $\alpha$ in some fixed ordering. %
}

\textbf{Definition. }\emph{If $R$ is any $\mathbf{Q}_{p}$-Banach
algebra and $s$ is a positive integer, the module $\mathbf{A}(\overline{N}^{1},R)^{s}$
of} $s$-locally analytic $R$-valued functions on $\overline{N}^{1}$\emph{
is the $R$-module of continuous functions $f:\overline{N}^{1}\to R$
such that \[
f\left(x\psi^{s}(z_{1},\dots,z_{d})\right):\mathbf{Z}_{p}^{d}\to R\]
is given by an element of the $d$-variable Tate algebra $T_{d,R}=R\left\langle z_{1},\dots,z_{d}\right\rangle $
for any fixed $x\in\overline{N}^{1}$.}

Letting $\left\Vert \bullet\right\Vert _{T_{d,R}}$ denote the canonical
norm on the Tate algebra, the norm $\left\Vert f(x\psi^{s})\right\Vert _{T_{d,R}}$
depends only on the image of $x$ in $\overline{N}^{1}/\overline{N}^{s}$,
and the formula\[
\left\Vert f\right\Vert _{s}=\mathrm{sup}_{x\in\overline{N}^{1}}\left\Vert f(x\psi^{s})\right\Vert _{T_{d,R}}\]
defines a Banach $R$-module structure on $\mathbf{A}(\overline{N}^{1},R)^{s}$,
with respect to which the canonical inclusion $\mathbf{A}(\overline{N}^{1},R)^{s}\subset\mathbf{A}(\overline{N}^{1},R)^{s+1}$
is compact.

Given a tame level group $K^{p}\subset\mathbf{G}(\mathbf{A}_{f}^{p})$,
let $\overline{Z(K^{p}I)}$ be the $p$-adic closure of $Z(K^{p}I)$
in $T(\mathbf{Z}_{p})$. The\emph{ weight space of level $K^{p}$
}is the rigid analytic space $\sw=\sw_{K^{p}}$ over $\mathbf{Q}_{p}$
such that for any $\mathbf{Q}_{p}$-affinoid algebra $A$, $\mathrm{Hom}(\mathrm{Sp}A,\sw_{K^{p}})$
represents the functor which associates with $A$ the set of $p$-adically
continuous characters $\chi:T(\mathbf{Z}_{p})\to A^{\times}$ trivial
on $\overline{Z(K^{p}I)}$. It's not hard to check that $\sw$ is
the rigid space associated with the formal scheme $\mathrm{Spf}\mathbf{Z}_{p}[[T(\mathbf{Z}_{p})/\overline{Z(K^{p}I)}]]$
via Raynaud's generic fiber functor (cf. §7 of \cite{deJ}). Given
an admissible affinoid open $\Omega\subset\sw$, we write $\chi_{\Omega}:T(\mathbf{Z}_{p})\to\sco(\Omega)^{\times}$
for the unique character it determines. We define $s[\Omega]$ as
the minimal integer such that $\chi_{\Omega}|_{T^{s[\Omega]}}$ is
analytic. For any integer $s\geq s[\Omega]$, we make the definition\[
\mathbf{A}_{\Omega}^{s}=\left\{ f:I\to\sco(\Omega),\, f\,\mathrm{analytic\, on\, each\,}I^{s}-\mathrm{coset},\, f(gtn)=\chi_{\Omega}(t)f(g)\,\forall n\in N(\mathbf{Z}_{p}),\, t\in T(\mathbf{Z}_{p}),\, g\in I\right\} .\]
By the Iwahori decomposition, restricting an element $f\in\mathbf{A}_{\Omega}^{s}$
to $\overline{N}^{1}$ induces an isomorphism\begin{eqnarray*}
\mathbf{A}_{\Omega}^{s} & \simeq & \mathbf{A}(\overline{N}^{1},\sco(\Omega))^{s}\\
f & \mapsto & f|_{\overline{N}^{1}},\end{eqnarray*}
and we regard $\mathbf{A}_{\Omega}^{s}$ as a Banach $\sco(\Omega)$-module
via pulling back the Banach module structure on $\mathbf{A}(N(\mathbf{Z}_{p}),\sco(\Omega))^{s}$
under this isomorphism. The rule $(f|\gamma)(g)=f(\gamma g)$ gives
$\mathbf{A}_{\Omega}^{s}$ the structure of a continuous right $\sco(\Omega)[I]$-module.
More generally, the formula\[
\delta\star(nB(\mathbf{Z}_{p}))=\delta n\delta^{-1}\sigma(\delta)B(\mathbf{Z}_{p}),\, n\in\overline{N}^{1}\simeq I/B(\mathbf{Z}_{p})\,\mathrm{and}\,\delta\in T^{+}\]
yields a left action of $\Delta$ on $I/B(\mathbf{Z}_{p})$ which
extends the natural left translation action by $I$ (cf.\emph{ }§2.5
of \cite{AS}) and induces a right $\Delta$-action on $\mathbf{A}_{\Omega}^{s}$
which we denote by $f\star\delta,\, f\in\mathbf{A}_{\Omega}^{s}$.
For any $\delta\in T^{++}$, the operator $\delta\star-\in\mathcal{L}_{\sco(\Omega)}(\mathbf{A}_{\Omega}^{s},\mathbf{A}_{\Omega}^{s})$
factors through the inclusion $\mathbf{A}_{\Omega}^{s-1}\hookrightarrow\mathbf{A}_{\Omega}^{s}$,
and so defines a compact operator on $\mathbf{A}_{\Omega}^{s}$. The
Banach dual \begin{eqnarray*}
\mathbf{D}_{\Omega}^{s} & = & \mathcal{L}_{\sco(\Omega)}(\mathbf{A}_{\Omega}^{s},\sco(\Omega))\\
 & \simeq & \mathcal{L}_{\sco(\Omega)}(\mathbf{A}(\overline{N}^{1},\mathbf{Q}_{p})^{s}\widehat{\otimes}_{\mathbf{Q}_{p}}\sco(\Omega),\sco(\Omega))\\
 & \simeq & \mathcal{L}_{\mathbf{Q}_{p}}(\mathbf{A}(\overline{N}^{1},\mathbf{Q}_{p})^{s},\sco(\Omega))\end{eqnarray*}
inherits a dual left action of $\Delta$, and the operator $-\star\delta$
for $\delta\in T^{++}$ likewise factors through the inclusion $\mathbf{D}_{\Omega}^{s+1}\hookrightarrow\mathbf{D}_{\Omega}^{s}$.

We define an ind-Banach module \[
\sca_{\Omega}=\lim_{s\to\infty}\mathbf{A}_{\Omega}^{s}\]
where the direct limit is taken with respect to the natural compact,
injective transition maps $\mathbf{A}_{\Omega}^{s}\to\mathbf{A}_{\Omega}^{s+1}$.
Note that $\sca_{\Omega}$ is topologically isomorphic to the module
of $\sco(\Omega)$-valued locally analytic functions on $\overline{N}^{1}$,
equipped with the finest locally convex topology for which the natural
maps $\mathbf{A}_{\Omega}^{s}\hookrightarrow\sca_{\Omega}$ are continuous.
The $\Delta$-actions on $\mathbf{A}_{\Omega}^{s}$ induce a continuous
$\Delta$-action on $\sca_{\Omega}$. Set\[
\scd_{\Omega}=\left\{ \mu:\sca_{\Omega}\to\sco(\Omega),\,\mu\,\mathrm{is}\,\sco(\Omega)-\mathrm{linear\, and\, continuous}\right\} ,\]
and topologize $\scd_{\Omega}$ via the coarsest locally convex topology
for which the natural maps $\scd_{\Omega}\to\mathbf{D}_{\Omega}^{s}$
are continuous. In particular, the canonical map\[
\scd_{\Omega}\to\lim_{\infty\leftarrow s}\mathbf{D}_{\Omega}^{s}\]
is a topological isomorphism of locally convex $\sco(\Omega)$-modules,
and $\scd_{\Omega}$ is compact and Fréchet. Note that the transition
maps $\mathbf{D}_{\Omega}^{s+1}\to\mathbf{D}_{\Omega}^{s}$ are \emph{injective},
so $\scd_{\Omega}=\cap_{s\gg0}\mathbf{D}_{\Omega}^{s}$.

Suppose $\Sigma\subset\Omega$ is a Zariski closed subspace; by Corollary
9.5.2/8 of \cite{BGR}, $\Sigma$ arises from a surjection $\sco(\Omega)\twoheadrightarrow\sco(\Sigma)$
with $\sco(\Sigma)$ an affinoid algebra. We make the definitions
$\mathbf{D}_{\Sigma}^{s}=\mathbf{D}_{\Omega}^{s}\otimes_{\sco(\Omega)}\sco(\Sigma)$
and $\scd_{\Sigma}=\scd_{\Omega}\otimes_{\sco(\Omega)}\sco(\Sigma)$.

\textbf{Proposition 2.2.1. }\emph{There are canonical topological
isomorphisms $\mathbf{D}_{\Sigma}^{s}\simeq\mathcal{L}_{\sco(\Omega)}(\mathbf{A}_{\Omega}^{s},\sco(\Sigma))$
and $\scd_{\Sigma}\simeq\mathcal{L}_{\sco(\Omega)}(\sca_{\Omega},\sco(\Sigma))$.}

\emph{Proof. }Set $\mathfrak{a}_{\Sigma}=\ker(\sco(\Omega)\to\sco(\Sigma))$,
so $\sco(\Sigma)\simeq\sco(\Omega)/\mathfrak{a}_{\Sigma}$. The definitions
immediately imply isomorphisms\begin{eqnarray*}
\mathbf{D}_{\Sigma}^{s} & \simeq & \mathcal{L}_{\sco(\Omega)}(\mathbf{A}_{\Omega}^{s},\sco(\Omega))/\mathfrak{a}_{\Sigma}\mathcal{L}_{\sco(\Omega)}(\mathbf{A}_{\Omega}^{s},\sco(\Omega))\\
 & \simeq & \mathcal{L}_{\sco(\Omega)}(\mathbf{A}_{\Omega}^{s},\sco(\Omega))/\mathcal{L}_{\sco(\Omega)}(\mathbf{A}_{\Omega}^{s},\mathfrak{a}_{\Sigma}),\end{eqnarray*}
so the first isomorphism will follow if we can verify that the sequence
\[
0\to\mathcal{L}_{\sco(\Omega)}(\mathbf{A}_{\Omega}^{s},\mathfrak{a}_{\Sigma})\to\mathcal{L}_{\sco(\Omega)}(\mathbf{A}_{\Omega}^{s},\sco(\Omega))\to\mathcal{L}_{\sco(\Omega)}(\mathbf{A}_{\Omega}^{s},\sco(\Sigma))\]
is exact on the right. Given a $\mathbf{Q}_{p}$-Banach space $E$,
write $b(E)$ for the Banach space of bounded sequences $\{(e_{i})_{i\in\mathbb{N}},\mathrm{sup}_{i\in\mathbb{N}}\left|e_{i}\right|_{E}<\infty\}$.
Choosing an orthonormal basis of $\mathbf{A}(\overline{N}^{1},\sco(\Omega))^{s}$
gives rise to an isometry $\mathcal{L}_{\sco(\Omega)}(\mathbf{A}_{\Omega}^{s},E)\simeq b(E)$
for $E$ any Banach $\sco(\Omega)$-module. Thus we need to show the
surjectivity of the reduction map $b(\sco(\Omega))\to b(\sco(\Sigma))$.
Choose a presentation $\sco(\Omega)=T_{n}/\mathfrak{b}_{\Omega}$,
so $\sco(\Sigma)=T_{n}/\mathfrak{b}_{\Sigma}$ with $\mathfrak{b}_{\Omega}\subseteq\mathfrak{b}_{\Sigma}$.
Quite generally for any $\mathfrak{b}\subset T_{n}$, the function\[
f\in T_{n}/\mathfrak{b}\mapsto\left\Vert f\right\Vert _{\mathfrak{b}}=\mathrm{inf}_{\tilde{f}\in f+\mathfrak{b}}\left\Vert \tilde{f}\right\Vert _{T_{n}}\]
defines a norm on $T_{n}/\mathfrak{b}$. By Proposition 3.7.5/3 of \cite{BGR},
there is a unique Banach algebra structure on any affinoid algebra.
Hence for any sequence $(f_{i})_{i\in\mathbb{N}}\in b(\sco(\Sigma))$,
we may choose a bounded sequence of lifts $(\widetilde{f}_{i})_{i\in\mathbb{N}}\in b(T_{n})$;
reducing the latter sequence modulo $\mathfrak{b}_{\Omega}$, we are
done.

Taking inverse limits in the sequence we just proved to be exact,
the second isomorphism follows. $\square$

Suppose $\lambda\in X_{+}^{\ast}\subset\scw$ is a dominant weight
for $B$, with $\scl_{\lambda}$ the corresponding irreducible left
$G(\mathbf{Q}_{p})$-representation of highest weight $\lambda$.
We may realize $\scl_{\lambda}$ explicitly as\[
\scl_{\lambda}(L)=\left\{ f:G\to L\,\mathrm{algebraic},\, f(n'tg)=\lambda(t)f(g)\,\mathrm{for}\, n'\in\overline{N},t\in T,g\in G\right\} \]
with $G$ acting by right translation. The function $f_{\lambda}(g)$
defined by $f_{\lambda}(n'tn)=\lambda(t)$ on the big cell extends
uniquely to an algebraic function on $G$, and is the highest weight
vector in $\scl_{\lambda}$. For $g\in G(\mathbf{Q}_{p})$ and $h\in I$,
the function $f_{\lambda}(gh)$ defines an element of $\scl_{\lambda}\otimes\sca_{\lambda}$,
and pairing it against $\mu\in\scd_{\lambda}$ defines a map $i_{\lambda}:\scd_{\lambda}\to\scl_{\lambda}$
which we notate suggestively as\[
i_{\lambda}(\mu)(g)=\int f_{\lambda}(gh)\mu(h).\]
The map $\mu\mapsto i_{\lambda}(\mu)(g)$ satisfies the following
intertwining relation for $\gamma\in I$:\begin{eqnarray*}
\gamma\cdot i_{\lambda}(\mu)(g) & = & i_{\lambda}(\mu)(g\gamma)\\
 & = & \int f_{\lambda}(g\gamma n^{\circ})\mu(n^{\circ})\\
 & = & \int f_{\lambda}(n^{\circ}g)(\gamma\cdot\mu)(n^{\circ})\\
 & = & i_{\lambda}(\gamma\cdot\mu)(g).\end{eqnarray*}

\subsubsection*{The case of $\mathrm{GL}_{n}/\mathbf{Q}_{p}$}

We examine the case when $G\simeq\mathrm{GL}_{n}/\mathbf{Q}_{p}$.
We choose $B$ and $\overline{B}$ as the upper and lower triangular
Borel subgroups, respectively, and we identify $T$ with diagonal
matrices. The splitting $\sigma$ is canonically induced from the
homomorphism\begin{eqnarray*}
\mathbf{Q}_{p}^{\times} & \to & \mathbf{Z}_{p}^{\times}\\
x & \mapsto & p^{-v_{p}(x)}x.\end{eqnarray*}
Since $T(\mathbf{Z}_{p})\simeq(\mathbf{Z}_{p}^{\times})^{n}$, we
canonically identify a character $\lambda:T(\mathbf{Z}_{p})\to R^{\times}$
with the $n$-tuple of characters $(\lambda_{1},\dots,\lambda_{n})$
where \begin{eqnarray*}
\lambda_{i}:\mathbf{Z}_{p}^{\times} & \to & R^{\times}\\
x & \mapsto & \lambda\circ\mathrm{diag}(\underset{i-1}{\underbrace{1,\dots,1}},x,1,\dots,1).\end{eqnarray*}
Dominant weights are identified with tuples of integers $(k_{1},\dots,k_{n})$
with $k_{1}\geq k_{2}\geq\dots\geq k_{n}$, by associating to such
a tuple the character with $\lambda_{i}(x)=x^{k_{i}}$.

We want to explain how to {}``twist away'' one dimension's worth
of weights in a canonical fashion. For any $\lambda\in\sw$, a simple
calculation shows that the $\star$-action of $\Delta$ on $\mathbf{A}_{\lambda}^{s}$
is given explicitly by the formula\[
(f\star\delta)(x)=\lambda(\sigma(\mathrm{t}(\delta))\mathrm{t}(\delta)^{-1}\mathrm{t}(\delta x))f(\overline{\mathrm{n}}(\delta x)),\,\delta\in\Delta,\, x\in\overline{N}^{1},\, f\in\mathbf{A}(\overline{N}^{1},k(\lambda))^{s}.\]
Given $1\leq i\leq n$, let $m_{i}(g)$ denote the determinant of
the upper-left $i$-by-$i$ block of $g\in\mathrm{GL}_{n}$. For any
$g\in\Delta$, a pleasant calculation left to the reader shows that\[
\mathrm{t}(g)=\mathrm{diag}(m_{1}(g),m_{1}(g)^{-1}m_{2}(g),\dots,m_{i}^{-1}(g)m_{i+1}(g),\dots,m_{n-1}(g)^{-1}\mathrm{det}g).\]
In particular, writing $\lambda^{0}=(\lambda_{1}\lambda_{n}^{-1},\lambda_{2}\lambda_{n}^{-1},\dots,\lambda_{n-1}\lambda_{n}^{-1},1)$
yields a canonical isomorphism\[
\mathbf{A}_{\lambda}^{s}\simeq\mathbf{A}_{\lambda^{0}}^{s}\otimes\lambda_{n}(\det\cdot|\det|_{p})\]
of $\Delta$-modules, and likewise for $\mathbf{D}_{\lambda}^{s}$.

In the case of $\mathrm{GL}_{2}$ we can be even more explicit. Here
$\Delta$ is generated by the center of $G(\mathbf{Q}_{p})$ and by
the monoid\[
\Sigma_{0}(p)=\left\{ g=\left(\begin{array}{cc}
a & b\\
c & d\end{array}\right)\in M_{2}(\mathbf{Z}_{p}),\, c\in p\mathbf{Z}_{p},\, a\in\mathbf{Z}_{p}^{\times},\, ad-bc\neq0\right\} .\]
Another simple calculation shows that the center of $G(\mathbf{Q}_{p})$
acts on $\mathbf{A}_{\lambda}^{s}$ through the character $z\mapsto\lambda(\sigma(z)),$
while the monoid $\Sigma_{0}(p)$ acts via\[
(g\cdot f)(x)=(\lambda_{1}\lambda_{2}^{-1})(a+bx)\lambda_{2}(\det g|\det g|_{p})f\left(\frac{c+dx}{a+bx}\right),\, f\in\mathbf{A}(\overline{N}^{1},k)^{s},\,\left(\begin{array}{cc}
1\\
x & 1\end{array}\right)\in\overline{N}^{1},\]
(almost) exactly as in \cite{Strigidsymbs}.

\paragraph*{Remarks. }

There are some subtle differences between the different modules we
have defined. The assignment $\Omega\mapsto\sca_{\Omega}$ describes
an LB sheaf over $\sw$, and the modules $\mathbf{A}_{\Omega}^{s}$
are orthonormalizable. On the other hand, the modules $\mathbf{D}_{\Omega}^{s}$
are \emph{not }obviously orthonormalizable. There are alternate modules
of distributions available, namely $\tilde{\mathbf{D}}_{\Omega}^{s}=\mathcal{L}(\mathbf{A}(\overline{N}^{1},\mathbf{Q}_{p})^{s},\mathbf{Q}_{p})\widehat{\otimes}_{\mathbf{Q}_{p}}\sco(\Omega)$
equipped with the unique action for which the inclusion $\tilde{\mathbf{D}}_{\Omega}^{s}\hookrightarrow\mathbf{D}_{\Omega}^{s}$
is equivariant. The module $\tilde{\mathbf{D}}_{\Omega}^{s}$ is orthonormalizable
but of course is not the continuous dual of $\mathbf{A}_{\Omega}^{s}$.
However, in the inverse limit these differences disappear: there is
a compact injective map $j_{s}:\mathbf{D}_{\Omega}^{s+1}\to\tilde{\mathbf{D}}_{\Omega}^{s}$
such that the diagram\[
\xymatrix{\tilde{\mathbf{D}}_{\Omega}^{s+1}\ar[r]\ar[d] & \tilde{\mathbf{D}}_{\Omega}^{s}\ar[d]\\
\mathbf{D}_{\Omega}^{s+1}\ar[r]\ar[ur]^{j_{s}} & \mathbf{D}_{\Omega}^{s}}
\]
commutes, and this easily yields a canonical topological isomorphism
\[
\scd_{\Omega}=\lim_{\infty\leftarrow s}\mathbf{D}_{\Omega}^{s}\cong\lim_{\infty\leftarrow s}\tilde{\mathbf{D}}_{\Omega}^{s}.\]
In particular, the assignment $\Omega\mapsto\scd_{\Omega}$ defines
a Fréchet sheaf on $\scw$ as well. One of our goals is to demonstrate
the feasibility of working successfully with the modules $\mathbf{D}_{\Omega}^{s}$
by treating the dual modules $\mathbf{A}_{\Omega}^{s}$ on an equal
footing.

\subsection{Slope decompositions of modules and complexes}

Here we review the very general notion of slope decomposition introduced
in \cite{AS}. Let $A$ be a $\mathbf{Q}_{p}$-Banach algebra, and
let $M$ be an $A$-module equipped with an $A$-linear endomorphism
$u:M\to M$ (for short, {}``an $A[u]$-module''). Fix a rational
number $h\in\mathbf{Q}_{\geq0}$. We say a polynomial $Q\in A[x]$
is \emph{multiplicative} if the leading coefficient of $Q$ is a unit
in $A$, and that $Q$ has \emph{slope $\leq h$} if every edge of
the Newton polygon of $Q$ has slope $\leq h$. Write $Q^{\ast}(x)=x^{\deg Q}Q(1/x)$.
An element $m\in M$ has slope $\leq h$ if there is a multiplicative
polynomial $Q\in A[T]$ of slope $\leq h$ such that $ $$Q^{\ast}(u)\cdot m=0$.
Let $M_{\leq h}$ be the set of elements of $M$ of slope $\leq h$;
according to Proposition 4.6.2 of \emph{loc. cit.}, $M_{\leq h}$
is an $A$-submodule of $M$.

\textbf{Definition 2.3.1. }A \emph{slope-$\leq h$ decomposition }of
$M$ is an $A[u]$-module isomorphism\[
M\simeq M_{\leq h}\oplus M_{>h}\]
such that $M_{\leq h}$ is a finitely generated $A$-module and the
map $Q^{\ast}(u):M_{>h}\to M_{>h}$ is an $A$-module isomorphism
for every multiplicative polynomial $Q\in A[T]$ of slope $\leq h$.

The following proposition summarizes the fundamental results on slope
decompositions.

\textbf{Proposition 2.3.2 (Ash-Stevens): }
\begin{description}
\item [{\emph{a)}}] \emph{Suppose $M$ and $N$ are both $A[u]$-modules
with slope-$\leq h$ decompositions. If $\psi:M\to N$ is a morphism
of $A[u]$-modules, then $\psi(M_{\leq h})\subseteq N_{\leq h}$ and
$\psi(M_{>h})\subseteq N_{>h}$; in particular, a module can have
at most one slope-$\leq h$ decomposition. Furthermore, $\ker\psi$
and $\mathrm{im}\psi$ inherit slope-$\leq h$ decompositions. Given
a short exact sequence \[
0\to M\to N\to L\to0\]
of $A[u]$-modules, if two of the modules admit slope-$\leq h$ decompositions
then so does the third. }
\item [{\emph{b)}}] \emph{If $C^{\bullet}$ is a complex of $A[u]$-modules,
all with slope-$\leq h$ decompositions, then \[
H^{n}(C^{\bullet})\simeq H^{n}(C_{\leq h}^{\bullet})\oplus H^{n}(C_{>h}^{\bullet})\]
is a slope-$\leq h$ decomposition of $H^{n}(C^{\bullet})$.}
\end{description}
\emph{Proof. }This is a rephrasing of (a specific case of) Proposition
4.1.2 of \cite{AS}. $\square$

Suppose now that $A$ is a reduced affinoid algebra, $M$ is an orthonormalizable
Banach $A$-module, and $u$ is a compact operator. Let \[
F(T)=\mathrm{det}(1-uT)|M\in A[[T]]\]
denote the Fredholm determinant for the $u$-action on $M$. We say
$F$ admits a \emph{slope-$\leq h$} \emph{factorization} if we can
write $F(T)=Q(T)\cdot R(T)$ where $Q$ is a multiplicative polynomial
of slope $\leq h$ and $R(T)\in A[[T]]$ is an entire power series
of slope $>h$. Theorem 3.3 of \cite{BuEigen} guarantees that $F$
admits a slope-$\leq h$ factorization if and only if $M$ admits
a slope-$\leq h$ decomposition. Furthermore, given a slope-$\leq h$
factorization $F(T)=Q(T)\cdot R(T)$, we obtain the slope-$\le h$
decomposition of $M$ upon setting $M_{\leq h}=\left\{ m\in M|Q^{\ast}(u)\cdot m=0\right\} $,
and $M_{\leq h}$ in this case is a finite flat $A$-module upon which
$u$ acts invertibly.%
\footnote{Writing $Q^{\ast}(x)=a+x\cdot r(x)$ with $r\in A[x]$ and $a\in A^{\times}$,
$u^{-1}$ on $M_{\leq h}$ is given explicitly by $-a^{-1}r(u)$.%
} Combining this with Theorem 4.5.1 of \cite{AS} and Proposition 2.3.1,
we deduce:

\textbf{Proposition 2.3.3. }\emph{If $C^{\bullet}$ is a bounded complex
of orthonormalizable Banach $A[u]$-modules, and $u$ acts compactly
on the total complex $\oplus_{i}C^{i}$, then for any $x\in\mathrm{Max}(A)$
and any $h\in\mathbf{Q}_{\geq0}$ there is an affinoid subdomain $\mathrm{Max}(A')\subset\mathrm{Max}(A)$
containing $x$ such that the complex $C^{\bullet}\widehat{\otimes}_{A}A'$
of $A'[u]$-modules admits a slope-$\leq h$ decomposition, and $(C^{\bullet}\widehat{\otimes}_{A}A')_{\leq h}$
is a complex of finite flat $A'$-modules.}

\textbf{Proposition 2.3.4. }\emph{If $M$ is an orthonormalizable
Banach $A$-module with a slope-$\leq h$ decomposition, and $A'$
is a Banach $A$-algebra, then $M\widehat{\otimes}_{A}A'$ admits
a slope-$\leq h$ decomposition and in fact\[
(M\widehat{\otimes}_{A}A')_{\leq h}\simeq M_{\leq h}\otimes_{A}A'.\]
}

\textbf{Proposition 2.3.5. }\emph{If $N\in\mathrm{Ban}_{A}$ is finite
and $M\in\mathrm{Ban}_{A}$ is an $A[u]$-module with a slope-$\leq h$
decomposition, the $A[u]$-modules $M\widehat{\otimes}_{A}N$ and
$\mathcal{L}_{A}(M,N)$ inherit slope-$\leq h$ decompositions.}

\emph{Proof. }This is an immediate consequence of the $A$-linearity
of the $u$-action and the fact that $-\widehat{\otimes}_{A}N$ and
$\mathcal{L}_{A}(-,N)$ commute with finite direct sums. $\square$

If $A$ is a field and $M$ is either an orthonormalizable Banach
$A$-module or the cohomology of a complex of such, then $M$ admits
a slope-$\leq h$ decomposition for every $h$, and if $h<h'$ there
is a natural decomposition\[
M_{\leq h'}\simeq M_{\leq h}\oplus\left(M_{>h}\right)_{\leq h'}\]
and in particular a projection $M_{\leq h'}\to M_{\leq h}$. We set
$M^{\mathrm{fs}}=\lim_{\infty\leftarrow h}M_{\leq h}$.

\section{Overconvergent cohomology}

Fix a connected, reductive group $\mathbf{G}/\mathbf{Q}$ as in the
introduction. For any tame level group $K^{p}\subset\mathbf{G}(\mathbf{A}_{f}^{p})$,
we abbreviate $H_{\ast}(K^{p}I,-)$ by $H_{\ast}(K^{p},-)$, and likewise
for cohomology.

\subsection{Basic results}

In this section we establish some foundational results on overconvergent
cohomology. These results likely follow from the formalism introduced
in Chapter 5 of \cite{AS}, but we give different proofs. The key
idea exploited here, namely the lifting of the $U_{t}$-action to
the level of chain complexes, is due to Ash. We use freely the notations
introduced in §2.1-§2.3.

Fix an augmented Borel-Serre complex $C_{\bullet}(K^{p},-)=C_{\bullet}(K^{p}I,-)$.
Fix an element $t\in\Lambda^{+}$, and let $\tilde{U}=\tilde{U}_{t}$
denote the lifting of $U_{t}=[ItI]$ to an endomorphism of the complex
$C_{\bullet}(K^{p},-)$ defined in §2.1. Given a connected admissible
open affinoid subset $\Omega\subset\sw_{K^{p}}$ and any integer $s\geq s[\Omega]$,
the endomorphism $\tilde{U}_{t}\in\mathrm{End}_{\sco(\Omega)}(C_{\bullet}(K^{p},\mathbf{A}_{\Omega}^{s}))$
is compact; let \[
F_{\Omega}^{s}(X)=\det(1-X\tilde{U}_{t})|C_{\bullet}(K^{p},\mathbf{A}_{\Omega}^{s})\in\sco(\Omega)[[X]]\]
denote its Fredholm determinant. We say $(U_{t},\Omega,h)$ is a \emph{slope
datum} if $C_{\bullet}(K^{p},\mathbf{A}_{\Omega}^{s})$ admits a slope-$\leq h$
decomposition for the $\tilde{U}_{t}$ action for some $s\geq s[\Omega]$.

\textbf{Proposition 3.1.1. }\emph{The function $F_{\Omega}^{s}(X)$
is independent of $s$.}

\emph{Proof. }For any integer $s\geq s[\Omega]$ we write $C_{\bullet}^{s}=C_{\bullet}(K^{p},\mathbf{A}_{\Omega}^{s})$
for brevity. By construction, the operator $\tilde{U}_{t}$ factors
into compositions $\rho_{s}\circ\check{U}_{t}$ and $\check{U}_{t}\circ\rho_{s+1}$
where $\check{U}_{t}:C_{\bullet}^{s}\to C_{\bullet}^{s-1}$ is continuous
and $\rho_{s}:C_{\bullet}^{s-1}\to C_{\bullet}^{s}$ is compact. Now,
considering the commutative diagram\[
\xymatrix{C_{\bullet}^{s}\ar[r]^{\check{U}_{t}}\ar[d]^{\tilde{U}_{t}} & C_{\bullet}^{s-1}\ar[dl]_{\rho_{s}}\ar[d]^{\tilde{U}_{t}}\\
C_{\bullet}^{s}\ar[r]_{\check{U}_{t}} & C_{\bullet}^{s-1}}
\]
we calculate\begin{eqnarray*}
\det(1-X\tilde{U}_{t})|C_{\bullet}^{s} & = & \det(1-X\rho_{s}\circ\check{U}_{t})|C_{\bullet}^{s}\\
 & = & \det(1-X\check{U}_{t}\circ\rho_{s})|C_{\bullet}^{s-1}\\
 & = & \det(1-X\tilde{U}_{t})|C_{\bullet}^{s-1},\end{eqnarray*}
where the second line follows from Lemma 2.7 of \cite{BuEigen}, so
$F_{\Omega}^{s}(X)=F_{\Omega}^{s-1}(X)$ for all $s>s[\Omega]$. $\square$

\textbf{Proposition 3.1.2. }\emph{The slope-$\leq h$ subcomplex $C_{\bullet}(K^{p},\mathbf{A}_{\Omega}^{s})_{\leq h}$,
if it exists, is independent of $s$. If $\Omega'$ is an affinoid
subdomain of $\Omega$, then the restriction map $\mathbf{A}_{\Omega}^{s}\to\mathbf{A}_{\Omega'}^{s}$
induces a canonical isomorphism\[
C_{\bullet}(K^{p},\mathbf{A}_{\Omega}^{s})_{\leq h}\otimes_{\sco(\Omega)}\sco(\Omega')\simeq C_{\bullet}(K^{p},\mathbf{A}_{\Omega'}^{s})_{\leq h}\]
for any $s\geq s[\Omega]$.}

\emph{Proof. }Since $F_{\Omega}^{s}(X)$ is independent of $s$, we
simply write $F_{\Omega}(X)$. Suppose we are given a slope-$\leq h$
factorization $F_{\Omega}(X)=Q(X)\cdot R(X)$; by the remarks in §2.3,
setting $C_{\bullet}(K^{p},\mathbf{A}_{\Omega}^{s})_{\leq h}=\ker Q^{\ast}(\tilde{U}_{t})$
yields a slope-$\leq h$ decomposition of $C_{\bullet}(K^{p},\mathbf{A}_{\Omega}^{s})$
for any $s\geq s[\Omega]$. By Proposition 2.3.1, the injection $\rho_{s}:C_{\bullet}(K^{p},\mathbf{A}_{\Omega}^{s-1})\hookrightarrow C_{\bullet}(K^{p},\mathbf{A}_{\Omega}^{s})$
gives rise to a canonical injection \[
\rho_{s}:C_{\bullet}(K^{p},\mathbf{A}_{\Omega}^{s-1})_{\leq h}\hookrightarrow C_{\bullet}(K^{p},\mathbf{A}_{\Omega}^{s})_{\leq h}\]
for any $s>s[\Omega]$. The operator $\tilde{U}_{t}$ acts invertibly
on $C_{\bullet}(K^{p},\mathbf{A}_{\Omega}^{s})_{\leq h}$, and its
image factors through $\rho_{s}$, so $\rho_{s}$ is surjective and
hence bijective. This proves the first claim.

For the second claim, by Proposition 2.3.3 we have\begin{eqnarray*}
C_{\bullet}(K^{p},\mathbf{A}_{\Omega}^{s})_{\leq h}\otimes_{\sco(\Omega)}\sco(\Omega') & \simeq & \left(C_{\bullet}(K^{p},\mathbf{A}_{\Omega}^{s})\widehat{\otimes}_{\sco(\Omega)}\sco(\Omega')\right)_{\leq h}\\
 & \simeq & C_{\bullet}(K^{p},\mathbf{A}_{\Omega'}^{s})_{\leq h},\end{eqnarray*}
so the result now follows from the first claim. $\square$

\textbf{Proposition 3.1.3. }\emph{Given a slope datum $(U_{t},\Omega,h)$
and an affinoid subdomain $\Omega'\subset\Omega$, there is a canonical
isomorphism\[
H_{\ast}(K^{p},\mathbf{A}_{\Omega}^{s})_{\leq h}\otimes_{\sco(\Omega)}\sco(\Omega')\simeq H_{\ast}(K^{p},\mathbf{A}_{\Omega'}^{s})_{\leq h}\]
for any $s\geq s[\Omega]$.}

\emph{Proof. }Since $\sco(\Omega')$ is $\sco(\Omega)$-flat, the
functor $-\otimes_{\sco(\Omega)}\sco(\Omega')$ commutes with taking
the homology of any complex of $\sco(\Omega)$-modules. Thus we calculate\begin{eqnarray*}
H_{\ast}(K^{p},\mathbf{A}_{\Omega}^{s})_{\leq h}\otimes_{\sco(\Omega)}\sco(\Omega') & \simeq & H_{\ast}\left(C_{\bullet}(K^{p},\mathbf{A}_{\Omega}^{s})_{\leq h}\right)\otimes_{\sco(\Omega)}\sco(\Omega')\\
 & \simeq & H_{\ast}\left(C_{\bullet}(K^{p},\mathbf{A}_{\Omega}^{s})_{\leq h}\otimes_{\sco(\Omega)}\sco(\Omega')\right)\\
 & \simeq & H_{\ast}\left(C_{\bullet}(K^{p},\mathbf{A}_{\Omega'}^{s})_{\leq h}\right)\\
 & \simeq & H_{\ast}(K^{p},\mathbf{A}_{\Omega'}^{s})_{\leq h},\end{eqnarray*}
where the third line follows from Proposition 2.3.4. $\square$

\textbf{Proposition 3.1.4. }\emph{Given a slope datum $(U_{t},\Omega,h)$,
the complex $C_{\bullet}(K^{p},\sca_{\Omega})$ and the homology module
$H_{\ast}(K^{p},\sca_{\Omega})$ admit slope-$\leq h$ decompositions,
and there is an isomorphism \[
H_{\ast}(K^{p},\sca_{\Omega})_{\leq h}\simeq H_{\ast}(K^{p},\mathbf{A}_{\Omega}^{s})_{\leq h}\]
for any $s\geq s[\Omega]$. Furthermore, given an affinoid subdomain
$\Omega'\subset\Omega$, there is a canonical isomorphism\[
H_{\ast}(K^{p},\sca_{\Omega})_{\leq h}\otimes_{\sco(\Omega)}\sco(\Omega')\simeq H_{\ast}(K^{p},\sca_{\Omega'})_{\leq h}.\]
}

\emph{Proof. }For any fixed $s\geq s[\Omega]$, we calculate\begin{eqnarray*}
C_{\bullet}(K^{p},\sca_{\Omega}) & \simeq & \lim_{\substack{\to\\
s'}
}C_{\bullet}(K^{p},\mathbf{A}_{\Omega}^{s'})\\
 & \simeq & \lim_{\substack{\to\\
s'}
}C_{\bullet}(K^{p},\mathbf{A}_{\Omega}^{s'})_{\leq h}\oplus C_{\bullet}(K^{p},\mathbf{A}_{\Omega}^{s'})_{>h}\\
 & \simeq & C_{\bullet}(K^{p},\mathbf{A}_{\Omega}^{s})_{\leq h}\oplus\lim_{\substack{\to\\
s'}
}C_{\bullet}(K^{p},\mathbf{A}_{\Omega}^{s'})_{>h}\end{eqnarray*}
with the third line following from Proposition 3.1.2. The two summands
in the third line naturally form the components of a slope-$\leq h$
decomposition, so passing to homology yields the first sentence of
the proposition, and the second sentence then follows immediately
from Proposition 2.3.3. $\square$

We're now in a position to prove the subtler cohomology analogue of
Proposition 3.1.4.

\textbf{Proposition 3.1.5. }\emph{Given a slope datum $(U_{t},\Omega,h)$
and a Zariski-closed subspace $\Sigma\subset\Omega$, the complex
$C^{\bullet}(K^{p},\scd_{\Sigma})$ and the cohomology module $H^{\ast}(K^{p},\scd_{\Sigma})$
admit slope-$\leq h$ decompositions, and there is an isomorphism
\[
H^{\ast}(K^{p},\scd_{\Sigma})_{\leq h}\simeq H^{\ast}(K^{p},\mathbf{D}_{\Sigma}^{s})_{\leq h}\]
for any $s\geq s[\Omega]$. Furthermore, given an affinoid subdomain
$\Omega'\subset\Omega$, there are canonical isomorphisms\[
C^{\bullet}(K^{p},\scd_{\Omega})_{\leq h}\otimes_{\sco(\Omega)}\sco(\Omega')\simeq C^{\bullet}(K^{p},\scd_{\Omega'})_{\leq h}\]
and\[
H^{\ast}(K^{p},\scd_{\Omega})_{\leq h}\otimes_{\sco(\Omega)}\sco(\Omega')\simeq H^{\ast}(K^{p},\scd_{\Omega'})_{\leq h}.\]
}

\emph{Proof. }By a topological version of the duality stated in §2.1,
we have a natural isomorphism\begin{eqnarray*}
C^{\bullet}(K^{p},\mathbf{D}_{\Sigma}^{s}) & = & C^{\bullet}(K^{p},\mathcal{L}_{\sco(\Omega)}(\mathbf{A}_{\Omega}^{s},\sco(\Sigma)))\\
 & \simeq & \mathcal{L}_{\sco(\Omega)}(C_{\bullet}(K^{p},\mathbf{A}_{\Omega}^{s}),\sco(\Sigma))\end{eqnarray*}
for any $s\geq s[\Omega]$. By assumption, $C_{\bullet}(K^{p},\mathbf{A}_{\Omega}^{s})$
admits a slope-$\leq h$ decomposition, so we calculate\begin{eqnarray*}
C^{\bullet}(K^{p},\mathbf{D}_{\Sigma}^{s}) & \simeq & \mathcal{L}_{\sco(\Omega)}(C_{\bullet}(K^{p},\mathbf{A}_{\Omega}^{s}),\sco(\Sigma))\\
 & \simeq & \mathcal{L}_{\sco(\Omega)}(C_{\bullet}(K^{p},\mathbf{A}_{\Omega}^{s})_{\leq h},\sco(\Sigma))\\
 &  & \oplus\mathcal{L}_{\sco(\Omega)}(C_{\bullet}(K^{p},\mathbf{A}_{\Omega}^{s})_{>h},\sco(\Sigma)).\end{eqnarray*}
By Proposition 3.1.2, passing to the inverse limit over $s$ in this
isomorphism yields a slope-$\leq h$ decomposition of $C^{\bullet}(K^{p},\scd_{\Sigma})$
together with a natural isomorphism \[
C^{\bullet}(K^{p},\scd_{\Sigma})_{\leq h}\simeq C^{\bullet}(K^{p},\mathbf{D}_{\Sigma}^{s})_{\leq h}\simeq\mathcal{L}_{\sco(\Omega)}(C^{\bullet}(K^{p},\mathbf{A}_{\Omega}^{s})_{\leq h},\sco(\Sigma))\]
for any $s\geq s[\Omega]$. This proves the first sentence of the
proposition.

For the second sentence, we first note that since $C_{\bullet}(K^{p},\mathbf{A}_{\Omega}^{s})_{\leq h}$
is a complex of finite $\sco(\Omega)$-modules, the natural map\[
\mathcal{L}_{\sco(\Omega)}(C_{\bullet}(K^{p},\mathbf{A}_{\Omega}^{s})_{\leq h},\sco(\Omega))\to\mathrm{Hom}_{\sco(\Omega)}(C_{\bullet}(K^{p},\mathbf{A}_{\Omega}^{s})_{\leq h},\sco(\Omega))\]
is an isomorphism by Lemma 2.2 of \cite{BuEigen}. Next, note that
if $R$ is a commutative ring, $S$ is a flat $R$-algebra, and $M,N$
are $R$-modules with $M$ finitely presented, the natural map $\mathrm{Hom}_{R}(M,N)\otimes_{R}S\to\mathrm{Hom}_{S}(M\otimes_{R}S,N\otimes_{R}S)$
is an isomorphism. With these two facts in hand, we calculate as follows:
\begin{eqnarray*}
C^{\bullet}(K^{p},\scd_{\Omega})_{\leq h}\otimes_{\sco(\Omega)}\sco(\Omega') & \simeq & \mathrm{Hom}_{\sco(\Omega)}(C_{\bullet}(K^{p},\mathbf{A}_{\Omega}^{s})_{\leq h},\sco(\Omega))\otimes_{\sco(\Omega)}\sco(\Omega')\\
 & \simeq & \mathrm{Hom}_{\sco(\Omega')}(C_{\bullet}(K^{p},\mathbf{A}_{\Omega}^{s})_{\leq h}\otimes_{\sco(\Omega)}\sco(\Omega'),\sco(\Omega'))\\
 & \simeq & \mathrm{Hom}_{\sco(\Omega')}(C_{\bullet}(K^{p},\mathbf{A}_{\Omega'}^{s})_{\leq h},\sco(\Omega'))\\
 & \simeq & C^{\bullet}(K^{p},\scd_{\Omega'})_{\leq h},\end{eqnarray*}
where the third line follows from Proposition 2.3.3. Passing to cohomology,
the result follows as in the proof of Proposition 3.1.3. $\square$

\subsection{Finite-slope eigenpackets and non-critical classes}

In this section we explain two results which are fundamental in our
analysis. First of all, we recall and summarize some of the work of
Eichler, Shimura, Matsushima, Borel-Wallach, Franke, and Li-Schwermer
on the cohomology of arithmetic groups. Next, we state a fundamental
theorem of Ash-Stevens and Urban (Theorem 6.4.1 of \cite{AS}, Proposition
4.3.10 of \cite{UrEigen}) relating overconvergent cohomology classes
of small slope with classical automorphic forms. The possibility of
such a result was largely the original \emph{raison d'etre} of overconvergent
cohomology; in the case $\mathbf{G}=\mathrm{GL}_{2}/\mathbf{Q}$,
Stevens proved this theorem in a famous preprint \cite{Strigidsymbs}.

Let $\lambda\in X_{+}^{\ast}$ be a $B$-dominant algebraic weight,
and let $K_{f}\subset\mathbf{G}(\mathbf{A}_{f})$ be any open compact
subgroup. By fundamental work of Franke, the cohomology $H^{\ast}(Y(K_{f}),\scl_{\lambda})_{\mathbf{C}}=H^{\ast}(Y(K_{f}),\scl_{\lambda})\otimes_{\mathbf{Q}_{p},\iota}\mathbf{C}$
admits an analytically defined splitting\[
H^{\ast}(Y(K_{f}),\scl_{\lambda})_{\mathbf{C}}\simeq H_{\mathrm{cusp}}^{\ast}(Y(K_{f}),\scl_{\lambda})_{\mathbf{C}}\oplus H_{\mathrm{Eis}}^{\ast}(Y(K_{f}),\scl_{\lambda})_{\mathbf{C}}\]
into $\mathbf{T}(\mathbf{G}(\mathbf{A}_{f}),K_{f})_{\mathbf{C}}$-stable
submodules, which we refer to as the cuspidal and Eisenstein cohomology,
respectively. The cuspidal cohomology admits an exact description
in terms of cuspidal automorphic forms \cite{BWcohom,FrankeL2,FrSchwermer}:

\textbf{Proposition 3.2.1. }\emph{There is a canonical isomorphism\[
H_{\mathrm{cusp}}^{\ast}(Y(K_{f}),\scl_{\lambda})_{\mathbf{C}}\simeq\bigoplus_{\pi\in L_{\mathrm{cusp}}^{2}\left(\mathbf{G}(\mathbf{Q})\backslash\mathbf{G}(\mathbf{A})\right)}m(\pi)\pi_{f}^{K_{f}}\otimes H^{\ast}(\mathfrak{g},K_{\infty};\pi_{\infty}\otimes\scl_{\lambda})\]
of graded $\mathbf{T}(\mathbf{G}(\mathbf{A}_{f}),K_{f})$-modules,
where $m(\pi)$ denotes the multiplicity of $\pi$ in $L_{\mathrm{cusp}}^{2}\left(\mathbf{G}(\mathbf{Q})\backslash\mathbf{G}(\mathbf{A})\right)$.
If $\lambda$ is a regular weight, the natural inclusion of $H_{\mathrm{cusp}}^{\ast}(Y(K_{f}),\scl_{\lambda})$
into $H_{!}^{\ast}(Y(K_{f}),\scl_{\lambda})$ is an isomorphism}.

Note that if $\pi$ contributes nontrivially to the direct sum decomposition
of the previous proposition, the central and infininitesimal characters
of $\pi_{\infty}$ are necessarily inverse to those of $\scl_{\lambda}$.

For any weight $\lambda\in\sw_{K^{p}}(\overline{\mathbf{Q}_{p}})$
and a given controlling operator $U_{t}$, we define $\mathbf{T}_{\lambda,h}(K^{p})$
as the subalgebra of $\mathrm{End}_{k_{\lambda}}\left(H^{\ast}(K^{p},\scd_{\lambda})_{\leq h}\right)$
generated by the image of $\mathbf{T}(K^{p})\otimes_{\mathbf{Q}_{p}}k_{\lambda}$,
and we set $\mathbf{T}_{\lambda}(K^{p})=\lim_{\infty\leftarrow h}\mathbf{T}_{\lambda,h}(K^{p})$.
The algebra $\mathbf{T}_{\lambda}(K^{p})$ is independent of the choice
of controlling operator used in its definition.

\textbf{Definition 3.2.2. }\emph{A }finite-slope eigenpacket of weight
$\lambda$ and level $K^{p}$ \emph{(or simply a }finite-slope eigenpacket\emph{)}
\emph{is an algebra homomorphism $\phi:\mathbf{T}_{\lambda}(K^{p})\to\overline{\mathbf{Q}_{p}}$.}

It's easy to check that this coincides with the definition given in
the introduction. If $\phi$ is a finite-slope eigenpacket, we shall
regard the contraction of $\ker\phi$ under the structure map $\mathbf{T}(K^{p})\to\mathbf{T}_{\lambda}(K^{p})$
as a maximal ideal in $\mathbf{T}(K^{p})$, which we also denote $ $by
$\ker\phi$. Note that $\mathbf{T}_{\lambda}(K^{p})$ is a countable
direct product of zero-dimensional Artinian local rings, and the factors
in this direct product are in natural bijection with the finite-slope
eigenpackets.

A weight $\lambda$ is \emph{arithmetic }if it factors as the product
of a finite-order character $\varepsilon$ of $T(\mathbf{Z}_{p})$
and an element $\lambda^{\mathrm{alg}}$ of $X^{\ast}$; if $\lambda^{\mathrm{alg}}\in X_{+}^{\ast}$
we say $\lambda$ is \emph{dominant arithmetic}. If $\lambda=\lambda^{\mathrm{alg}}\varepsilon$
is a dominant arithmetic weight, we are going to formulate some comparisons
between $H^{\ast}(K^{p},\scd_{\lambda})$ and $H^{\ast}(Y(K^{p}I_{1}^{s}),\scl_{\lambda^{\mathrm{alg}}})$.
In order to do this, we need to twist the natural Hecke action on
the latter module slightly. More precisely, if $M$ is any $\mathcal{A}_{p}^{+}$-module
and $\lambda\in X^{\ast}$, we define the \emph{$\star$-action in
weight $\lambda$ }by $U_{t}\star_{\lambda}m=\lambda(t)^{-1}U_{t}m$ \cite{AS}.
The map $i_{\lambda}$ defined in §2.2 induces a morphism\[
i_{\lambda}:H^{\ast}(K^{p},\scd_{\lambda})\to H^{\ast}(K^{p}I_{1}^{s},\scl_{\lambda^{\mathrm{alg}}})\]
for any $s\geq s[\varepsilon]$ which intertwines the standard action
of $\mathbf{T}(K^{p})$ on the source with the $\star$-action on
the target. 

\textbf{Definition 3.2.3. }\emph{If $\lambda=\lambda^{\mathrm{alg}}\varepsilon$
is an arithmetic weight and $s=s[\varepsilon]$, a finite-slope eigenpacket
is }classical \emph{if the module $H^{\ast}(Y(K^{p}I_{1}^{s}),\scl_{\lambda^{\mathrm{alg}}})$
is nonzero after localization at $\ker\phi$, and }noncritical \emph{if
the map\[
i_{\lambda}:H^{\ast}(K^{p},\scd_{\lambda})\to H^{\ast}(K^{p}I_{1}^{s},\scl_{\lambda^{\mathrm{alg}}})\]
becomes an isomorphism after localization at $\ker\phi$. A classical
eigenpacket is }interior \emph{if $H_{\partial}^{\ast}(K^{p}I_{1}^{s},\scl_{\lambda^{\mathrm{alg}}})$
vanishes after localization at $\ker\phi$, and }strongly interior\emph{
if $H_{\partial}^{\ast}(K^{p},\scd_{\lambda})$ vanishes after localization
at $\ker\phi$ as well.}

Next we formulate a result which generalizes Stevens's control theorem \cite{Strigidsymbs}.
Given $\lambda\in X^{\ast}$, we define an action of the Weyl group
$W$ by $w\cdot\lambda=(\lambda+\rho)^{w}-\rho$.

\textbf{Definition 3.2.4. }\emph{Fix a controlling operator $U_{t}$,
$t\in\Lambda^{+}$. Given an arithmetic weight $\lambda=\lambda^{\mathrm{alg}}\varepsilon$,
a rational number $h$ is a }small slope for $\lambda$ \emph{if \[
h<\inf_{w\in W\smallsetminus\{\mathrm{id}\}}v_{p}(w\cdot\lambda^{\mathrm{alg}}(t))-v_{p}(\lambda^{\mathrm{alg}}(t)).\]
}

\textbf{Theorem 3.2.5 (Ash-Stevens, Urban). }\emph{Fix an arithmetic
weight $\lambda=\lambda^{\mathrm{alg}}\varepsilon$ and a controlling
operator $U_{t}$. If $h$ is a small slope for $\lambda$, there
is a natural isomorphism of Hecke modules\[
H^{\ast}(K^{p},\scd_{\lambda})_{\leq h}\simeq H^{\ast}(Y(K^{p}I_{1}^{s}),\scl_{\lambda^{\mathrm{alg}}})_{\leq h}^{T(\mathbf{Z}/p^{s}\mathbf{Z})=\varepsilon}\]
for any $s\geq s[\varepsilon]$.}

\emph{Proof (sketch). }Suppose $\lambda=\lambda^{\mathrm{alg}}$ for
simplicity. Urban constructs a second quadrant spectral sequence\[
E_{1}^{i,j}=\bigoplus_{w\in W,\,\ell(w)=-i}H^{j}(K^{p},\scd_{w\cdot\lambda})^{\mathrm{fs}}\Rightarrow H^{i+j}(K^{p}I,\scl_{\lambda})^{\mathrm{fs}}\]
which is equivariant for $U_{t}$ if we twist the action as follows:
$U_{t}$ acts through the $\star$-action in weight $\lambda$ on
the target and $(\lambda(t)^{-1}w\cdot\lambda)(t))U_{t}$ acts on
the $w$-summand of the $E_{1}$-page. In particular, taking the slope-$\leq h$
part yields a spectral sequence\[
E_{1}^{i,j}=\bigoplus_{w\in W,\,\ell(w)=-i}H^{j}(K^{p},\scd_{w\cdot\lambda^{\mathrm{alg}}})_{\leq h-v_{p}\left(w\cdot\lambda(t)\right)+v_{p}(\lambda(t))}\Rightarrow H^{i+j}(K^{p}I,\scl_{\lambda})_{\leq h}.\]
But any element of $\mathcal{A}_{p}^{+}$ is contractive on $H^{j}(K^{p},\scd_{\lambda})$,
so the $w$-summand of the $E_{1}$-page is now empty for $w\neq\mathrm{id}$.
$\square$

\subsection{The spectral sequences}

In this section we introduce our main technical tool for analyzing
overconvergent cohomology. Fix a choice of tame level $K^{p}$ and
an augmented Borel-Serre complex $C_{\bullet}(K^{p},-)$. 

\textbf{Theorem 3.3.1. }\emph{Fix a slope datum $(U_{t},\Omega,h)$,
and let $\Sigma\subseteq\Omega$ be an arbitrary rigid Zariski closed
subspace. Then $H^{\ast}(K^{p},\scd_{\Sigma})$ admits a slope-$\leq h$
decomposition, and there is a convergent first quadrant spectral sequence\[
E_{2}^{i,j}=\mathrm{Ext}_{\sco(\Omega)}^{i}(H_{j}(K^{p},\sca_{\Omega})_{\leq h},\sco(\Sigma))\Rightarrow H^{i+j}(K^{p},\scd_{\Sigma})_{\leq h}.\]
Furthermore, there is a convergent second quadrant spectral sequence\[
E_{2}^{i,j}=\mathrm{Tor}_{-i}^{\sco(\Omega)}(H^{j}(K^{p},\scd_{\Omega})_{\leq h},\sco(\Sigma))\Rightarrow H^{i+j}(K^{p},\scd_{\Sigma})_{\leq h}.\]
In addition, there are analogous spectral sequences relating Borel-Moore
homology with compactly supported cohomology, and boundary homology
with boundary cohomology, and there are morphisms between the spectral
sequences compatible with the morphisms between these different cohomology
theories. Finally, the spectral sequences and the morphisms between
them are equivariant for the natural Hecke actions on their $E_{2}$
pages and abutments; more succinctly, they are spectral sequences
of $\mathbf{T}(K^{p})$-modules.}

When no ambiguity is likely, we will refer to the two spectral sequences
of Theorem 3.3.1 as {}``the Ext spectral sequence'' and {}``the
Tor spectral sequence.'' The Hecke-equivariance of these spectral
sequences is crucial for applications, since it allows one to localize
the entire spectral sequence at any ideal in the Hecke algebra.

\emph{Proof of Theorem 3.3.1. }By the isomorphisms proved in §3.1,
it suffices to construct a Hecke-equivariant spectral sequence $\mathrm{Ext}_{\sco(\Omega)}^{i}(H_{j}(K^{p},\mathbf{A}_{\Omega}^{s})_{\leq h},\sco(\Sigma))\Rightarrow H^{i+j}(K^{p},\mathbf{D}_{\Sigma}^{s})_{\leq h}$
for some $s\geq s[\Omega]$.

Consider the hyperext group $\mathbf{E}\mathrm{xt}_{\sco(\Omega)}^{n}(C_{\bullet}^{ad}(K^{p},\mathbf{A}_{\Omega}^{s}),\sco(\Sigma))$.
Since $C_{\bullet}^{ad}(K^{p},\mathbf{A}_{\Omega}^{s})$ is a complex
of $\mathbf{T}(K^{p})$-modules, this hyperext group is naturally
a $\mathbf{T}(K^{p})$-module, and the hypercohomology spectral sequence
\[
E_{2}^{i,j}=\mathrm{E}\mathrm{xt}_{\sco(\Omega)}^{i}(H_{j}(K^{p},\mathbf{A}_{\Omega}^{s}),\sco(\Sigma))\Rightarrow\mathbf{E}\mathrm{xt}_{\sco(\Omega)}^{i+j}(C_{\bullet}^{ad}(K^{p},\mathbf{A}_{\Omega}^{s}),\sco(\Sigma))\]
is a spectral sequence of $\mathbf{T}(K^{p})$-modules. On the other
hand, the quasi-isomorphism $C_{\bullet}^{ad}(K^{p},\mathbf{A}_{\Omega}^{s})\simeq C_{\bullet}(K^{p},\mathbf{A}_{\Omega}^{s})$
in $\mathbf{D}^{b}(A(\Omega))$ together with the slope-$\leq h$
decomposition $C_{\bullet}(K^{p},\mathbf{A}_{\Omega}^{s})\simeq C_{\bullet}(K^{p},\mathbf{A}_{\Omega}^{s})_{\leq h}\oplus C_{\bullet}(K^{p},\mathbf{A}_{\Omega}^{s})_{>h}$
induces Hecke-stable slope-$\leq h$-decompositions of the abutment
and of the entries on the $E_{2}$ page. By Proposition 2.3.2, the
slope decomposition of the $E_{2}$ page induces slope decompositions
of all entries on all higher pages of the spectral sequence. In other
words, we may pass to the {}``slope-$\leq h$ part'' of the hypercohomology
spectral sequence in a Hecke-equivariant way, getting a spectral sequence\[
'E_{2}^{i,j}=\mathrm{Ext}_{\sco(\Omega)}^{i}(H_{j}(K^{p},\mathbf{A}_{\Omega}^{s}),\sco(\Sigma))_{\leq h}\Rightarrow\mathbf{E}\mathrm{xt}_{\sco(\Omega)}^{i+j}(C_{\bullet}^{ad}(K^{p},\mathbf{A}_{\Omega}^{s}),\sco(\Sigma))_{\leq h}\]
of $\mathbf{T}(K^{p})$-modules. But $\mathrm{Ext}_{\sco(\Omega)}^{i}(H_{j}(K^{p},\mathbf{A}_{\Omega}^{s}),\sco(\Sigma))_{\leq h}\simeq\mathrm{Ext}_{\sco(\Omega)}^{i}(H_{j}(K^{p},\mathbf{A}_{\Omega}^{s})_{\leq h},\sco(\Sigma))$,
and\begin{eqnarray*}
\mathbf{E}\mathrm{xt}_{\sco(\Omega)}^{n}(C_{\bullet}^{ad}(K^{p},\mathbf{A}_{\Omega}^{s}),\sco(\Sigma))_{\leq h} & \simeq & \mathbf{E}\mathrm{xt}_{\sco(\Omega)}^{n}(C_{\bullet}(K^{p},\mathbf{A}_{\Omega}^{s})_{\leq h},\sco(\Sigma))\\
 & \simeq & H^{n}\left(\mathbf{R}\mathrm{Hom}_{\sco(\Omega)}(C_{\bullet}(K^{p},\mathbf{A}_{\Omega}^{s})_{\leq h},\sco(\Sigma))\right)\\
 & \simeq & H^{n}\left(\mathrm{Hom}_{\sco(\Omega)}(C_{\bullet}(K^{p},\mathbf{A}_{\Omega}^{s})_{\leq h},\sco(\Sigma))\right)\\
 & \simeq & H^{n}(K^{p},\mathbf{D}_{\Sigma}^{s})_{\le h},\end{eqnarray*}
where the third line follows from the projectivity of each $C_{i}(K^{p},\mathbf{A}_{\Omega}^{s})_{\leq h}$
and the fourth line follows from the proof of Proposition 3.1.5. 

For the Tor spectral sequence, the isomorphism\[
C^{\bullet}(K^{p},\mathbf{D}_{\Omega}^{s})_{\leq h}\otimes_{\sco(\Omega)}\sco(\Sigma)\simeq C^{\bullet}(K^{p},\mathbf{D}_{\Sigma}^{s})_{\leq h}\]
yields an isomorphism\[
C^{\bullet}(K^{p},\mathbf{D}_{\Omega}^{s})_{\leq h}\otimes_{\sco(\Omega)}^{\mathbf{L}}\sco(\Sigma)\simeq C^{\bullet}(K^{p},\mathbf{D}_{\Sigma}^{s})_{\leq h}\]
of $\mathbf{T}(K^{p})$-module complexes in $\mathbf{D}^{b}(\sco(\Omega))$,
and the result follows analogously from the hypertor spectral sequence\[
\mathrm{Tor}_{-i}^{R}(H^{j}(C^{\bullet}),N)\Rightarrow\mathbf{T}\mathrm{or}_{-i-j}^{R}(C^{\bullet},N).\]

\textbf{Remark 3.3.2. }If $(\Omega,h)$ is a slope datum, $\Sigma_{1}$
is Zariski-closed in $\Omega$, and $\Sigma_{2}$ is\textbf{ }Zariski-closed
in $\Sigma_{1}$, the transitivity of the derived tensor product yields
an isomorphism\begin{eqnarray*}
C^{\bullet}(K^{p},\scd_{\Sigma_{2}})_{\leq h} & \simeq & C^{\bullet}(K^{p},\scd_{\Omega})_{\leq h}\otimes_{\sco(\Omega)}^{\mathbf{L}}\sco(\Sigma_{2})\\
 & \simeq & C^{\bullet}(K^{p},\scd_{\Omega})_{\leq h}\otimes_{\sco(\Omega)}^{\mathbf{L}}\sco(\Sigma_{1})\otimes_{\sco(\Sigma_{1})}^{\mathbf{L}}\sco(\Sigma_{2})\\
 & \simeq & C^{\bullet}(K^{p},\scd_{\Sigma_{1}})_{\leq h}\otimes_{\sco(\Sigma_{1})}^{\mathbf{L}}\sco(\Sigma_{2})\end{eqnarray*}
which induces a relative version of the Tor spectral sequence, namely\[
E_{2}^{i,j}=\mathrm{Tor}_{\sco(\Sigma_{1})}^{-i}(H^{j}(K^{p},\scd_{\Sigma_{1}})_{\leq h},\sco(\Sigma_{2}))\Rightarrow H^{i+j}(K^{p},\scd_{\Sigma_{2}})_{\leq h}.\]
This spectral sequence plays an important role in Newton's proof of
Theorem 1.1.6.

\paragraph{The boundary and Borel-Moore/compactly supported spectral sequences}

Notation as in §2.1, let $\overline{D_{\infty}}$ denote the Borel-Serre
bordification of $D_{\infty}$, and let $\partial\overline{D_{\infty}}=\overline{D_{\infty}}\smallsetminus D_{\infty}$.
Setting $\overline{D_{\mathbf{A}}}=\overline{D_{\infty}}\times\mathbf{G}(\mathbf{A}_{f})$
and $\partial\overline{D_{\mathbf{A}}}=\partial\overline{D_{\infty}}\times\mathbf{G}(\mathbf{A}_{f})$,
the natural map $C_{\bullet}(D_{\mathbf{A}})\to C_{\bullet}(\overline{D_{\mathbf{A}}})$
induces a functorial isomorphism $H_{\ast}(K_{f},M)\simeq H_{\ast}\left(C_{\bullet}(\overline{D_{\mathbf{A}}})\otimes_{\mathbf{Z}[\mathbf{G}(\mathbf{Q})\times K_{f}]}M\right)$,
so we may redefine $C_{\bullet}^{ad}(K_{f},M)$ as $C_{\bullet}(\overline{D_{\mathbf{A}}})\otimes_{\mathbf{Z}[\mathbf{G}(\mathbf{Q})\times K_{f}]}M$.
Setting $C_{\bullet}^{\partial,ad}(K_{f},M)=C_{\bullet}(\partial\overline{D_{\mathbf{A}}})\otimes_{\mathbf{Z}[\mathbf{G}(\mathbf{Q})\times K_{f}]}M$,
the natural inclusion induces a map $C_{\bullet}^{\partial,ad}(K_{f},M)\to C_{\bullet}^{ad}(K_{f},M)$,
and we define $C_{\bullet}^{\mathrm{BM},ad}(K_{f},M)$ as the cone
of this map. Not surprisingly, the homology of $C_{\bullet}^{\partial}$
(resp. $C_{\bullet}^{\mathrm{BM}}$) computes boundary (resp. Borel-Moore)
homology. Choosing a triangulation of $\overline{Y(K_{f})}$ induces
a triangulation on the boundary, and yields a complex $C_{\bullet}^{\partial}(K_{f},,M)$
together with a map $ $$C_{\bullet}^{\partial}(K_{f},,M)\to C_{\bullet}^{\partial}(K_{f},,M)$;
defining $C_{\bullet}^{\mathrm{BM}}$ as the cone of this map, these
complexes all fit into a big diagram\[
\xymatrix{C_{\bullet}^{\partial,ad}(K_{f},-)\ar[r] & C_{\bullet}^{ad}(K_{f},-)\ar[r] & C_{\bullet}^{\mathrm{BM},ad}(K_{f},-)\\
C_{\bullet}^{\partial}(K_{f},-)\ar[u]\ar[r] & C_{\bullet}(K_{f},-)\ar[u]\ar[r] & C_{\bullet}^{\mathrm{BM}}(K_{f},-)\ar[u]}
\]
in which the rows are exact triangles functorial in $M$, and the
vertical arrows are quasi-isomorphisms. We make analogous definitions
of $C_{c,ad}^{\bullet}(K_{f},M)$, etc.

The boundary and Borel-Moore/compactly supported sequences, and the
morphisms between them, follow from {}``taking the slope-$\leq h$
part'' of the diagram\[
\xymatrix{\mathbf{R}\mathrm{Hom}_{\sco(\Omega)}(C_{\bullet}^{\mathrm{BM},ad}(K^{p},\mathbf{A}_{\Omega}^{s}),\sco(\Sigma))\ar[r]\ar[d] & C_{c,ad}^{\bullet}(K^{p},\mathrm{Hom}_{\sco(\Omega)}(\mathbf{A}_{\Omega}^{s},\sco(\Sigma)))\ar[d]\\
\mathbf{R}\mathrm{Hom}_{\sco(\Omega)}(C_{\bullet}^{ad}(K^{p},\mathbf{A}_{\Omega}^{s}),\sco(\Sigma))\ar[r]\ar[d] & C_{ad}^{\bullet}(K^{p},\mathrm{Hom}_{\sco(\Omega)}(\mathbf{A}_{\Omega}^{s},\sco(\Sigma)))\ar[d]\\
\mathbf{R}\mathrm{Hom}_{\sco(\Omega)}(C_{\bullet}^{\partial,ad}(K^{p},\mathbf{A}_{\Omega}^{s}),\sco(\Sigma))\ar[r] & C_{\partial,ad}^{\bullet}(K^{p},\mathrm{Hom}_{\sco(\Omega)}(\mathbf{A}_{\Omega}^{s},\sco(\Sigma)))}
\]
in which the horizontal arrows are quasi-isomorphisms, the columns
are exact triangles in $\mathbf{D}^{b}(\sco(\Omega))$, and the diagram
commutes for the natural action of $\mathbf{T}(K^{p})$.

\section{Eigenvarieties}

In this section we begin to use global rigid analytic geometry in
a more serious way; the main reference for this topic is \cite{BGR},
and \cite{ConSeveral} is a nice survey of the main ideas. We shall
repeatedly and tacitly use the fact that if $\Omega'$ is an affinoid
subdomain of an affinoid space $\Omega$, $\sco(\Omega')$ is a flat
$\sco(\Omega)$-module; this is an easy consequence of the universal
property of an affinoid subdomain together with the local criterion
for flatness.

\subsection{Fredholm hypersurfaces}

Let $A$ be an affinoid integral domain. We say that such an $ $$A$
is \emph{relatively factorial }if for any $f=\sum_{n=0}^{\infty}a_{n}X^{n}\in A\left\langle X\right\rangle $
with $a_{0}=1$, $(f)$ factors uniquely as a product of principal
prime ideals $(f_{i})$ where each $f_{i}$ may be chosen with constant
term $1$. A rigid analytic space $\sw$ is relatively factorial if
it has an admissible covering by relatively factorial affinoids. Throughout
the remainder of this article, we reserve the letter $\sw$ for a
relatively factorial rigid analytic space.

\textbf{Definition 4.1.1. }\emph{A }Fredholm series \emph{is a global
section $f\in\sco(\sw\times\mathbf{A}^{1})$ such that under the map
$\sco(\sw\times\mathbf{A}^{1})\overset{j^{\ast}}{\to}\sco(\sw)$ induced
by $j:\sw\times\{0\}\to\sw\times\mathbf{A}^{1}$ we have $j^{\ast}f=1$.
A }Fredholm hypersurface \emph{is a closed immersion $\sz\subset\sw\times\mathbf{A}^{1}$
such that the ideal sheaf of $\sz$ is generated by a Fredholm series
$f$, in which case we write $\sz=\sz(f)$.}

Note that the natural projection $\sw\times\mathbf{A}^{1}\to\sw$
induces a map $w:\sz\to\sw$. Let $\sco(\sw)\{\{X\}\}$ denote the
subring of $\sco(\sw)[[X]]$ consisting of series $\sum_{n=0}^{\infty}a_{n}X^{n}$
such that $|a_{n}|_{\Omega}r^{n}\to0$ as $n\to\infty$ for any affinoid
$\Omega\subset\sw$ and any $r\in\mathbf{R}_{>0}$. The natural injection
$\sco(\sw\times\mathbf{A}^{1})\simeq\sco(\sw)\{\{X\}\}\hookrightarrow\sco(\sw)[[X]]$
identifies the monoid of Fredholm series with elements of $\sco(\sw)\{\{X\}\}$
such that $a_{0}=1$. When $\sw$ is relatively factorial, the ring
$\sco(\sw)\{\{X\}\}$ admits a good factorization theory, and we may
speak of irreducible Fredholm series without ambiguity. We say a collection
of distinct irreducible Fredholm series $\{f_{i}\}_{i\in I}$ is \emph{locally
finite }if $\sz(f_{i})\cap U=\emptyset$ for all but finitely many
$i\in I$ and any quasi-compact admissible open subset $U\subset\sw\times\mathbf{A}^{1}$.

\textbf{Proposition 4.1.2 (Coleman-Mazur, Conrad). }\emph{If $\sw$
is relatively factorial, any Fredholm series $f$ admits a factorization
$f=\prod_{i\in I}f_{i}^{n_{i}}$ as a product of irreducible Fredholm
series with $n_{i}\geq1$; any such factorization is unique up to
reordering the terms, the collection $\{f_{i}\}_{i\in I}$ is locally
finite, and the irreducible components of $\sz(f)$ are exactly the
Fredholm hypersurfaces $\sz(f_{i}^{n_{i}})$. The nilreduction of
$\sz(f)$ is $\sz(\prod_{i\in I}f_{i})$.}

\emph{Proof. }See §1 of \cite{CMeigencurve} and §4 of \cite{Conirred}
(especially Theorems 4.2.2 and 4.3.2). $\square$

\textbf{Proposition 4.1.3. }\emph{If $\sz$ is a Fredholm hypersurface,
the image $w(\sz)$ is Zariski-open in $\sw$.}

\emph{Proof. }We may assume $\sz=\sz(f)$ with $f=1+\sum_{n=1}^{\infty}a_{n}X^{n}$
irreducible. By Lemma 1.3.2 of \cite{CMeigencurve}, the fiber of
$\sz$ over $\lambda\in\sw(\overline{\mathbf{Q}_{p}})$ is empty if
and only if $a_{n}\in\mathfrak{m}_{\lambda}$ for all $n$, if and
only if $\mathscr{I}=(a_{1},a_{2},a_{3},\dots)\subset\mathfrak{m}_{\lambda}$.
The ideal $\mathscr{I}$ is naturally identified with the global sections
of a coherent ideal sheaf, which cuts out a closed immersion $V(\mathscr{I})\hookrightarrow\sw$
in the usual way, and $w(\sz)$ is exactly the complement of $V(\mathscr{I})$.
$\square$

Given a Fredholm hypersurface $\sz=\sz(f)$, a rational number $h\in\mathbf{Q}$,
and an affinoid $\Omega\subset\sw$, we define $\sz_{\Omega,h}=\sco(\Omega)\left\langle p^{h}X\right\rangle /(f(X))$
regarded as an admissible affinoid open subset of $\sz$. The natural
map $\sz_{\Omega,h}\to\Omega$ is flat but not necessarily finite,
and we define an affinoid of the form $\sz_{\Omega,h}$ to be \emph{slope-adapted
}if $\sz_{\Omega,h}\to\Omega$ is a finite flat map. The affinoid
$\sz_{\Omega,h}$ is slope-adapted if and only if $f_{\Omega}=f|_{\sco(\Omega)\{\{X\}\}}$
admits a slope-$\leq h$ factorization $Q(X)\cdot R(X)$, in which
case $\sco(\sz_{\Omega,h})\simeq\sco(\Omega)[X]/(Q(X))$.

\textbf{Proposition 4.1.4. }\emph{For any Fredholm hypersurface $\sz$,
the collection of slope-adapted affinoids forms an admissible cover
of $\sz$.}

\emph{Proof. }See §4 of \cite{BuEigen}.

\subsection{Eigenvariety data}

\textbf{Definition 4.2.1. }\emph{An }eigenvariety datum \emph{is a
tuple $\mathfrak{D}=(\sw,\sz,\sm,\mathbf{T},\psi)$ where $\sw$ is
a separated, reduced, equidimensional, relatively factorial rigid
analytic space, $\sz\subset\sw\times\mathbf{A}^{1}$ is a Fredholm
hypersurface, $\sm$ is a coherent analytic sheaf on $\sz$, $\mathbf{T}$
is a commutative $\mathbf{Q}_{p}$-algebra, and $\psi$ is a $\mathbf{Q}_{p}$-algebra
homomorphism $\psi:\mathbf{T}\to\mathrm{End}_{\sco_{\sz}}(\sm)$.}

In practice $\mathbf{T}$ will be a Hecke algebra, $\sz$ will be
a {}``spectral variety'' parametrizing the eigenvalues of some distinguished
operator $U\in\mathbf{T}$ on some graded module $M^{\ast}$ of \emph{p-}adic
automorphic forms or on a complex whose cohomology yields $M^{\ast}$,
and $\sm$ will be the natural {}``spreading out'' of $M^{\ast}$
to a coherent sheaf over $\sz$. We do \emph{not }require that $\sm$
be locally free on $\sz$.

\textbf{Theorem 4.2.2. }\emph{Given an eigenvariety datum $\mathfrak{D}$,
there is a rigid analytic space $\sx=\sx(\mathfrak{D})$ together
with a finite morphism $\pi:\sx\to\sz$, a morphism $w:\sx\to\sw$,
an algebra homomorphism $\phi_{\sx}:\mathbf{T}\to\sco(\sx)$, and
a coherent sheaf $\sm^{\dagger}$ on $\sx$ together with a canonical
isomorphism $\sm\cong\pi_{\ast}\sm^{\dagger}$ compatible with the
actions of $\mathbf{T}$ on $\scm$ and $\scm^{\dagger}$ (via $\psi$ and $\phi_{\sx}$, respectively). The points
of $\sx$ lying over $z\in\sz$ are in bijection with the generalized
eigenspaces for the action of $\mathbf{T}$ on $\sm(z)$.}

\emph{Proof. }Let $\scc\mathrm{ov}=\left\{ \Omega_{i}\right\} _{i\in I}$
be an admissible affinoid cover of $\sz$; we abbreviate $\Omega_{i}\cap\Omega_{j}$
by $\Omega_{ij}$. For any $\Omega_{i}$ we let $\mathbf{T}_{\Omega_{i}}$
be the finite $\sco(\Omega_{i})$-subalgebra of $\mathrm{End}_{\sco(\Omega_{i})}\left(\sm(\Omega_{i})\right)$
generated by $\mathrm{im}\psi$, with structure map $\phi_{\Omega_{i}}:\mathbf{T}\to\mathbf{T}_{\Omega_{i}}$.
Let $\sx_{\Omega_{i}}$ be the affinoid rigid space $\mathrm{Sp}\mathbf{T}_{\Omega_{i}}$,
with $\pi:\sx_{\Omega_{i}}\to\Omega_{i}$ the natural morphism. The
canonical morphisms $\mathbf{T}_{\Omega_{i}}\otimes_{\sco(\Omega_{i})}\sco(\Omega_{ij})\to\mathbf{T}_{\Omega_{ij}}$
are isomorphisms, and so we may glue the affinoid rigid spaces $\sx_{\Omega_{i}}$
together via their overlaps $\sx_{\Omega_{ij}}$ into a rigid space
$\sx$ together with a finite map $\pi:\sx\to\sz$. The $\mathbf{T}_{\Omega_{i}}$-module
structure on $\sm(\Omega_{i})$ is compatible with the transition
maps, and so these modules glue to a coherent sheaf $\sm^{\dagger}$.
The structure maps $\phi_{\Omega_{i}}$ glue to a map $\phi:\mathbf{T}\to\sco(\sx)$
which is easily seen to be an algebra homomorphism. The remainder
of the theorem is straightforward from the construction. $\square$

The space $\sx$ is the \textbf{eigenvariety}\emph{ }associated with
the given eigenvariety datum. For any point $x\in\sx(\overline{\mathbf{Q}_{p}})$,
we write $\phi_{\sx,x}:\mathbf{T}\to k_{x}$ for the composite map\[
\left(\sco(\sx)\to\sco_{\sx,x}\twoheadrightarrow k_{x}\right)\circ\phi_{\sx},\]
and we say $\phi_{\sx,x}$ is the \textbf{eigenpacket }parametrized
by the point $x$. If the map $x\mapsto\phi_{\sx,x}$ determines a
bijection of $\sx(\overline{\mathbf{Q}_{p}})$ with a set of eigenpackets
of a certain type, we write $\phi\mapsto x_{\phi}$ for the inverse
map.

Maintaining the notation of the previous theorem, we note the following
useful tautology: if $\mathscr{Y}\hookrightarrow\scx$ is a closed
immersion with associated ideal sheaf $\mathscr{I}_{\mathscr{Y}}\subset\sco_{\scx}$,
then $\mathscr{Y}$ can be interpreted as the eigenvariety associated
with the eigenvariety datum\[
\mathfrak{D}_{\mathscr{Y}}=\left(\scw,\scz,\pi_{\ast}(\scm^{\dagger}\otimes_{\sco_{\scx}}\sco_{\scx}/\mathscr{I}_{\mathscr{Y}}),\mathbf{T},\phi_{\scx}\,\mathrm{mod}\,\mathscr{I}_{\mathscr{Y}}\right).\]

\subsection{Eigenvariety data from overconvergent cohomology}

Fix $\mathbf{G}$, $K^{p}$, a controlling operator $U_{t}$, and
an augmented Borel-Serre complex $C_{\bullet}(K^{p},-)$. For $\Omega\subset\sw_{K^{p}}$
an affinoid open, the Fredholm series $f_{\Omega}(X)=\det(1-\tilde{U}_{t}X)|C_{\bullet}(K^{p},\mathbf{A}_{\Omega}^{s})$
is well-defined independently of $s\geq s[\Omega]$ by Proposition
3.1.1, and if $\Omega'\subset\Omega$ is open then $f_{\Omega}(X)|_{\Omega'}=f_{\Omega'}(X)$.
By Tate's acyclicity theorem, there is a unique $f(X)\in\sco(\sw)\{\{X\}\}$
with $f(X)|_{\Omega}=f_{\Omega}(X)$ for all $\Omega$. Set $\sz=\sz_{f}$.
If $\sz_{\Omega,h}\subset\sz$ is a slope-adapted affinoid, then $C^{\bullet}(K^{p},\scd_{\Omega})$
admits a slope-$\leq h$ decomposition, with $C^{\bullet}(K^{p},\scd_{\Omega})_{\leq h}\simeq\mathrm{Hom}_{\sco(\Omega)}(C_{\bullet}(K^{p},\mathbf{A}_{\Omega}^{s})_{\leq h},\sco(\Omega))$
for any $s\geq s[\Omega]$, and $C^{\bullet}(K^{p},\scd_{\Omega})_{\leq h}$
is naturally a (graded) module over $\sco(\sz_{\Omega,h})\cong\sco(\Omega)[X]/(Q_{\Omega,h}(X))$
via the map $X\mapsto\tilde{U}_{t}^{-1}$; here $Q_{\Omega,h}(X)$
denotes the slope-$\leq h$ factor of $f_{\Omega}$.

\textbf{Proposition 4.3.1. }\emph{There is a unique complex $\mathscr{K}^{\bullet}$
of coherent analytic sheaves on $\sz$ such that $\mathscr{K}^{\bullet}(\sz_{\Omega,h})\cong C^{\bullet}(K^{p},\scd_{\Omega})_{\leq h}$
for any slope-adapted affinoid $\sz_{\Omega,h}$.}

\emph{Proof. }For $\sz_{\Omega,h}$ a slope-adapted affinoid, we simply
set\emph{ }$\mathscr{K}^{\bullet}(\sz_{\Omega,h})=C^{\bullet}(K^{p},\scd_{\Omega})_{\leq h}$,
with $\mathscr{K}^{\bullet}(\sz_{\Omega,h})$ regarded as an $\sco(\sz_{\Omega,h})$-module
as in the previous paragraph. $ $We are going to show that the formation
of $\mathscr{K}^{\bullet}(\sz_{\Omega,h})$ is compatible with overlaps
of slope-adapted affinoids; since slope-adapted affinoids form a base
for the ambient G-topology on $\sz$, this immediately implies that
the $\mathscr{K}^{\bullet}(\sz_{\Omega,h})$'s glue together into
a sheaf over $\sz$.

If $\sz_{\Omega,h}\in\mathscr{C}\mathrm{ov}$ and $\Omega'\subset\Omega$
with $\Omega'$ connected, a calculation gives $\sco(\sz_{\Omega',h})\simeq\sco(\sz_{\Omega,h})\otimes_{\sco(\Omega)}\sco(\Omega')$,
so then $\sz_{\Omega',h}\in\mathscr{C}\mathrm{ov}$. Fix $\sz_{\Omega',h'}\subseteq\sz_{\Omega,h}\in\mathscr{C}\mathrm{ov}$
with $\sz_{\Omega',h'}\in\mathscr{C}\mathrm{ov}$; we necessarily
have $\Omega'\subseteq\Omega$, and we may assume $h'\leq h$. Set
$C_{\Omega,h}=C^{\bullet}(K^{p},\scd_{\Omega})_{\leq h}$. We now
trace through the following sequence of canonical isomorphisms:\begin{eqnarray*}
C_{\Omega,h}\otimes_{\sco(\sz_{\Omega,h})}\sco(\sz_{\Omega',h'}) & \simeq & C_{\Omega,h}\otimes_{\sco(\sz_{\Omega,h})}\sco(\sz_{\Omega',h})\otimes_{\sco(\sz_{\Omega',h})}\sco(\sz_{\Omega',h'})\\
 & \simeq & \left(C_{\Omega,h}\otimes_{\sco(\sz_{\Omega,h})}\sco(\sz_{\Omega,h})\otimes_{\sco(\Omega)}\sco(\Omega')\right)\otimes_{\sco(\sz_{\Omega',h})}\sco(\sz_{\Omega',h'})\\
 & \simeq & \left(C_{\Omega,h}\otimes_{\sco(\Omega)}\sco(\Omega')\right)\otimes_{\sco(\sz_{\Omega',h})}\sco(\sz_{\Omega',h'})\\
 & \simeq & C_{\Omega',h}\otimes_{\sco(\sz_{\Omega',h})}\sco(\sz_{\Omega',h'})\\
 & \simeq & C_{\Omega',h'}.\end{eqnarray*}
The fourth line here follows from Proposition 3.1.5. $\square$

Taking the cohomology of $\mathscr{K}^{\bullet}$ yields a graded
sheaf $\sm^{\ast}$ on $\sz$ together with canonical isomorphisms
$\sm^{\ast}(\sz_{\Omega,h})\cong H^{\ast}(K^{p},\scd_{\Omega})_{\leq h}$
for any slope-adapted affinoid $\sz_{\Omega,h}$. By Proposition 3.1.5
the natural maps $\mathbf{T}(K^{p})\to\mathrm{End}_{\sco(\sz_{\Omega,h})}\left(H^{\ast}(K^{p},\scd_{\Omega})_{\leq h}\right)$
glue together into a degree-preserving algebra homomorphism $\psi:\mathbf{T}(K^{p})\to\mathrm{End}_{\sco_{\sz}}(\sm^{\ast})$.

\textbf{Definition 4.3.2. }\emph{The eigenvariety $\sx_{\mathbf{G},K^{p}}$
is the eigenvariety associated with the eigenvariety datum $(\sw_{K^{p}},\sz,\sm^{\ast},\mathbf{T}(K^{p}),\psi)$.
For $n$ a given integer, $\sx_{\mathbf{G},K^{p},n}$ is the eigenvariety
associated with the eigenvariety datum }$(\sw_{K^{p}},\sz,\sm^{n},\mathbf{T}(K^{p}),\psi)$.

Note that $\sz$ is highly noncanonical, depending as it does on a
choice of augmented Borel-Serre complex. However, $\sx_{\mathbf{G},K^{p}}$
is completely canonical and independent of this choice. Furthermore,
setting $\mathscr{J}=\mathrm{ann}_{\sco_{\sz}}\sm^{\ast}\subset\sco_{\sz}$,
the closed immersion $\scz_{t}\hookrightarrow\sw\times\mathbf{A}^{1}$
cut out by \[
\sco_{\sw\times\mathbf{A}^{1}}\twoheadrightarrow\sco_{\sz}\twoheadrightarrow\sco_{\sz}/\mathscr{J}=\sco_{\scz_{t}}\]
is independent of all choices, and is the {}``true'' spectral variety
for $U_{t}$. Note also that in practice, the eigenvarieties $\sx_{\mathbf{G},K^{p}}$
carry some extra structure which we don't really exploit in this article:
in particular, the sheaves $\sm^{\dagger}$ are sheaves of $\mathbf{T}_{\mathrm{ram}}(K^{p})$-modules.
Our first main result on the eigenvarieties $\sx_{\mathbf{G},K^{p}}$
is the following.

\textbf{Theorem 4.3.3. }\emph{The points $x\in\sx_{\mathbf{G},K^{p}}(\overline{\mathbf{Q}_{p}})$
lying over a given weight $\lambda\in\sw_{K^{p}}(\overline{\mathbf{Q}_{p}})$
are in bijection with the finite-slope eigenpackets for $\mathbf{G}$
of weight $\lambda$ and level $K^{p}$, and this bijection is realized
by sending $x\in\sx(\overline{\mathbf{Q}_{p}})$ to the eigenpacket
$\phi_{\sx,x}$.}

This theorem is due to Ash and Stevens (Theorem 6.2.1 of \cite{AS}),
but the following proof is new.

\emph{Proof. }Given a finite-slope eigenpacket $\phi$ of weight $\lambda$,
fix a slope-adapted affinoid $\sz_{\Omega,h}$ with $\lambda\in\Omega$
and $h>v_{p}(\phi(U_{t}))$, and let $\mathbf{T}_{\Omega,h}=\mathbf{T}_{\Omega,h}(K^{p})$
be the $\sco(\Omega)$-subalgebra of $\mathrm{End}_{\sco(\Omega)}\left(H^{\ast}(K^{p},\scd_{\Omega})_{\leq h}\right)$
generated by $\mathbf{T}(K^{p})\otimes_{\mathbf{Q}_{p}}\sco(\Omega)$.
Let $\mathfrak{M}$ be the maximal ideal of $\mathbf{T}(K^{p})\otimes_{\mathbf{Q}_{p}}\sco(\Omega)$
defined by \[
\mathfrak{M}=(T\otimes1+1\otimes x),\, T\in\ker\phi\,\mathrm{and}\, x\in\mathfrak{m}_{\lambda}.\]
After localizing the spectral sequence\[
E_{2}^{i,j}=\mathrm{Tor}_{-i}^{\sco(\Omega)}(H^{j}(K^{p},\scd_{\Omega})_{\leq h},k(\lambda))\Rightarrow H^{i+j}(K^{p},\scd_{\lambda})_{\leq h}\]
at $\mathfrak{M}$, the abutment is nonzero by assumption, so the
source must be nonzero as well. Therefore $\mathfrak{M}$ determines
a maximal ideal of $\mathbf{T}_{\Omega,h}$ lying over $\mathfrak{m}_{\lambda}$,
or equivalently a point $x\in\sx_{\Omega,h}$ with $w(x)=\lambda$.

On the other hand, given a point $x\in\sx_{\Omega,h}$ with $w(x)=\lambda$,
let $\mathfrak{M}=\mathfrak{M}_{x}\subset\mathbf{T}_{\Omega,h}$ be
the maximal ideal associated with $x$, and let $d$ be the largest
degree for which $H^{j}(K^{p},\scd_{\Omega})_{\leq h,\mathfrak{M}}\neq0$.
Localizing the spectral sequence at $\mathfrak{M}$, the entry $E_{2}^{0,d}$
is nonzero and stable, so the spectral sequence induces an isomorphism\[
0\neq H^{d}(K^{p},\scd_{\Omega})_{\leq h,\mathfrak{M}}\otimes_{\sco(\Omega)}k(\lambda)\simeq H^{d}(K^{p},\scd_{\lambda})_{\leq h,\mathfrak{M}},\]
and thus $\mathfrak{M}$ induces a finite-slope eigenpacket in weight
$\lambda$ as desired. $\square$

\subsection{The support of overconvergent cohomology modules}

As in the previous section, fix $\mathbf{G}$, $K^{p}$, and an augmented
Borel-Serre complex $C_{\bullet}(K^{p},-)$. We are going to prove
the following theorem.

\textbf{Theorem 4.4.1. }\emph{Fix a slope datum $(U_{t},\Omega,h)$.}
\begin{description}
\item [{i.}] \emph{For any $i$, $H_{i}(K^{p},\sca_{\Omega})_{\leq h}$
is a faithful $\sco(\Omega)$-module if and only if $H^{i}(K^{p},\scd_{\Omega})_{\leq h}$
is faithful.}
\item [{ii.}] \emph{If the derived group of $\mathbf{G}$ is $\mathbf{Q}$-anisotropic,
$H_{i}(K^{p},\sca_{\Omega})_{\leq h}$ and $H^{i}(K^{p},\scd_{\Omega})_{\leq h}$
are torsion $\sco(\Omega)$-modules for all $i$, unless $\mathbf{G}^{\mathrm{der}}(\mathbf{R})$
has a discrete series, in which case they are torsion for all $i\neq\frac{1}{2}\mathrm{dim}\mathbf{G}(\mathbf{R})/K_{\infty}$.}
\end{description}
Let $R$ be a Noetherian ring, and let $M$ be a finite $R$-module.
We say $M$ has \emph{full support }if $\mathrm{Supp}(M)=\mathrm{Spec}(R)$,
and that $M$ is \emph{torsion }if $\mathrm{ann}(M)\neq0$. We shall
repeatedly use the following basic result.

\textbf{Proposition 4.4.2. }\emph{If $\mathrm{Spec}(R)$ is reduced
and irreducible, the following are equivalent:}

\emph{i) $M$ is faithful (i.e. $\mathrm{ann}(M)=0$),}

\emph{ii) $M$ has full support,}

\emph{iii)} \emph{$M$ has nonempty open support,}

\emph{iv) $\mathrm{Hom}_{R}(M,R)\neq0$,}

\emph{v) $M\otimes_{R}K\ne0$, $K=\mathrm{Frac}(R)$.}

\emph{Proof. }Since $M$ is finite, $\mathrm{Supp}(M)$ is the underlying
topological space of $\mathrm{Spec}(R/\mathrm{ann}(M))$, so i) obviously
implies ii). If $\mathrm{Spec}(R/\mathrm{ann}(M))=\mathrm{Spec}(R)$
as topological spaces, then $\mathrm{ann}(M)\subset\surd(0)=(0)$
since $R$ is reduced, so ii) implies i). The set $\mathrm{Supp}(M)=\mathrm{Spec}(R/\mathrm{ann}(M))$
is \emph{a priori }closed; since $\mathrm{Spec}(R)$ is irreducible
by assumption, the only nonempty simultaneously open and closed subset
of $\mathrm{Spec}(R)$ is all of $\mathrm{Spec}(R)$, so ii) and iii)
are equivalent. By finiteness, $M$ has full support if and only if
$(0)$ is an associated prime of $M$, if and only if there is an
injection $R\hookrightarrow M$; tensoring with $K$ implies the equivalence
of ii) and v). Finally, $\mathrm{Hom}_{R}(M,R)\otimes_{R}K\simeq\mathrm{Hom}_{K}(M\otimes_{R}K,K)$,
so $M\otimes_{R}K\neq0$ if and only if $\mathrm{Hom}_{R}(M,R)\neq0$,
whence iv) and v) are equivalent. $\square$

\emph{Proof of Theorem 4.4.1.i}. (I'm grateful to Jack Thorne for
suggesting this proof.) Tensoring the Ext spectral sequence with $K(\Omega)=\mathrm{Frac}(\sco(\Omega))$,
it degenerates to isomorphisms \[
\mathrm{Hom}_{K(\Omega)}(H_{i}(K^{p},\sca_{\Omega})_{\leq h}\otimes_{\sco(\Omega)}K(\Omega),K(\Omega))\simeq H^{i}(K^{p},\scd_{\Omega})_{\leq h}\otimes_{\sco(\Omega)}K(\Omega),\]
so the claim is immediate from the preceding proposition. $\square$

\emph{Proof of Theorem 4.4.1.ii. }We give the proof in two steps,
with the first step naturally breaking into two cases. In the first
step, we prove the result assuming $\Omega$ contains an arithmetic
weight. In the second step, we eliminate this assumption via analytic
continuation.

\paragraph*{Step One, Case One: G doesn't have a discrete series.}

Let $\sw^{\mathrm{sd}}$ be the rigid Zariski closure in $\sw_{K^{p}}$
of the arithmetic weights whose algebraic parts are the highest weights
of irreducible $\mathbf{G}$-representations with nonvanishing $(\mathfrak{g},K_{\infty})$-cohomology.
A simple calculation using §II.6 of \cite{BWcohom} shows that $\sw^{\mathrm{sd}}$
is the union of its countable set of irreducible components, each
of dimension $<\mathrm{dim}\sw_{K^{p}}$. An arithmetic weight is
\emph{non-self-dual }if $\lambda\notin\sw^{\mathrm{sd}}$.

Now, by assumption $\Omega$ contains an arithmetic weight, so $\Omega$
automatically contains a Zariski dense $ $set $\mathcal{N}_{h}\subset\Omega\smallsetminus\Omega\cap\sw^{\mathrm{sd}}$
of non-self-dual arithmetic weights for which $h$ is a small slope.
By Theorem 3.2.5, $H^{\ast}(K^{p},\scd_{\lambda})_{\leq h}$ vanishes
identically for any $\lambda\in\mathcal{N}_{h}$. For any fixed $\lambda\in\mathcal{N}_{h}$,
suppose $\mathfrak{m}_{\lambda}\in\mathrm{Supp}_{\Omega}H^{\ast}(K^{p},\scd_{\Omega})_{\leq h}$;
let $d$ be the largest integer with $\mathfrak{m}_{\lambda}\in\mathrm{Supp}_{\Omega}H^{d}(K^{p},\scd_{\Omega})_{\leq h}$.
Taking $\Sigma=\lambda$ in the Tor spectral sequence gives\[
E_{2}^{i,j}=\mathrm{Tor}_{-i}^{\sco(\Omega)}(H^{j}(K^{p},\scd_{\Omega})_{\leq h},k(\lambda))\Rightarrow H^{i+j}(K^{p},\scd_{\lambda})_{\leq h}.\]
The entry $E_{2}^{0,d}=H^{d}(K^{p},\scd_{\Omega})_{\leq h}\otimes_{\sco(\Omega)}k(\lambda)$
is nonzero by Nakayama's lemma, and is stable since every row of the
$E_{2}$-page above the \emph{d}th row vanishes by assumption. In
particular, $E_{2}^{0,d}$ contributes a nonzero summand to the grading
on $H^{d}(K^{p},\scd_{\lambda})_{\leq h}$ - but this module is zero,
contradicting our assumption that $\mathfrak{m}_{\lambda}\in\mathrm{Supp}_{\Omega}H^{\ast}(K^{p},\scd_{\Omega})$.
Therefore, $H^{\ast}(K^{p},\scd_{\Omega})_{\leq h}$ does \emph{not
}have full support, so is not a faithful $\sco(\Omega)$-module.

\paragraph*{Step One, Case Two: G has a discrete series.}

The idea is the same as Case One, but with $\mathcal{N}_{h}$ replaced
by $\mathcal{R}_{h}$, the set of arithmetic weights with \emph{regular
}algebraic part for which $h$ is a small slope. For these weights,
Proposition 3.2.5 together with known results on $(\mathfrak{g},K_{\infty})$-cohomology
(see e.g. Sections 4-5 of \cite{LSEiscohom}) implies that $H^{i}(K^{p},\scd_{\lambda})_{\leq h},\,\lambda\in\mathcal{R}_{h}$\emph{
}vanishes for $i\neq d_{\mathbf{G}}=\frac{1}{2}\mathrm{dim}\mathbf{G}(\mathbf{R})/K_{\infty}$.
The Tor spectral sequence with $\Sigma=\lambda\in\mathcal{R}_{h}$
then shows that $\mathcal{R}_{h}$ doesn't meet $\mathrm{Supp}_{\Omega}H^{i}(K^{p},\scd_{\Omega})_{\leq h}$
for any $i>d_{\mathbf{G}}$. On the other hand, the Ext spectral sequence
with $\Sigma=\lambda\in\mathcal{R}_{h}$ then shows that $\mathcal{R}_{h}$
doesn't meet $\mathrm{Supp}_{\Omega}H_{i}(K^{p},\sca_{\Omega})_{\leq h}$
for any $i<d_{\mathbf{G}}$, whence the Ext spectral sequence with
$\Sigma=\Omega$ shows that $\mathcal{R}_{h}$ doesn't meet $\mathrm{Supp}_{\Omega}H^{i}(K^{p},\scd_{\Omega})_{\leq h}$
for any $i<d_{\mathbf{G}}$. The result follows.

\paragraph*{Step Two.}

We maintain the notation of §4.3. As in that subsection, $H^{n}(K^{p},\scd_{\Omega})_{\leq h}$
glues together over the slope-adapted affinoids $\sz_{\Omega,h}\subseteq\sz$
into a coherent $\sco_{\sz}$-module sheaf $\sm^{n}$, and in particular,
the support of $\sm^{n}$ is a closed analytic subset of $\sz$. Let
$w:\sz\to\sw$ denote the natural projection. For any $\sz_{\Omega,h}\in\mathscr{C}\mathrm{ov}$,
we have \[
w_{\ast}\mathrm{Supp}_{\sz_{\Omega,h}}\sm^{n}(\sz_{\Omega,h})=\mathrm{Supp}_{\Omega}H^{n}(K^{p},\scd_{\Omega})_{\leq h}.\]
Suppose $\mathrm{Supp}_{\Omega}H^{n}(K^{p},\scd_{\Omega})_{\leq h}=\Omega$
for some $\sz_{\Omega,h}\in\mathscr{C}\mathrm{ov}$. This implies
that $\mathrm{Supp}_{\sz_{\Omega,h}}\sm^{n}(\sz_{\Omega,h})$ contains
a closed subset of dimension equal to $\mathrm{dim}\sz$, so contains
an irreducible component of $\sz_{\Omega,h}$. Since $\mathrm{Supp}_{\sz}\sm^{n}$
is a priori closed, Corollary 2.2.6 of \cite{Conirred} implies that
$\mathrm{Supp}_{\sz}\sm^{n}$ contains an entire irreducible component
of $\sz$, say $\sz_{0}$. By Proposition 4.1.3, the image of $\sz_{0}$
is Zariski-open in $\sw$, so we may choose an arithmetic weight $\lambda_{0}\in w_{\ast}\sz_{0}$.
For some sufficiently large $h_{0}$ and some affinoid $\Omega_{0}$
containing $\lambda_{0}$, $\sz_{\Omega_{0},h_{0}}$ will contain
$\sz_{\Omega_{0},h_{0}}\cap\sz_{0}$ as a nonempty union of irreducible
components, and the latter intersection will be finite flat over $\Omega_{0}$.
Since $\sm^{n}(\sz_{\Omega_{0},h_{0}})\simeq H^{n}(K^{p},\scd_{\Omega_{0}})_{\leq h_{0}}$,
we deduce that $\mathrm{Supp}_{\Omega_{0}}H^{n}(K^{p},\scd_{\Omega_{0}})_{\leq h_{0}}=\Omega_{0}$,
whence $H^{n}(K^{p},\scd_{\Omega_{0}})_{\leq h_{0}}$ is faithful,
so by Step One $\mathbf{G}^{\mathrm{der}}(\mathbf{R})$ has a discrete
series and $n=\frac{1}{2}\mathrm{dim}\mathbf{G}(\mathbf{R})/K_{\infty}$.

\subsection{Eigenvarieties at noncritical interior points.}

In this section we prove the following result; part i. of this theorem
is a generalization of {}``Coleman families''.

\textbf{Theorem 4.5.1. }\emph{Let $x=x(\phi)\in\sx_{\mathbf{G},K^{p}}(\overline{\mathbf{Q}_{p}})$
be a point associated with a classical, noncritical, interior eigenpacket
$\phi$ such that $w(x)$ has regular algebraic part.}
\begin{description}
\item [{\emph{i.}}] \emph{If $l(\mathbf{G})=0$, every irreducible component
of $\sx$ containing $x$ has dimension equal to $\dim\sw_{K^{p}}$.}
\item [{\emph{ii.}}] \emph{If $l(\mathbf{G})\geq1$ and $\phi$ is strongly
interior, then every irreducible component of $\sx_{\mathbf{G},K^{p}}$
containing $x$ has dimension $\leq\mathrm{dim}\sw_{K^{p}}-1$, with
equality if $l(\mathbf{G})=1$.}
\end{description}
\emph{Proof. }By the basic properties of irreducible components together
with the construction given in §4.2-4.3, it suffices to work locally
over a fixed $\sz_{\Omega,h}$. Suppose $x\in\sx_{\Omega,h}$ is as
in the theorem, with $\phi:\mathbf{T}_{\Omega,h}\to\overline{\mathbf{Q}_{p}}$
the corresponding eigenpacket. Set $\mathfrak{M}=\ker\phi$, and let
$\mathfrak{m}=\mathfrak{m}_{\lambda}$ be the contraction of $\mathfrak{M}$
to $\sco(\Omega)$. Let $\mathfrak{P}\subset\mathbf{T}_{\Omega,h}$
be any minimal prime contained in $\mathfrak{M}$, and let $\mathfrak{p}$
be its contraction to a prime in $\sco(\Omega)$. The ring $\mathbf{T}_{\Omega,h}/\mathfrak{P}$
is a finite integral extension of $\sco(\Omega)/\mathfrak{p}$, so
both rings have the same dimension. Note also that $\mathfrak{p}$
is an associated prime of $H^{\ast}(K^{p},\scd_{\Omega})_{\leq h}$.

\textbf{Proposition 4.5.2. }\emph{The largest degrees for which $\phi$
occurs in $H^{\ast}(K^{p},\scd_{\Omega})_{\leq h}$ and $H^{\ast}(K^{p},\scd_{\lambda})_{\leq h}$
coincide, and the smallest degrees for which $\phi$ occurs in $H^{\ast}(K^{p},\scd_{\lambda})_{\leq h}$
and $H_{\ast}(K^{p},\sca_{\Omega})_{\leq h}$ coincide. Finally, the
smallest degree for which $\phi$ occurs in $H^{\ast}(K^{p},\scd_{\Omega})_{\leq h}$
is greater than or equal to the smallest degree for which $\phi$
occurs in $H_{\ast}(K^{p},\sca_{\Omega})_{\leq h}$.}

\emph{Proof. }For the first claim, localize the Tor spectral sequence
at $\mathfrak{M}$, with $\Sigma=\lambda$. If $\phi$ occurs in $H^{i}(K^{p},\scd_{\lambda})_{\leq h}$
then it occurs in a subquotient of $\mathrm{Tor}_{j}^{\sco(\Omega)}(H^{i+j}(K^{p},\scd_{\Omega})_{\leq h},k(\lambda))$
for some $j\geq0$. On the other hand, if $d$ is the largest degree
for which $\phi$ occurs in $H^{d}(K^{p},\scd_{\Omega})_{\leq h}$,
the entry $E_{2}^{0,d}$ of the Tor spectral sequence is stable and
nonzero after localizing at $\mathfrak{M}$, and it contributes to
the grading on $H^{d}(K^{p},\scd_{\lambda})_{\leq h,\mathfrak{M}}$.
The second and third claims follow from an analogous treatment of
the Ext spectral sequence. $\square$

First we treat the case where $l(\mathbf{G})=0$, so $\mathbf{G}^{\mathrm{der}}(\mathbf{R})$
has a discrete series. By the noncriticality of $\phi$ together with
the results recalled in §3.2, the only degree for which $\phi$ occurs
in $H^{i}(K^{p},\scd_{\lambda})_{\leq h}$ is the middle degree $d=\frac{1}{2}\mathrm{dim}\mathbf{G}(\mathbf{R})/K_{\infty}$,
so Proposition 4.5.2 implies that the only degree for which $\phi$
occurs in $H^{\ast}(K^{p},\scd_{\Omega})_{\leq h}$ is the middle
degree as well. The Tor spectral sequence localized at $\mathfrak{M}$
now degenerates, and yields\[
\mathrm{Tor}_{i}^{\sco(\Omega)_{\mathfrak{m}}}(H^{d}(K^{p},\scd_{\Omega})_{\leq h,\mathfrak{M}},k(\lambda))=0\,\mathrm{for\, all}\, i\geq1,\]
so $H^{d}(K^{p},\scd_{\Omega})_{\leq h,\mathfrak{M}}$ is a \emph{free
}module over $\sco(\Omega)_{\mathfrak{m}}$ by Proposition A.3. Since
$\sco(\Omega)_{\mathfrak{m}}$ is a domain and $\mathfrak{p}$ is
(locally at $\mathfrak{m}$) an associated prime of a free module,
$\mathfrak{p}=0$ and thus $\mathrm{dim}\mathbf{T}_{\Omega,h}/\mathfrak{P}=\mathrm{dim}\sco(\Omega)_{\mathfrak{m}}=\mathrm{dim}\sw$.

Now we turn to the case $l(\mathbf{G})\geq1$. First we demonstrate
the existence of an affinoid open $\mathscr{Y}\subset\mathscr{X}_{\Omega,h}$
containing $x$, and meeting every component of $\mathscr{X}_{\Omega,h}$
containing $x$, such that every regular classical non-critical point
in $\mathscr{Y}$ is cuspidal. By our assumptions, $\phi$ does not
occur in $H_{\partial}^{\ast}(K^{p},\scd_{\lambda})_{\leq h}$, so
by the boundary spectral sequence $\phi$ does not occur in $H_{\partial}^{\ast}(K^{p},\scd_{\Omega})_{\leq h}$.
Since $\mathrm{Supp}_{\mathbf{T}_{\Omega,h}}H_{\partial}^{\ast}(K^{p},\scd_{\Omega})_{\leq h}$
is closed in $\mathscr{X}_{\Omega,h}$ and does not meet $x$, the
existence of a suitable $\mathscr{Y}$ now follows easily. Shrinking
$\Omega$ and $\mathscr{Y}$ as necessary, we may assume that $\sco(\mathscr{Y})$
is finite over $\sco(\Omega)$, and thus $\mathscr{M}^{\ast}(\mathscr{Y})=H^{\ast}(K^{p},\scd_{\Omega})_{\leq h}\otimes_{\mathbf{T}_{\Omega,h}}\sco(\mathscr{Y})$
is finite over $\sco(\Omega)$ as well. Exactly as in the proof\emph{
}of Theorem 4.4.1, the Tor spectral sequence shows that $\mathrm{Supp}_{\Omega}\mathscr{M}^{\ast}(\mathscr{Y})$
doesn't contain any regular non-self-dual weights for which $h$ is
a small slope, so $\mathscr{M}^{\ast}(\mathscr{Y})$ and $\sco(\mathscr{Y})$
are torsion $\sco(\Omega)$-modules.

Finally, suppose $l(\mathbf{G})=1$. Set $d=q(\mathbf{G})$, so $\phi$
occurs in $H^{\ast}(K^{p},\scd_{\lambda})_{\mathfrak{M}}\simeq H^{\ast}(K^{p},\scl_{\lambda})_{\mathfrak{M}}$
only in degrees $d$ and $d+1$. By the argument of the previous paragraph,
$H^{\ast}(K^{p},\scd_{\Omega})_{\leq h,\mathfrak{M}}$ is a torsion
$\sco(\Omega)_{\mathfrak{m}}$-module. Taking $\Sigma=\lambda$ in
the Ext spectral sequence and localizing at $\mathfrak{M}$, Proposition
4.5.2 yields\[
H^{d}(K^{p},\scd_{\Omega})_{\leq h,\mathfrak{M}}\simeq\mathrm{Hom}_{\sco(\Omega)_{\mathfrak{m}}}(H_{d}(K^{p},\sca_{\Omega})_{\leq h,\mathfrak{M}},\sco(\Omega)_{\mathfrak{m}}).\]
Since the left-hand term is a torsion $\sco(\Omega)_{\mathfrak{m}}$-module,
Proposition 4.4.2 implies that both modules vanish identically. Proposition
4.5.2 now shows that $d+1$ is the only degree for which $\phi$ occurs
in $H^{\ast}(K^{p},\scd_{\Omega})_{\leq h}$. Taking $\Sigma=\lambda$
in the Tor spectral sequence and localizing at $\mathfrak{M}$, the
only nonvanishing entries are $E_{2}^{0,d+1}$ and $E_{2}^{-1,d+1}$.
In particular, $\mathrm{Tor}_{i}^{\sco(\Omega)_{\mathfrak{m}}}(H^{d+1}(K^{p},\scd_{\Omega})_{\leq h,\mathfrak{M}},\sco(\Omega)/\mathfrak{m})=0$
for all $i\geq2$, so $H^{d+1}(K^{p},\scd_{\Omega})_{\leq h,\mathfrak{M}}$
has projective dimension at most one by Proposition A.3. Summarizing,
we've shown that $H^{i}(K^{p},\scd_{\Omega})_{\leq h,\mathfrak{M}}$
vanishes in degrees $\neq q(\mathbf{G})+1$, and that $H^{q(\mathbf{G})+1}(K^{p},\mathcal{D}_{\Omega})_{\le h,\mathfrak{M}}$
is a torsion $\sco(\Omega)_{\mathfrak{m}}$-module of projective dimension
one, so $\mathrm{ht}\mathfrak{p}=1$ by Proposition A.6. $\square$

We remark here that results towards Newton's Theorem 1.1.6 were established by Stevens and Urban in unpublished
work: however, in the notation of the present subsection, their proofs required an a priori assumption that the
ideal $\mathfrak{p}$ is generated by a regular sequence locally in
$\sco(\Omega)_{\mathfrak{m}}$, which seems quite hard to check.

\subsection{General linear groups}

In this section we examine the special case when $\mathbf{G}\simeq\mathrm{Res}_{F/\mathbf{Q}}\mathbf{H}$
for some number field $F$ and some $F$-inner form $\mathbf{H}$
of $\mathrm{GL}_{n}/F$, introducing notation which will remain in
effect throughout Chapters 5 and 6. In particular, we often work with
a canonical family of tame level subgroups suggested by the theory
of new vectors. More precisely, given an integral ideal $\mathfrak{n}=\prod\mathfrak{p}_{v}^{e_{v}(\mathfrak{n})}\subset\mathcal{O}_{F}$
with $e_{v}=0$ if $v|p$ or if $\mathbf{H}(F_{v})\ncong\mathrm{GL}_{n}(F_{v})$,
set \[
K_{1}(\mathfrak{n})=\prod_{v\,\mathrm{with}\, e_{v}(\mathfrak{n})>0}K_{v}(\varpi_{v}^{e_{v}(\mathfrak{n})})\prod_{v\,\mathrm{with}\, e_{v}(\mathfrak{n})=0}K_{v}\]
where $K_{v}(\varpi_{v}^{e})$ denotes the open compact subgroup of
$\mathrm{GL}_{n}(\mathcal{O}_{v})$ consisting of matrices with lowest
row congruent to $(0,\dots,0,1)\,\mathrm{mod}\,\varpi_{v}^{e}$, and
$K_{v}$ denotes a fixed maximal compact subgroup of $\mathbf{H}(F_{v})$.
The Hecke algebra $\mathbf{T}(K_{1}(\mathfrak{n}))$ then contains
the usual operators $T_{v,i}$ corresponding to the double cosets
of the matrices $\mathrm{diag}\left(\underbrace{\varpi_{v},\dots,\varpi_{v}}_{i},1,\dots,1\right)$
for $1\leq i\leq n$ and $v$ a place of $F$ such that $e_{v}(\mathfrak{n})=0$
and $\mathbf{H}(F_{v})\simeq\mathrm{GL}_{n}(F_{v})$. For a place
$v|p$ we write $U_{v,i}$ for the element of $\mathcal{A}_{p}^{+}$
corresponding to the double coset of $\mathrm{diag}\left(1,\dots,1,\underbrace{\varpi_{v},\dots,\varpi_{v}}_{i}\right)$,
and we set $u_{v,i}=U_{v,i-1}^{-1}U_{v,i}\in\mathcal{A}_{p}$. The
operator $U_{p}=\prod_{v|p}\prod_{i=1}^{n-1}U_{v,i}$ is a canonical
controlling operator. If $n=2$ we adopt the more classical notation,
writing $T_{v}=T_{v,1}$ and $S_{v}=T_{v,2}$. We write $\mathbf{T}_{\lambda}(\mathfrak{n})$
for the finite-slope Hecke algebra of weight $\lambda$ and tame level
$K_{1}(\mathfrak{n})$ as in §3.2, and we abbreviate $\scx_{\mathrm{Res}_{F/\mathbf{Q}}\mathrm{GL}_{n},K_{1}(\mathfrak{n})}$
by $\scx_{\mathrm{GL}_{n}/F,\mathfrak{n}}$.

If $F=\mathbf{Q}$, we define $\sw_{K^{p}}^{0}$ as the subspace of
$\sw_{K^{p}}$ parametrizing characters trivial on the one-parameter
subgroup $\mathrm{diag}(1,\dots,1,t_{n})$, and we set $\sx_{\mathbf{G},K^{p}}^{0}=\sx_{\mathbf{G},K^{p}}\cap w^{-1}\left(\sw_{K^{p}}^{0}\right)$.
By the remarks in §2.2, $\sx_{\mathbf{G},K^{p}}^{0}$ is a disc bundle
over $\sx_{\mathbf{G},K^{p}}$: for any point $x\in\sx_{\mathbf{G},K^{p}}$
with $w(x)=(\lambda_{1},\dots,\lambda_{n})\in\scw_{K^{p}}$, there
is a unique point $x^{0}\in\sx_{\mathbf{G},K^{p}}^{0}$ with $\lambda^{0}=(\lambda_{1}\lambda_{n}^{-1},\dots,\lambda_{n-1}\lambda_{n}^{-1},1)$
such that $\phi(x)(T_{\ell,i})=\lambda_{n}(\ell)^{i}\phi(x^{0})(T_{\ell,i})$.
Restricting attention to $\sx_{\mathbf{G},K^{p}}^{0}$ amounts to
factoring out {}``wild twists'', and $\sx_{\mathrm{GL}_{2}/\mathbf{Q}_{p},K^{p}}^{0}$
is canonically isomorphic to the Coleman-Mazur-Buzzard eigencurve
of tame level $K^{p}$ (this is a theorem of Bellaïche \cite{BeCrit}).

\section{\emph{p-}adic Langlands functoriality}

\subsection{An interpolation theorem}

\textbf{Definition 5.1.1. }\emph{Given an eigenvariety datum $\mathfrak{D}=(\sw,\sz,\sm,\mathbf{T},\psi)$
with associated eigenvariety $\sx$, the }core \emph{of $\sx$, denoted
$\sx^{\circ}$, is the union of the $\dim\sw$-dimensional irreducible
components of the nilreduction $\sx^{\mathrm{red}}$, regarded as
a closed subspace of $\sx$. An eigenvariety $\sx$ is }unmixed \emph{if
$\sx^{\circ}\simeq\sx$.}

Let $\sz^{\circ}$ denote the subspace of points in $\sz$ whose preimage
in $\sx$ meets the core of $\sx$, with its reduced rigid subspace
structure; $\sz^{\circ}$ is naturally a union of irreducible components
of $\sz^{\mathrm{red}}$. We will see below that $\sx^{\circ}$ really
is an eigenvariety, in the sense of being associated with an eigenvariety
datum.

Suppose we are given two eigenvariety data $\mathfrak{D}_{i}=(\sw_{i},\sz_{i},\sm_{i},\mathbf{T},\psi_{i})$
for $i=1,2$, together with a closed immersion $\jmath:\sw_{1}\hookrightarrow\sw_{2}$;
we write $j$ for the natural extension of $\jmath$ to a closed immersion
$\jmath\times\mathrm{id}:\sw_{1}\times\mathbf{A}^{1}\hookrightarrow\sw_{2}\times\mathbf{A}^{1}$.
Given a point $z\in\sz_{i}$ and any $T\in\mathbf{T}$, we write $D_{i}(T,X)(z)\in k(z)[X]$
for the characteristic polynomial $\det(1-\psi_{i}(T)X)|\sm_{i}(z)$.

\textbf{Theorem 5.1.2. }\emph{Notation and assumptions as in the previous
paragraph, suppose there is some very Zariski-dense set $\sz^{\mathrm{cl}}\subset\sz_{1}^{\circ}$
with $j(\sz^{\mathrm{cl}})\subset\sz_{2}$ such that the polynomial
$D_{1}(T,X)(z)$ divides $D_{2}(T,X)(j(z))$ in $k(z)[X]$ for all
$T\in\mathbf{T}$ and all $z\in\sz^{\mathrm{cl}}$. Then $j$ induces
a closed immersion $\zeta:\sz_{1}^{\circ}\hookrightarrow\sz_{2}$,
and there is a canonical closed immersion $i:\sx_{1}^{\circ}\hookrightarrow\sx_{2}$
such that the diagrams\[
\xymatrix{\sx_{1}^{\circ}\ar[d]_{w_{1}}\ar@{^{(}->}[r]^{i} & \sx_{2}\ar[d]^{w_{2}}\\
\sw_{1}\ar@{^{(}->}[r]_{\jmath} & \sw_{2}}
\]
and\[
\xymatrix{\mathbf{T}\ar[rr]^{\phi_{1}^{\circ}}\ar[dr]_{\phi_{2}} &  & \sco(\sx_{1}^{\circ})\\
 & \sco(\sx_{2})\ar[ur]_{i^{\ast}}}
\]
commute.}

Before proving this result, we establish three lemmas.

\textbf{Lemma 5.1.3. }\emph{Suppose $A$ is an affinoid algebra, $B$
is a module-finite $A$-algebra, and $S$ is a Zariski-dense subset
of $\mathrm{Max}A$. Then\[
I=\bigcap_{\mathfrak{m}\in\mathrm{Max}B\,\mathrm{with}\,\mathfrak{m}\,\mathrm{lying}\,\mathrm{over}\,\mathrm{some}\,\mathfrak{n}\in S}\mathfrak{m}\subset B\]
is contained in every minimal prime $\mathfrak{p}$ of $B$ with $\mathrm{dim}B/\mathfrak{p}=\mathrm{dim}A$.}

\emph{Proof. }Translated into geometric language, this is the self-evident
statement that the preimage of $S$ under $\pi:\mathrm{Sp}B\to\mathrm{Sp}A$
is Zariski-dense in any irreducible component of $\mathrm{Sp}B$ with
$\mathrm{dim}A$-dimensional image in $\mathrm{Sp}A$. $\square$

\textbf{Lemma 5.1.4. }\emph{Suppose $A\to B$ are affinoid algebra,
with $\mathrm{Max}B$ an affinoid subdomain of $\mathrm{Max}A$. Let
$A^{\circ}$ be the maximal reduced quotient of $A$ which is equidimensional
of dimension $\mathrm{dim}A$. Then $A^{\circ}\otimes_{A}B$ is the
maximal reduced quotient of $B$ which is equidimensional of dimension
$\mathrm{dim}A$ if $\dim B=\dim A$, and is zero if $\mathrm{dim}B<\dim A$.}

\emph{Proof. }Set $d=\mathrm{dim}A$. The kernel of $A\to A^{\circ}$
is the ideal $I^{\circ}=\cap_{\mathfrak{p}\in\mathrm{Spec}A,\mathrm{coht}\mathfrak{p}=d}\mathfrak{p}$,
so tensoring the sequence \[
0\to I^{\circ}\to A\to A^{\circ}\to0\]
with $B$ yields\[
0\to I^{\circ}\otimes_{A}B\to B\to A^{\circ}\otimes_{A}B\to0\]
by the $A$-flatness of $B$. It thus suffices to prove an isomorphism
\[
I^{\circ}\otimes_{A}B\simeq\bigcap_{\mathfrak{p}\in\mathrm{Spec}B,\,\mathrm{coht}\mathfrak{p}=d}\mathfrak{p}.\]
By the Jacobson property of affinoid algebras, we can rewrite $I^{\circ}$
as follows:\begin{eqnarray*}
I^{\circ} & = & \bigcap_{\mathrm{coht}\mathfrak{p}=d}\mathfrak{p}\\
 & = & \bigcap_{\substack{\mathfrak{m}\in\mathrm{Max}(A)\,\mathrm{with}\\
\mathfrak{m}\supset\mathfrak{p}\,\mathrm{and}\,\mathrm{coht}\mathfrak{p}=d}
}\mathfrak{m}\\
 & = & \bigcap_{\substack{\mathfrak{m}\in\mathrm{Max}(A)\,\mathrm{with}\\
\mathrm{ht}\mathfrak{m}=d}
}\mathfrak{m}.\end{eqnarray*}
Since $B$ is $A$-flat, we have $(I_{1}\cap I_{2})\otimes_{A}B=I_{1}B\cap I_{2}B$
for any ideals $I_{i}\subset A$, so\begin{eqnarray*}
I^{\circ}\otimes_{A}B & = & \bigcap_{\substack{\mathfrak{m}\in\mathrm{Max}(A)\,\mathrm{with}\\
\mathrm{ht}\mathfrak{m}=d}
}\mathfrak{m}B\\
 & = & \bigcap_{\substack{\mathfrak{m}\in\mathrm{Max}(B)\,\mathrm{with}\\
\mathrm{ht}\mathfrak{m}=d}
}\mathfrak{m}\\
 & = & \bigcap_{\mathfrak{p}\in\mathrm{Spec}B,\,\mathrm{coht}\mathfrak{p}=d}\mathfrak{p},\end{eqnarray*}
and we're done if $\dim B=\mathrm{dim}A$. But if $\mathrm{dim}B<\mathrm{dim}A$,
then $\mathfrak{m}B=B$ for all $\mathfrak{m}\in\mathrm{Max}(A)$
of height $d$, so $I^{\circ}\otimes_{A}B=B$ as desired. $\square$

\textbf{Lemma 5.1.5. }\emph{Let $\mathbf{T}$ be a commutative $\mathbf{Q}_{p}$-algebra,
$L/\mathbf{Q}_{p}$ a finite extension, and $M_{1},M_{2}$ a pair
of $\mathbf{T}\otimes_{\mathbf{Q}_{p}}L$-modules finite over $L$.
For any $t\in\mathbf{T}$ set $D_{i}(X,t)=\det(1-Xt)|M_{i}\in L[X]$.
If $D_{1}(X,t)|D_{2}(X,t)$ in $L[X]$ for any $t\in\mathbf{T}$,
then the same divisibility holds for any $t\in\mathbf{T}\otimes_{\mathbf{Q}_{p}}L$.}

Proof. Given an arbitrary element $t=\sum_{i=1}^{n}t_{i}\otimes a_{i}\in\mathbf{T}\otimes_{\mathbf{Q}_{p}}L$
with $t_{i}\in\mathbf{T}$ and $a_{i}\in L$, we set\begin{eqnarray*}
P_{i}(X,X_{1},\dots,X_{n}) & = & \det\left(1-X(X_{1}t_{1}+X_{2}t_{2}+\cdots+X_{n}t_{n})\right)|M_{i}\\
 & \in & L[X,X_{1},\dots,X_{n}],\end{eqnarray*}
so in particular $P_{i}(X,a_{1},\dots,a_{n})=D_{i}(X,t)$. Consider
the meromorphic quotient \begin{eqnarray*}
Q(X,X_{1},\dots,X_{n}) & = & \frac{P_{2}(X,X_{1},\dots,X_{n})}{P_{1}(X,X_{1},\dots,X_{n})}.\end{eqnarray*}
Let $U\subset\mathrm{Spec}L[X_{1},\dots,X_{n}]^{\mathrm{an}}$ be
a sufficiently small connected affinoid open neighborhood of the origin
and let $C\gg0$ be an integer such that on the domain $\{|X|\leq p^{-C}\}\times U$,
$Q$ admits a convergent power series expansion\[
Q(X,X_{1},\dots,X_{n})=\sum_{j=0}^{\infty}f_{j}(X_{1},\dots,X_{n})X^{j}\in\sco(U)\left\langle p^{-C}X\right\rangle .\]
By hypotheis, for any $(b_{1},\dots,b_{n})\in\mathbf{Q}_{p}^{n}$
the specialization $Q(X,b_{1},\dots,b_{n})$ is an element of $L[X]$,
necessarily of degree at most $d=\mathrm{dim}_{L}M_{2}$. In particular,
for any $j>d$ the function $f_{j}$ vanishes on the Zariski-dense
set $U\cap\mathbf{Q}_{p}^{n}$ and therefore vanishes identically.
We now have the identity $Q(X,X_{1},\dots,X_{n})=\sum_{0\leq j\leq d}f_{j}(X_{1},\dots,X_{n})X^{j}\in\sco(U)[X]$.
For some $N\gg0$ we have $(p^{N}a_{1},\dots,p^{N}a_{n})\in U$; specializing
$Q$ at this point and making the change of variables $X\mapsto p^{-N}X$
yields\begin{eqnarray*}
\frac{D_{2}(X,t)}{D_{1}(X,t)} & = & Q(X,a_{1},\dots,a_{n})\\
 & = & Q(p^{-N}X,p^{N}a_{1},\dots,p^{N}a_{n})\\
 & = & \sum_{0\leq j\leq d}f_{j}(p^{N}a_{1},\dots,p^{N}a_{n})p^{-Nj}X^{j}\\
 & \in & L[X].\end{eqnarray*}
$\square$

\emph{Proof of Theorem 5.1.2.} Set $d=\mathrm{dim}\sw_{1}$. First,
we establish the theorem in the special case when $\sw_{1}=\sw_{2}$,
$\jmath=\mathrm{id}$, $\sz_{1}^{\circ}\simeq\sz_{1}\simeq\sz_{2}$,
and $\sx_{1}^{\circ}\simeq\sx_{1}$; we refer to this as the \emph{narrow
case.} For brevity we write $\sw=\sw_{1}$ and $\sz=\sz_{1}$. As
in §4.2, let $\scc\mathrm{ov}=\left\{ \Omega_{i}\right\} _{i\in I}$
be an admissible affinoid covering of $\sz$. For any $\Omega\in\scc\mathrm{ov}$
and $i\in\{1,2\}$, let $\mathbf{T}_{\Omega,i}$ denote the $\sco(\Omega)$-subalgebra
of $\mathrm{End}_{\sco(\Omega)}\left(\sm_{i}(\Omega)\right)$ generated
by the image of $\psi_{i}$, and let $I_{\Omega,i}$ be the kernel
of the natural surjection \[
\phi_{\Omega,i}:\mathbf{T}\otimes_{\mathbf{Q}_{p}}\sco(\Omega)\twoheadrightarrow\mathbf{T}_{\Omega,i}.\]
We are going to establish an inclusion $I_{\Omega,2}\subseteq I_{\Omega,1}$
for all $\Omega\in\scc\mathrm{ov}$. Granting this inclusion, let
$\mathscr{I}_{\Omega}\subset\mathbf{T}_{\Omega,2}$ be the kernel
of the induced surjection $\mathbf{T}_{\Omega,2}\twoheadrightarrow\mathbf{T}_{\Omega,1}$.
If $\Omega'\subset\Omega$ is an affinoid subdomain, then applying
$-\otimes_{\sco(\Omega)}\sco(\Omega')$ to the sequence\[
0\to I_{\Omega,i}\to\mathbf{T}\otimes_{\mathbf{Q}_{p}}\sco(\Omega)\overset{\phi_{\Omega,i}}{\to}\mathbf{T}_{\Omega,i}\to0\]
yields a canonical isomorphism $I_{\Omega,i}\otimes_{\sco(\Omega)}\sco(\Omega')\cong I_{\Omega',i}$,
so applying $-\otimes_{\sco(\Omega)}\sco(\Omega')$ to the canonical
isomorphism $\mathscr{I}_{\Omega}\cong I_{\Omega,2}/I_{\Omega,1}$
yields an isomorphism $\mathscr{I}_{\Omega}\otimes_{\sco(\Omega)}\sco(\Omega')\cong\mathscr{I}_{\Omega'}$.
Therefore the assignments $\Omega\mapsto\mathscr{I}_{\Omega}$ glue
together into a coherent ideal subsheaf of the structure sheaf of
$ $$\sx_{2}$ cutting out $\sx_{1}$; equivalently, the surjections
$\mathbf{T}_{\Omega,2}\twoheadrightarrow\mathbf{T}_{\Omega,1}$ glue
together over $\Omega\in\scc\mathrm{ov}$ into the desired closed
immersion.

It remains to establish the inclusion $I_{\Omega,2}\subseteq I_{\Omega,1}$.
Let $\sz^{\mathrm{reg}}$ be the maximal subset of $\sz$ such that
$\sco_{\sz,z}$ is regular for all $z\in\sz^{\mathrm{reg}}$ and the
sheaves $\sm_{1}$ and $\sm_{2}$ are locally free after restriction
to $\sz^{\mathrm{reg}}$; since $\sz^{\mathrm{reg}}$ is naturally
the intersection \[
\mathrm{Reg}(\sz)\bigcap\left(\sz\smallsetminus\mathrm{Supp}\oplus_{i=1}^{d}\se\mathrm{xt}_{\sco_{\sz}}^{i}(\sm_{1}\oplus\sm_{2},\sco_{\sz})\right),\]
and $\mathrm{Reg}(\sz)$ is Zariski-open by the excellence of affinoid
algebras, $\sz^{\mathrm{reg}}$ is naturally a Zariski-open and Zariski-dense
rigid subspace of $\sz$. For any $T\in\mathbf{T}$, let $D_{i}(T,X)\in\sco(\sz^{\mathrm{reg}})[X]$
be the determinant of $1-\psi_{i}(T)X$ acting on $\sm_{i}|_{\sz^{\mathrm{reg}}}$,
defined in the usual way (this is why we need local freeness). For
any $z\in\sz^{\mathrm{reg}}$, the image of $D_{i}(T,X)$ in the residue
ring $k(z)[X]$ is simply $D_{i}(T,X)(z)$. By our hypotheses, the
formal quotient \[
Q(T,X)=D_{2}(T,X)/D_{1}(T,X)=\sum_{n\geq0}a_{n}X^{n}\in\sco(\sz^{\mathrm{reg}})[[X]]\]
reduces for any $z\in\sz^{\mathrm{reg}}\cap\sz^{\mathrm{cl}}$ to
an element of $k(z)[X]$ with degree bounded uniformly above as a
function of $z$ on any given irreducible component of $\sz^{\mathrm{reg}}$.
In particular, the restriction of $a_{n}$ to any given irreducible
component of $\sz^{\mathrm{reg}}$ is contained in a Zariski-dense
set of maximal ideals for sufficiently large $n$, and so is zero.
Thus $D_{1}(T,X)(z)$ divides $D_{2}(T,X)(z)$ in $k(z)[X]$ for \emph{any
}$z\in\sz^{\mathrm{reg}}$ and any $T\in\mathbf{T}$. This extends
by Lemma 5.1.5 to the same divisibility for any $T\in\mathbf{T}\otimes_{\mathbf{Q}_{p}}k(z)$.

Suppose now that $T\in\mathbf{T}\otimes_{\mathbf{Q}_{p}}\sco(\Omega)$
is contained in $I_{\Omega,2}$. Since $D_{2}(T,X)(z)=1$ for any
$z\in\Omega\cap\sz^{\mathrm{reg}}$, the deduction in the previous
paragraph shows that $D_{1}(T,X)(z)=1$ for any $z\in\Omega\cap\sz^{\mathrm{reg}}$.
But then \[
\phi_{1}(T)\in\bigcap_{x\in\mathrm{Sp}\mathbf{T}_{\Omega,1}\,\mathrm{with}\,\pi(x)\in\Omega\cap\sz^{\mathrm{reg}}}\mathfrak{m}_{x}\subseteq\bigcap_{\mathfrak{p}\in\mathrm{Spec}\mathbf{T}_{\Omega,1},\,\mathrm{coht}\mathfrak{p}=d}\mathfrak{p}=0,\]
where the middle inclusion follows from Lemma 5.1.3 and the rightmost
equality follows since $\mathbf{T}_{\Omega,1}$ is reduced and equidimensional
of dimension $d$ by assumption. This establishes the narrow case.

It remains to establish the general case. By the hypotheses of the
theorem, $j(\sz^{\mathrm{cl}})\in j(\sz_{1}^{\circ})\cap\sz_{2}$,
so $j$ induces the closed immersion $\zeta:\sz_{1}^{\circ}\hookrightarrow\sz_{2}$
by the Zariski-density of $\sz^{\mathrm{cl}}$ in $\sz_{1}^{\circ}$.
Let $\sx_{2}^{'}$ denote the fiber product $\sx_{2}\times_{\sz_{2},\zeta}\sz_{1}^{\circ}$;
note that $\sx_{2}'$ is the eigenvariety associated with the eigenvariety
datum $\mathfrak{D}_{2}'=(\sw_{1},\sz_{1}^{\circ},\zeta^{\ast}\sm_{2},\mathbf{T},\zeta^{\sharp}\psi_{2})$,
and there is a canonical closed immersion $i':\sx_{2}'\hookrightarrow\sx_{2}$
by construction. For $\Omega\subseteq\sz_{1}$ an affinoid open we
define an ideal $\sj(\Omega)\subset\mathbf{T}_{\Omega,1}$ by the
rule\[
\sj(\Omega)=\begin{cases}
\mathbf{T}_{\Omega,1} & \mathrm{if}\,\mathrm{dim}\mathbf{T}_{\Omega,1}<d\\
\bigcap_{x\in\mathrm{Sp}\mathbf{T}_{\Omega,1},\,\mathrm{ht}\mathfrak{m}_{x}=d}\mathfrak{m}_{x} & \mathrm{if}\,\mathrm{dim}\mathbf{T}_{\Omega,1}=d.\end{cases}\]
By Lemma 5.1.4 the ideals $\sj(\Omega)$ glue into a coherent ideal
sheaf $\sj\subset\sco_{\sx_{1}}$. The support of $\sco_{\sx_{1}}/\sj$
in $\sz_{1}$ is exactly $\sz_{1}^{\circ}$ , and in fact the closed
immersion cut out by$\sj$ is exactly the core of $\sx_{1}$. In particular,
the core of $\sx_{1}$ is the eigenvariety associated with the (somewhat
tautological) eigenvariety datum $\mathfrak{D}_{1}^{\circ}=\left(\sw_{1},\sz_{1}^{\circ},\pi_{\ast}\left(\sm_{1}^{\dagger}\otimes_{\sco_{\sx_{1}}}\sco_{\sx_{1}}/\sj\right),\mathbf{T},\psi\,\mathrm{mod}\,\sj\right)$.
The narrow case of Theorem 5.1.2 applies to the pair of eigenvariety
data $\mathfrak{D}_{1}^{\circ}$ and $\mathfrak{D}_{2}'$, producing
a closed immersion $i'':\sx_{1}^{\circ}\hookrightarrow\sx_{2}'$,
and the general case follows upon setting $i=i'\circ i''$. $ $$\square$

From Theorem 5.1.2, it's easy to deduce the following more flexible
interpolation theorem.

\textbf{Theorem 5.1.6. }\emph{Suppose we are given two eigenvariety
data $\mathfrak{D}_{i}=(\sw_{i},\sz_{i},\sm_{i},\mathbf{T}_{i},\psi_{i})$
for $i=1,2$, together with the following additional data:}

\emph{i) A closed immersion $\jmath:\sw_{1}\hookrightarrow\sw_{2}$.}

\emph{ii) An algebra homomorphism $\sigma:\mathbf{T}_{2}\to\mathbf{T}_{1}$.}

\emph{iii) A very Zariski-dense set $\sz^{\mathrm{cl}}\subset\sz_{1}^{\circ}$
with $j(\sz^{\mathrm{cl}})\subset\sz_{2}$ such that $D_{1}(\sigma(T),X)(z)$
divides $D_{2}(T,X)(j(z))$ in $k(z)[X]$ for all $z\in\sz^{\mathrm{cl}}$
and all $T\in\mathbf{T}_{2}$.}\\
\emph{Then there exists a morphism $i:\sx_{1}^{\circ}\to\sx_{2}$
such that the diagrams\[
\xymatrix{\sx_{1}^{\circ}\ar[d]_{w_{1}}\ar[r]^{i} & \sx_{2}\ar[d]^{w_{2}}\\
\sw_{1}\ar@{^{(}->}[r]_{\jmath} & \sw_{2}}
\]
and\[
\xymatrix{\sco(\sx_{2})\ar[r]^{i^{\ast}} & \sco(\sx_{1}^{\circ})\\
\mathbf{T}_{2}\ar[r]_{\sigma}\ar[u]^{\phi_{2}} & \mathbf{T}_{1}\ar[u]_{\phi_{1}^{\circ}}}
\]
commute, and $i$ may be written as a composite $i_{c}\circ i_{f}$
where $i_{f}$ is a finite morphism and $i_{c}$ is a closed immersion.}

\emph{Proof. }Let $\mathfrak{D}_{1}^{\sigma}$ be the eigenvariety
datum $(\sw_{1},\sz_{1},\sm_{1},\mathbf{T}_{2},\psi_{1}\circ\sigma)$.
Theorem 5.1.2 produces a closed immersion $i_{c}':\sx(\mathfrak{D}_{1}^{\sigma})^{\circ}\hookrightarrow\sx(\mathfrak{D}_{2})$.
The inclusion $\mathrm{im}(\psi_{1}\circ\sigma)(\mathbf{T}_{2})\subset\mathrm{im}\psi_{1}(\mathbf{T}_{1})\subset\mathrm{End}_{\sco(\Omega)}(\sm_{1}(\Omega))$
induces a finite morphism $i_{f}'=\sx(\mathfrak{D}_{1})\to\sx(\mathfrak{D}_{1}^{\sigma})$.
The ideal subsheaf of $\sco_{\sx(\mathfrak{D}_{1}^{\sigma})}$ cut
out by the kernel of the composite $\sco_{\sx(\mathfrak{D}_{1}^{\sigma})}\to\sco_{\sx(\mathfrak{D}_{1})}\twoheadrightarrow\sco_{\sx(\mathfrak{D}_{1})^{\circ}}$
determines a closed immersion $\mathscr{Y}\hookrightarrow\sx(\mathfrak{D}_{1}^{\sigma})$
fitting into a diagram\[
\xymatrix{\sx(\mathfrak{D}_{1})^{\circ}\ar[d]_{i_{f}}\ar@{^{(}->}[r] & \sx(\mathfrak{D}_{1})\ar[dr]^{i_{f}'}\\
\mathscr{Y}\ar@{^{(}->}[r]_{i_{c}''} & \sx(\mathfrak{D}_{1}^{\sigma})^{\circ}\ar@{^{(}->}[r] & \sx(\mathfrak{D}_{1}^{\sigma}),}
\]
and taking $i_{c}=i_{c}'\circ i_{c}''$ concludes. $\square$

The template for applying these results is as follows. Consider a
pair of connected, reductive groups $\mathbf{G}$ and $\mathbf{H}$
over $\mathbf{Q}$, together with an \emph{L-}homomorphism $^{L}\sigma:{}^{L}\mathbf{G}\to^{L}\!\!\mathbf{H}$
which is \emph{known }to induce a Langlands functoriality map. The
\emph{L}-homomorphism induces a morphism $\mathbf{T}_{\mathbf{H}}\to\mathbf{T}_{\mathbf{G}}$
of unramified Hecke algebras in the usual way. When $\mathbf{G}$
and $\mathbf{H}$ are inner forms of each other, $^{L}\sigma$ is
an isomorphism, and Theorem 5.1.2 (with $\sw_{1}=\sw_{2}$ and $\jmath=\mathrm{id}$)
gives rise to closed immersions of eigenvarieties interpolating correspondences
of Jacquet-Langlands type and/or comparing different theories of overconvergent
automorphic forms. In the general case, the homomorphism $X^{\ast}(\widehat{T}_{\mathbf{H}})\to X^{\ast}(\widehat{T}_{\mathbf{G}})$
together with the natural identification $T(A)\simeq X^{\ast}(\widehat{T})\otimes_{\mathbf{Z}}A$
induces a homomorphism $\tau:T_{\mathbf{H}}(\mathbf{Z}_{p})\to T_{\mathbf{G}}(\mathbf{Z}_{p})$,
and $\jmath$ is given by sending a character $\lambda$ of $T_{\mathbf{G}}$
to the character $(\lambda\circ\tau)\cdot\xi_{\sigma}$ for some \emph{fixed
}character $\xi_{\sigma}$ of $T_{\mathbf{H}}$ which may or may not
be trivial. In this case, Theorem 5.1.6 then induces morphisms of
eigenvarieties interpolating the functoriality associated with $^{L}\sigma$.
In practice, one must carefully choose the character $\xi_{\sigma}$,
compatible tame levels for $\mathbf{G}$ and $\mathbf{H}$, and an
extension of the map on unramified Hecke algebras to include the Atkin-Lehner
operators.

\subsection{Refinements of unramified representations}

In order to apply the interpolation theorems of the previous section,
we need a systematic way of producing sets $\sz^{\mathrm{cl}}$ such
that $\sm_{1}(z)$ consists entirely of classical automorphic forms
for $z\in\sz^{\mathrm{cl}}$. The key is Theorem 3.2.5 together with
Proposition 5.2.1 below.

Let $G\simeq\mathrm{GL}_{n}/\mathbf{Q}_{p}$, with $B$ the upper-triangular
Borel and $\overline{B}$ the lower-triangular Borel. In this case
we may canonically parametrize $L$-valued characters of $\mathcal{A}_{p}$
and unramified characters of $T(\mathbf{Q}_{p})$ by ordered $n$-tuples
$\mathbf{a}=(a_{1},\dots,a_{n})\in(L^{\times})^{n}$, the former via
the map \[
\mathbf{a}\mapsto\chi_{\mathbf{a}}(u_{p,i})=a_{n+1-i}\]
and the latter via the map\[
\mathbf{a}\mapsto\chi_{\mathbf{a}}(t_{1},\dots,t_{n})=\prod_{i=1}^{n}a_{i}^{v_{p}(t_{i})}.\]
Let $\pi$ be an unramified generic irreducible representation of
$G=\mathrm{GL}_{n}(\mathbf{Q}_{p})$ defined over $L$, and let $r:W_{\mathbf{Q}_{p}}\to\mathrm{GL}_{n}(L)$
be the unramified Weil-Deligne representation satisfying $r\simeq\mathrm{rec}(\pi\otimes|\det|^{\frac{1-n}{2}})$.
Let $\varphi_{1},\dots,\varphi_{n}$ be any fixed ordering on the
eigenvalues of $r(\mathrm{Frob}_{p})$, and let $\chi_{\sigma}$,
$\sigma\in S_{n}$ be the character of $\mathcal{A}_{p}$ defined
by $\mathcal{A}_{p}(u_{i})=p^{1-i}\varphi_{\sigma(i)}$.

\textbf{Proposition 5.2.1. }\emph{For every $\sigma\in S_{n}$, the
module $\pi_{p}^{I}$ contains a nonzero vector $v_{\sigma}$ such
that $\mathcal{A}_{p}$ acts on $v_{\sigma}$ through the character
$\chi_{\sigma}$.}

\emph{Proof. }Assembling some results of Casselman (cf. §3.2.2 of \cite{Taibi}),
there is a natural isomorphism of $\Lambda$-modules\[
\pi^{I}\overset{\sim}{\to}(\pi_{\overline{N}})^{T(\mathbf{Z}_{p})}\otimes\delta_{\overline{B}}^{-1}\]
where $t\in\Lambda$ acts on the left-hand side by $U_{t}$. By Theorem
4.2 of \cite{BZ} and Satake's classification of unramified representations,
we may write $\pi$ as the full normalized induction \[
\pi=\mathrm{Ind}_{\overline{B}}^{G}\chi\]
where $\chi$ is the character of $T$ associated with the tuple $(p^{\frac{1-n}{2}}\varphi_{\sigma(n)},\dots,p^{\frac{1-n}{2}}\varphi_{\sigma(1)})$.
By Frobenius reciprocity, there is an embedding of $T$-modules \[
L(\chi\delta_{\overline{B}}^{\frac{1}{2}})\hookrightarrow(\mathrm{Ind}_{\overline{B}}^{G}\chi_{\varphi})_{\overline{N}},\]
so $L(\chi\delta_{\overline{B}}^{-\frac{1}{2}})\hookrightarrow\pi^{I}$,
and $\chi_{\sigma}=\chi\delta_{\overline{B}}^{-\frac{1}{2}}$ upon
noting that $\delta_{\overline{B}}^{-\frac{1}{2}}$ corresponds to
the tuple $(p^{\frac{1-n}{2}},\dots p^{\frac{n-1}{2}})$. $\square$

\subsection{Some quaternionic eigencurves}

Fix a squarefree positive integer $\delta\geq2$, a positive integer
$N$ prime to $\delta$, and a prime $p$ with $p\nmid N\delta$.
Let $D$ be the quaternion division algebra over $\mathbf{Q}$ ramified
at exactly the finite places dividing $\delta$, and ramified or split
over $\mathbf{R}$ according to whether $\delta$ has an odd or even
number of distinct prime divisors. Let $\mathbf{G}$ be the inner
form of $\mathrm{GL}_{2}/\mathbf{Q}$ associated with $D$, and let
$\sx_{D}$ be the eigenvariety $\sx_{\mathbf{G},K_{1}(N)}^{0}$ as
defined in §4.2 and §4.6. Let $\sx$ be the eigenvariety $\sx_{\mathrm{GL}_{2}/\mathbf{Q},K_{1}(N\delta)}^{0}$.
The eigenvarieties $\sx$ and $\sx_{D}$ are both unmixed of dimension
one.

\textbf{Theorem 5.3.1. }\emph{There is a canonical closed immersion
$\iota_{\mathrm{JL}}:\sx_{D}\hookrightarrow\sx$ interpolating the
Jacquet-Langlands correspondence on non-critical classical points.}

\emph{Proof. }Let $\mathfrak{D}_{1}$ and $\mathfrak{D}_{2}$ be the
eigenvariety data giving rise to $\scx_{D}$ and $\scx$ as in Definition
4.3.2, and let $\sm^{\dagger}$ and $\sm_{D}^{\dagger}$ be the sheaves
of automorphic forms carried by $\sx$ and $\sx_{D}$, respectively.
Let $\sz^{\mathrm{cl}}\subset\scz_{1}^{\circ}(\overline{\mathbf{Q}_{p}})\subset(\sw\times\mathbf{A}^{1})(\overline{\mathbf{Q}_{p}})$
be the set of points $z=(\lambda,\alpha^{-1})$ with $\lambda\in\sw(\mathbf{Q}_{p})$
of the form $\lambda(x)=x^{k},\, k\in\mathbf{Z}_{\geq1}$, and with
$v_{p}(\alpha)<k+1$. For any $z\in\sz^{\mathrm{ncc}}$, Theorem 3.2.5
together with the classical Eichler-Shimura isomorphism induces isomorphisms
of Hecke modules\[
\sm^{\dagger}(z)\simeq\left(S_{k+2}(\Gamma_{1}(N\delta)\cap\Gamma_{0}(p))\oplus M_{k+2}(\Gamma_{1}(N\delta)\cap\Gamma_{0}(p))\right)^{U_{p}=\iota(\alpha)}\]
and\[
\sm_{D}^{\dagger}(z)\simeq\left(S_{k+2}^{D}(\Gamma_{1}(N)\cap\Gamma_{0}(p))^{\oplus\epsilon}\right)^{U_{p}=\iota(\alpha)}\]
where $\epsilon=1$ or $2$ according to whether $D$ is ramified
or split over $\mathbf{R}$. The set $\sz^{\mathrm{cl}}$ forms a
Zariski-dense accumulation subset of $\sz_{1}^{\circ}$, and $D_{1}(T,X)(z)$
divides $D_{2}(T,X)(z)$ in $k(z)[X]$ for any $z\in\sz^{\mathrm{cl}}$
by the classical Jacquet-Langlands correspondence. Theorem 5.1.2 now
applies.

\subsection{A symmetric square lifting}

Let $\scc_{0}(N)$ be the cuspidal locus of the Coleman-Mazur-Buzzard
eigencurve of tame level $N$. Given a non-CM cuspidal modular eigenform
$f\in S_{k}(\Gamma_{1}(N))$, Gelbart and Jacquet constructed an cuspidal
automorphic representation $\Pi(\mathrm{sym}^{2}f)$ of $\mathrm{GL}_{3}/\mathbf{Q}$
characterized by the isomorphism \[
\iota\mathrm{WD}(\mathrm{sym}^{2}\rho_{f,\iota}|G_{\mathbf{Q}_{\ell}})\simeq\mathrm{rec}\left(\Pi(\mathrm{sym}^{2}f)_{\ell}\otimes|\det|_{\ell}^{-1}\right)\]
for all primes $\ell$. We are going to interpolate this map into
a morphism $\mathbf{s}:\scc_{0}^{\mathrm{ncm}}(N)\to\sx$, where $\scc_{0}^{\mathrm{ncm}}(N)$
is the Zariski-closure inside $\scc_{0}(N)$ of the classical points
associated with non-CM eigenforms and $\sx$ is an eigenvariety arising
from overconvergent cohomology on $\mathrm{GL}_{3}$. For our intended
application to Galois representations, we need a more precise result.

\textbf{Definition. }
\begin{description}
\item [{i.}] \emph{For a prime $\ell\neq p$, an inertial Weil-Deligne
representation $\tau_{\ell}$ is a pair $(r_{\ell},N_{\ell})$ consisting
of a continuous semisimple representation $r_{\ell}:I_{\mathbf{Q}_{\ell}}\to\mathrm{GL}_{2}(\overline{\mathbf{Q}_{p}})$
and a nilpotent matrix $N_{\ell}\in M_{2}(\overline{\mathbf{Q}_{p}})$
such that $r_{\ell}N_{\ell}=N_{\ell}r_{\ell}$.}
\item [{ii.}] \emph{A }global inertial representation \emph{is a formal
tensor product $\tau=\otimes_{\ell\neq p}\tau_{\ell}$ of inertial
Weil-Deligne representations such that $r_{\ell}(I_{\mathbf{Q}_{\ell}})=\{I\}$
and $N_{\ell}=0$ for all but finitely many $\ell$.}
\end{description}
A global inertial representation $\tau$ has a well-defined conductor
$N(\tau)$, given formally as a product $\prod\ell^{f(\tau_{\ell})}$
of local conductors. 

\textbf{Definition. }\emph{The non-CM cuspidal eigencurve of inertial
type $\tau$, denoted $\scc_{0}^{\mathrm{ncm}}(\tau)$, is the union
of the irreducible components of $\scc(N(\tau))$ which contain a
Zariski-dense set $Z$ of non-CM classical cuspidal points such that
$\mathrm{WD}(\rho_{x}|I_{\mathbf{Q}_{\ell}})\sim\tau_{\ell}$ for
all $\ell\neq p$ and all $x\in Z$.}

To construct $\scc_{0}^{\mathrm{ncm}}(\tau)$ in a manner compatible
with our goals, set $\mathbf{G}=\mathrm{GL_{2}/\mathbf{Q}}$, $N=N(\tau)$
and $\mathbf{T}_{1}=\mathbf{T}_{\mathbf{G}}(K_{1}(N))$. The eigencurve
$\scc(N)$ arises from an eigenvariety datum $\mathfrak{D}_{1}'=(\sw_{1},\sz_{1},\sm_{1},\mathbf{T}_{1},\psi_{1})$
with $\sw_{1}=\mathrm{Hom}_{\mathrm{cts}}(\mathbf{Z}_{p}^{\times},\mathbf{G}_{m})$,
$\sz_{1}$ the Fredholm series of $U_{p}^{4}$ acting on overconvergent
modular forms of tame level $N$, and $\sm_{1}$ the natural spreading
out of overconvergent modular forms of tame level $N$. Let $\mathscr{I}\subset\sco_{\scc(N)}$
be the coherent ideal sheaf cutting out the closed immersion $\scc_{0}^{\mathrm{ncm}}(\tau)\hookrightarrow\scc(N)$,
and set \[
\mathfrak{D}_{1}=\left(\sw_{1},\sz_{1},\pi_{\ast}(\sm_{1}\otimes\sco_{\scc(N)}/\mathscr{I}),\mathbf{T}_{1},\psi\right),\]
so $\scc_{0}^{\mathrm{ncm}}(\tau)$ arises from the eigenvariety datum
$\mathfrak{D}_{1}$. The eigencurve $\scc_{0}^{\mathrm{ncm}}(\tau)$
is unmixed of dimension one. For any classical eigenform $f\in S_{k+2}^{\mathrm{ncm}}(\Gamma_{1}(N))$
of inertial type $\tau$ and $\alpha$ either root of the Hecke polynomial
$X^{2}-a_{f}(p)+p^{k+1}\varepsilon(p)$, we define $\phi_{f,\alpha}:\mathbf{T}_{1}\to\overline{\mathbf{Q}_{p}}$
as the eigenpacket associated with the point $x_{f,\alpha}\in\scc_{0}^{\mathrm{ncm}}(\tau)(\overline{\mathbf{Q}_{p}})$.
We normalize the weight map $\scc_{0}^{\mathrm{ncm}}(\tau)\to\sw_{1}$
so that for $f$ a classical cuspidal eigenform of weight $k+2$ whose
nebentype has $p$-part $\varepsilon$, $w(x_{f})$ corresponds to
the character $t\mapsto t^{k}\varepsilon(t)$.

Now take $\mathbf{H}=\mathrm{GL}_{3}/\mathbf{Q}$ and $\mathbf{T}_{2}=\mathbf{T}_{\mathbf{H}}(K_{1}(N(\mathrm{sym}^{2}\tau)))$.
Let $\mathfrak{D}_{2}$ be the eigenvariety datum from Definition
4.3.2, with $\sx=\sx_{\mathbf{H},K_{1}(N(\mathrm{sym}^{2}\tau))}$
the associated eigenvariety. 

\textbf{Theorem 5.4.1. }\emph{Under the hypotheses above, there is
a morphism $\mathbf{s}:\scc_{0}^{\mathrm{ncm}}(\tau)\to\sx$ interpolating
the symmetric square lift on classical points.}

Let $\jmath:\sw_{1}\hookrightarrow\sw_{2}$ be the closed immersion
sending a character $\lambda$ to the character $\jmath(\lambda)(t_{1},t_{2},t_{3})=\lambda(t_{1}^{2}t_{2})$.
Let $\sigma:\mathbf{T}_{2}\to\mathbf{T}_{1}$ be the map defined on
generators by\begin{eqnarray*}
\sigma(T_{\ell,1}) & = & T_{\ell}^{2}-\ell S_{\ell},\\
\sigma(T_{\ell,2}) & = & T_{\ell}^{2}S_{\ell}-\ell S_{\ell}^{2},\\
\sigma(T_{\ell,3}) & = & S_{\ell}^{3},\\
\sigma(U_{p,1}) & = & U_{p}^{2},\\
\sigma(U_{p,2}) & = & U_{p}^{2}S_{p},\\
\sigma(U_{p,3}) & = & S_{p}^{3}.\end{eqnarray*}

\textbf{Lemma 5.4.2. }\emph{If $f\in S_{k+2}^{\mathrm{ncm}}(\Gamma_{1}(N))$
has inertial type $\tau$ and nebentype $\varepsilon_{f}$, and $X^{2}-a_{f}(p)X+p^{k+1}\varepsilon_{f}(p)$
has a root $\alpha$ with $v_{p}(\alpha)<\frac{k+1}{4}$, then $\Pi(\mathrm{sym}^{2}f)$
contributes to $H^{\ast}(K_{1}(N(\mathrm{sym}^{2}\tau)),\scl_{(2k,k,0)})$,
and $H^{\ast}(K_{1}(N(\mathrm{sym}^{2}\tau)),\scd_{(2k,k,0)})$ contains
a nonzero vector $v$ such that every $T\in\mathbf{T}_{2}$ acts on
$v$ through the scalar $(\phi_{f,\alpha}\circ\sigma)(T)$.}

\emph{Proof. }Fix $f$ and $\alpha$ as in the lemma, and let $\lambda$
be the highest weight $(2k,k,0)$. By the local Langlands correspondence,
$\Pi(\mathrm{sym}^{2}f)$ has conductor exactly $N(\mathrm{sym}^{2}\tau)$.
Since $f$ is non-CM, $\Pi(\mathrm{sym}^{2}f)$ is cuspidal. The Hecke
module $\Pi(\mathrm{sym}^{2}f)^{K_{1}(N(\mathrm{sym}^{2}\tau))I}$
occurs in $H^{\ast}(K_{1}(N(\mathrm{sym}^{2}\tau)),\scl_{\lambda})$
by the Gelbart-Jacquet lifting and the calculations in \cite{Clmotifs}.
At primes $\ell\nmid Np$, $\Pi(\mathrm{sym}^{2}f)_{\ell}$ is unramified,
and $T_{\ell,i}$ acts on the unramified line via the scalar $(\phi_{f}\circ\sigma)(T_{\ell,i})$.
A simple calculation using Proposition 5.2.1 shows that $\Pi(\mathrm{sym}^{2}f)^{K_{1}(N(\mathrm{sym}^{2}\tau))I}$
contains a vector on which $\mathcal{A}_{p}$ acts through the character
associated with the tuple $(p^{-2}\beta^{2},p^{-1}\alpha\beta,\alpha^{2})$,
so there is a vector $v'$ in the $\Pi(\mathrm{sym}^{2}f)$-isotypic
component of $H^{\ast}(K_{1}(N(\mathrm{sym}^{2}\tau)),\scl_{\lambda})$
such that the $\star$-action of $\mathcal{A}_{p}$ is given by the
character associated with the tuple $(\alpha^{-2}\varepsilon_{f}(p)^{2},\varepsilon_{f}(p),\alpha^{2})$.
In particular, $U_{p}$ acts on $v'$ through the scalar $\alpha^{4}\varepsilon(p)$.
By Proposition 3.2.5, the integration map $i_{\lambda}:H^{\ast}(K_{1}(N(\mathrm{sym}^{2}\tau)),\scd_{\lambda})\to H^{\ast}(K_{1}(N(\mathrm{sym}^{2}\tau)),\scl_{\lambda})$
is an isomorphism on the subspace where $U_{p}$ acts with slope $<k+1$,
so $v=i_{\lambda}^{-1}(v')$ does the job. $\square$

Now we take $\sz^{\mathrm{cl}}$ to be the set of points in $\sz_{1}(\overline{\mathbf{Q}_{p}})$
of the form $(\lambda,\alpha^{-1})$, where $\lambda$ is a character
of the form $\lambda(x)=x^{k},k\in\mathbf{Z}_{\geq3}$ and $\alpha$
satisfies $v_{p}(\alpha)<\frac{k+1}{4}$. This is a Zariski-dense
accumulation subset of $\sz_{1}$. By Coleman's classicality criterion,
there is a natural injection $\sm_{1}(z)\to S_{k+2}^{\mathrm{ncm}}(\Gamma_{1}(N)\cap\Gamma_{0}(p))^{U_{p}=\iota(\alpha)}$
of $\mathbf{T}_{1}$-modules, so Theorem 5.1.6 now applies, with the
divisibility hypothesis following from Lemma 5.4.2. We thus conclude.
$\square$

It's not hard to show that the image of $\mathrm{s}$ is actually
a union of irreducible components of $\sx^{0}$.

\subsection{A Rankin-Selberg lifting}

Let $f$ and $g$ be a pair of level one holomorphic cuspidal eigenforms
of weights $k_{f}+2$, $k_{g}+2$ with associated Galois representations
$V_{f,\iota}$ and $V_{g,\iota}$. By a deep theorem of Ramakrishnan \cite{RamRS},
there is a unique isobaric automorphic representation $\Pi(f\otimes g)$
of $\mathrm{GL}_{4}(\mathbf{A}_{\mathbf{Q}})$ characterized by the
equality\[
\mathrm{rec}_{\ell}\left(\Pi(f\otimes g)_{\ell}\otimes|\det|_{\ell}^{-\frac{3}{2}}\right)\simeq\iota\mathrm{WD}(V_{f,\iota}\otimes V_{g,\iota}|G_{\mathbf{Q}_{\ell}})\]
for all primes $\ell$. We are going to interpolate the map $(f,g)\mapsto\Pi(f\otimes g)$
into a morphism of eigenvarieties $\scc_{0}\times\scc_{0}\to\sx$,
where $\scc_{0}$ denotes the cuspidal locus of the Coleman-Mazur
eigencurve and $\sx$ denotes an eigenvariety associated with overconvergent
cohomology on $\mathrm{GL}_{4}/\mathbf{Q}$.

Set $\mathbf{G}=\mathrm{GL_{2}/\mathbf{Q}}$, $\mathbf{T}_{1}=\mathbf{T}_{\mathbf{G}}(K_{1}(1))\otimes\mathbf{T}_{\mathbf{G}}(K_{1}(1))$,
and $\sw_{1}=\mathrm{Hom}_{\mathrm{cts}}(\mathbf{Z}_{p}^{\times}\times\mathbf{Z}_{p}^{\times},\mathbf{G}_{m})$;
we regard an $A$-point of $\sw_{1}$ as a pair of characters $\lambda_{1},\lambda_{2}:\mathbf{Z}_{p}^{\times}\to A^{\times}$
in the obvious way. The product $\scc_{0}\times\scc_{0}$ arises from
an eigenvariety datum $\mathfrak{D}_{1}=(\sw_{1},\sz_{1},\sm_{1},\mathbf{T}_{1},\psi_{1})$,
where $(\lambda_{1},\lambda_{2},\alpha^{-1})\in\sz_{1}(\overline{\mathbf{Q}_{p}})$
if and only if there exist cuspidal overconvergent eigenforms $f_{1}$
and $f_{2}$ of weights $\lambda_{1}$ and $\lambda_{2}$ such that
$U_{p}^{4}\otimes U_{p}^{2}-\alpha$ annihilates $f_{1}\otimes f_{2}$.

Set $\mathbf{H}=\mathrm{GL}_{4}/\mathbf{Q}$, $\mathbf{T}_{2}=\mathbf{T}_{\mathbf{H}}(K_{1}(1))$,
and let $\mathfrak{D}_{2}$ be the eigenvariety datum from Definition
4.3.2, with $\sx$ the associated eigenvariety.

\textbf{Theorem 5.5.1. }\emph{Under the hypotheses above, there is
a morphism $\mathbf{t}:\scc_{0}\times\scc_{0}\to\sx$ interpolating
the Rankin-Selberg lift on classical points.}

Let $\jmath:\sw_{1}\hookrightarrow\sw_{2}$ be the closed immersion
defined by sending a character $\lambda\in\sw_{1}$ to the character\[
\jmath(\lambda)(t)=(t_{1}t_{2})^{-1}\lambda_{1}(t_{1}t_{2})\lambda_{2}(t_{1}t_{3}),\, t=\mathrm{diag}(t_{1},t_{2},t_{3},t_{4})\in T_{\mathbf{H}}.\]
Define a map $\sigma:\mathbf{T}_{2}\to\mathbf{T}_{1}$ on generators
by\begin{eqnarray*}
\sigma(T_{\ell,1}) & = & T_{\ell}\otimes T_{\ell},\\
\sigma(T_{\ell,2}) & = & S_{\ell}\otimes T_{\ell}^{2}+T_{\ell}^{2}\otimes S_{\ell}-2\ell S_{\ell}\otimes S_{\ell},\\
\sigma(T_{\ell,3}) & = & \ell^{-1}S_{\ell}T_{\ell}\otimes S_{\ell}T_{\ell},\\
\sigma(T_{\ell,4}) & = & \ell^{-2}S_{\ell}^{2}\otimes S_{\ell}^{2}\\
\sigma(U_{p,1}) & = & U_{p}\otimes U_{p},\\
\sigma(U_{p,2}) & = & U_{p}^{2}\otimes S_{p},\\
\sigma(U_{p,3}) & = & U_{p}S_{p}\otimes S_{p},\\
\sigma(U_{p,4}) & = & S_{p}^{2}\otimes S_{p}^{2}\end{eqnarray*}
Let $(f,g)$ be an ordered pair as above with\emph{ $k_{f}-1>k_{g}>0$.
}Set \emph{$\lambda=\lambda(f\otimes g):(x_{1},x_{2})\mapsto x_{1}^{k_{f}}x_{2}^{k_{g}}\in\sw_{1}$}.
Let $\alpha_{f},\beta_{f}$ be the roots of the Hecke polynomial $X^{2}-a_{f}(p)X+p^{k_{f}+1}$,
and likewise for $g$.

\textbf{Lemma 5.5.2. }\emph{The module $\Pi(f\otimes g)_{p}^{I}$
contains a vector $v$ on which $\mathcal{A}_{p}$ acts through the
character associated with the tuple }$(p^{-3}\beta_{f}\beta_{g},p^{-2}\beta_{f}\alpha_{g},p^{-1}\alpha_{f}\beta_{g},\alpha_{f}\alpha_{g})$;\emph{
in particular, $U_{p}$ acts via the scalar $p^{2k_{g}+k_{f}-1}\alpha_{f}^{4}\alpha_{g}^{2}$.}

\emph{Proof. }This is a direct consequence of Proposition 5.2.1, together
with the characterization of $\Pi(f\otimes g)$ given above. $\square$

\textbf{Lemma 5.5.3.}\emph{ If $\alpha_{f}$ and $\alpha_{g}$ satisfy
$v_{p}(\alpha_{f}^{4}\alpha_{g}^{2})<\mathrm{min}(k_{f}-k_{g},k_{g}+1)$,
then $H^{\ast}(K_{1}(1),\scd_{\jmath(\lambda)})$ contains a nonzero
vector $v$ such that every $T\in\mathbf{T}_{2}$ acts on $v$ through
the scalar $(\phi_{f,\alpha_{f}}\otimes\phi_{g,\alpha_{g}})(\sigma(T))$.}

\emph{Proof. }Dominance of $\jmath(\lambda)$ is obvious, so $\Pi(f\otimes g)$
is cohomological in the weight $\jmath(\lambda)$. For primes $\ell\nmid p$,
$\Pi(f\otimes g)_{\ell}$ is unramified and $T_{\ell,i}$ acts on
the unramified line via the scalar $(\phi_{f,\alpha_{f}}\otimes\phi_{g,\alpha_{g}})(\sigma(T_{\ell,i}))$. 

Next, recall that the $\star$-action of $\mathcal{A}_{p}$ on $\Pi(f\otimes g)_{p}^{I}$
is simply the usual action rescaled by $\jmath(\lambda)(1,p,p^{2},p^{3})^{-1}$,
and $\jmath(\lambda)$ corresponds to the highest weight $(k_{f}+k_{g}-1,k_{f}-1,k_{g},0)$,
so $\lambda(1,p,p^{2},p^{3})^{-1}=p^{1-2k_{g}-k_{f}}$. In particular,
$\Pi(f\otimes g)_{p}^{I}$ contains a vector on which $U_{p}$ acts
through the scalar $\alpha_{f}^{4}\alpha_{g}^{2}$ by Lemma 5.5.2.
Writing $\kappa=\mathrm{min}(k_{f}-k_{g},k_{g}+1)$, the integration
map $i_{\jmath(\lambda)}$ induces an isomorphism\[
H^{\ast}(K_{1}(1),\scd_{\jmath(\lambda)})_{<\kappa}\overset{\sim}{\to}H^{\ast}(K_{1}(1),\scl_{\jmath(\lambda)})_{<\kappa},\]
and the target contains a vector satisfying the claim of the theorem.
$\square$

Finally, we take $\sz^{\mathrm{cl}}$ to be the set of points in $\sz_{1}$
of the form $(\lambda_{1},\lambda_{2},\alpha^{-1})$ where $\lambda_{1}(x)=x^{k_{1}}$
and $\lambda_{2}(x)=x^{k_{2}}$ with $k_{i}\in\mathbf{Z}$ and $0<k_{2}<k_{1}-1$,
and $\alpha$ satisfies $\alpha<\mathrm{min}(k_{1}-k_{2},k_{2}+1)$.
This is a Zariski-dense accumulation subset, and Theorem 5.1.6 applies
with the divisibility hypothesis following from Lemma 5.5.3. This
proves Theorem 5.5.1. $\square$

\section{Galois representations}

\subsection{Background on trianguline representations and $(\varphi,\Gamma)$-modules}

In this section we give some background and context on $(\varphi,\Gamma)$-modules
and trianguline representations. Our primary references here are the
articles \cite{BCast,BerTri,KPX}. Throughout this section, we let
$K\subset\overline{\mathbf{Q}_{p}}$ denote a finite extension of
$\mathbf{Q}_{p}$, with $K_{0}$ the maximal unramified subfield of
$K$. Let $K_{0}'$ be the maximal unramified subfield of $K(\zeta_{p^{\infty}})$.
We identify $G_{K}\simeq\mathrm{Gal}(\overline{\mathbf{Q}_{p}}/K)$
without comment. Set $H_{K}=\mathrm{Gal}(\overline{\mathbf{Q}_{p}}/K(\zeta_{p^{\infty}}))$
and $\Gamma_{K}=\mathrm{Gal}(K(\zeta_{p^{\infty}})/K)=G_{K}/H_{K}$;
the cyclotomic character $\chi_{\mathrm{cyc}}:\Gamma_{\mathbf{Q}_{p}}\overset{\sim}{\to}\mathbf{Z}_{p}^{\times}$
identifies $\Gamma_{K}$ with an open subgroup of $\mathbf{Z}_{p}^{\times}$.

Let $\mathbf{B}_{\mathrm{rig}}^{\dagger}$ and $\widetilde{\mathbf{B}}_{\mathrm{rig}}^{\dagger}$
be the topological $\mathbf{Q}_{p}$-algebras defined in \cite{BergerInvent}.
These rings are equipped with a continuous action of $G_{\mathbf{Q}_{p}}$
and a commuting operator $\varphi$. The \emph{Robba ring }is the
ring $\srr_{K}=(\mathbf{B}_{\mathrm{rig}}^{\dagger})^{H_{K}}$, with
its natural actions of $\varphi$ and $\Gamma_{K}$. There is an isomorphism\[
\srr_{K}\simeq\left\{ f(\pi_{K})\mid f(T)\in K_{0}'[[T,T^{-1}]]\,\mathrm{with}\, f\,\mathrm{convergent\, on}\, r\leq|T|<1\,\mathrm{for\, some}\, r=r_{f}<1\right\} .\]
Here $\pi_{K}$ is a certain indeterminate arising from the theory
of the field of norms; $\varphi$ and $\Gamma_{K}$ act on the coefficients
of $f$ through the absolute Frobenius and the natural map $\Gamma_{K}\to\mathrm{Gal}(K_{0}'/K_{0})$,
respectively, but the actions on $\pi_{K}$ are noncanonical in general.
The topological ring $\srr_{K}$ is naturally an \emph{LF-space, }i.e.
a strict inductive limit of Fréchet spaces: setting\[
\srr_{K}^{r,s}\simeq\left\{ f(\pi_{K})\mid f(T)\in K_{0}'[[T,T^{-1}]]\,\mathrm{with}\, f\,\mathrm{convergent\, on}\, r\leq|T|\leq s\right\} \]
with its natural affinoid structure, $\srr_{K}=\cup_{r<1}\cap_{s<1}\srr_{K}^{r,s}$.
In particular, for any affinoid algebra $A$, the completed tensor
product $\srr_{K}\widehat{\otimes}A$ is well-defined, and the $\varphi$-
and $\Gamma_{K}$-actions extend naturally.

\textbf{Definition 6.1.1. }\emph{A $(\varphi,\Gamma_{K})$-module
over $L$ is a finite free $\srr_{K}\otimes L$-module $D$ equipped
with commuting $L$-linear and $\srr_{K}$-semilinear actions of $\varphi$
and $\Gamma_{K}$ such that $\srr_{K}\cdot\varphi(D)=D$ and such
that the map $\Gamma_{K}\to\mathrm{End}(D)$ is continuous.}

We write $\mathrm{Mod}(\varphi,\Gamma_{K})_{/L}$ for the category
of $(\varphi,\Gamma_{K})$-modules over $L$. On the other hand, let
$\mathrm{Rep}(G_{K})_{/L}$ denote the category of pairs $(V,\rho)$
where $V$ is a finite-dimensional $L$-vector space and $\rho:G_{K}\to\mathrm{End}_{L}(V)$
is a continuous group homomorphism. The significance of $(\varphi,\Gamma_{K})$-modules
arises from the following beautiful theorem of Berger, Cherbonnier-Colmez,
and Fontaine.

\textbf{Theorem 6.1.2. }\emph{The functor\begin{eqnarray*}
\mathrm{Rep}(G_{K})_{/L} & \to & \mathrm{Mod}(\varphi,\Gamma_{K})_{/L}\\
\rho & \mapsto & \mathbf{D}_{\mathrm{rig}}^{\dagger}(\rho)\end{eqnarray*}
defined by \[
\mathbf{D}_{\mathrm{rig}}^{\dagger}(\rho)=\left(\rho\otimes_{\mathbf{Q}_{p}}\mathbf{B}_{\mathrm{rig}}^{\dagger}\right)^{H_{K}}\]
is additive, tensor-exact and fully faithful, with a natural quasi-inverse
given by\[
D\mapsto\mathbf{V}_{\mathrm{rig}}^{\dagger}(D)=\left(D\otimes_{\srr_{K}}\widetilde{\mathbf{B}}_{\mathrm{rig}}^{\dagger}\right)^{\varphi=1}.\]
}

Historically, this is actually the \emph{third }flavor of $(\varphi,\Gamma)$-module
introduced. Fontaine initiated the theory \cite{FonPG} by constructing
a functor $\rho\mapsto\mathbf{D}(\rho)$ to a similar category, but
with the role of $\srr_{K}$ played by the ring \[
\mathbf{B}_{K}=\left(\lim_{\leftarrow n}(\mathcal{O}_{K_{0}'}/p^{n})[[\pi_{K}]][\pi_{K}^{-1}]\right)[\tfrac{1}{p}].\]
Cherbonnier and Colmez \cite{CCover} then showed, in a difficult
piece of work, that $\mathbf{D}(\rho)$ arises via base change from
a module $\mathbf{D}^{\dagger}(\rho)$ defined over the subring \[
\mathbf{B}_{K}^{\dagger}=\left\{ f(\pi_{K})\in\mathbf{B}_{K}\mid f(T)\,\mathrm{convergent}\,\mathrm{on}\, r(f)<|T|<1\,\mathrm{for}\,\mathrm{some}\, r(f)<1\right\} .\]
Berger \cite{BergerInvent} then defined $\mathbf{D}_{\mathrm{rig}}^{\dagger}$ as the
base change $\mathbf{D}_{\mathrm{rig}}^{\dagger}(\rho)=\mathbf{D}^{\dagger}(\rho)\otimes_{\mathbf{B}_{K}^{\dagger}}\srr_{K}$.
The change from $\mathbf{B}_{K}^{\dagger}$ to $\srr_{K}$ may look
slight, but in fact has enormous consequences. For example, the ring
$\srr_{K}$ is intimately connected with the theory of \emph{p-}adic
differential equations, and Berger used this link to give the first
proof of Fontaine's \emph{p-}adic monodromy conjecture. Berger also
showed how to functorially recover the usual Fontaine modules $\mathbf{D}_{\bullet}(\rho)$
for $\bullet\in\left\{ \mathrm{crys},\mathrm{st,dR},\mathrm{Sen},\mathrm{dif}\right\} $
from $\mathbf{D}_{\mathrm{rig}}^{\dagger}(\rho)$; in particular,
there is a natural equality\[
\mathbf{D}_{\mathrm{crys}}(\rho)=\left(\mathbf{D}_{\mathrm{rig}}^{\dagger}(\rho)[\tfrac{1}{t}]\right)^{\Gamma_{K}}\]
of filtered $\varphi$-modules over $K_{0}\otimes L$, and a $(\varphi,\Gamma_{K})$-equivariant
inclusion\[
\mathbf{D}_{\mathrm{crys}}(\rho)\otimes_{K_{0}}\srr_{K}[\tfrac{1}{t}]\subseteq\mathbf{D}_{\mathrm{rig}}^{\dagger}(\rho)[\tfrac{1}{t}]\]
which is an isomorphism if and only if $\rho$ is crystalline. In
retrospect we might say the Robba ring has {}``just the right amount
of flexibility'': the ring $\mathbf{B}_{K}$ is a field and so its
abstract module theory has little structure, while the ring $\mathbf{B}_{K}^{\dagger}$
is too small to contain the ubiquitous $t$ of $p$-adic Hodge theory,
its elements being bounded as analytic functions on their annuli of
definition.

For our purposes, the following classification of rank one $(\varphi,\Gamma_{K})$-modules
is indispensable \cite{KPX,Nakamura}:

\textbf{Theorem 6.1.3. }\emph{The rank one $(\varphi,\Gamma_{K})$-modules
over $L$ are naturally parametrized by the continuous characters
$\delta:K^{\times}\to L^{\times}$. Writing $\srr_{K}(\delta)$ for
the module corresponding to a character $\delta$, there is a natural
isomorphism\[
\srr_{K}(\delta_{1})\otimes_{\srr_{K}\otimes L}\srr_{K}(\delta_{2})\simeq\srr_{K}(\delta_{1}\delta_{2}).\]
If $\delta(\varpi_{K})\in\mathcal{O}_{L}^{\times}$ for some uniformizer
of $K$, then \[
\srr_{K}(\delta)\simeq\mathbf{D}_{\mathrm{rig}}^{\dagger}(\mathrm{rec}(\delta))\]
where $\mathrm{rec}(\delta)$ is the unique continuous character of
$G_{K}$ such that}\[
\mathrm{rec}(\delta)\circ\mathrm{Art}_{K}=\delta.\]

One very significant point in the study of Galois representations
via their $(\varphi,\Gamma)$-modules is that the latter may be highly
reducible even when $\rho$ is irreducible. To explain this, recall
that Kedlaya \cite{Kedlayaslopes} associated with any $\varphi$-module
$D$ over $\srr_{K}$ a finite set of slopes $\mathbf{s}(D)\subset\mathbf{Q}$,
and proved that $D$ admits a unique filtration $0=D_{0}\subset D_{1}\subset\cdots D_{j}\subset D_{j+1}=D$
by $\varphi$-submodules such that each $\mathbf{s}(D_{i+1}/D_{i})$
is a singleton, say $\mathbf{s}(D_{i+1}/D_{i})=\{s_{i}(D)\}$, and
such that $s_{i}(D)>s_{i-1}(D)$ for all $1\leq i\leq j$. One says
that each $D_{i+1}/D_{i}$ is \emph{pure of slope $s_{i}$. }Berger \cite{BergerInvent}
then proved that the essential image of $\mathbf{D}_{\mathrm{rig}}^{\dagger}$
consists exactly of those $(\varphi,\Gamma_{K})$-modules whose underlying
$\varphi$-module is pure of slope zero, or \emph{étale}. However,
Kedlaya's slope filtration theorem does not preclude étale $(\varphi,\Gamma_{K})$-modules
from containing subobjects whose slopes are \emph{positive}.

Let us now formalize the notion of a Galois representation whose $(\varphi,\Gamma_{K})$-module
is totally reducible. Given a continuous representation $\rho:G_{K}\to\mathrm{GL}_{n}(L)$,
an ordered $n$-tuple $\delta=(\delta_{1},\dots,\delta_{n})$ of continuous
characters $\delta_{i}:K^{\times}\to L^{\times}$ is a \emph{parameter
of $\rho$ }if $\mathbf{D}_{\mathrm{rig}}^{\dagger}(\rho)$ admits
a filtration\[
0=\mathrm{Fil}^{0}\subset\mathrm{Fil}^{1}\subset\dots\subset\mathrm{Fil}^{n}=\mathbf{D}_{\mathrm{rig}}^{\dagger}(\rho)\]
by $(\varphi,\Gamma_{K})$-stable $\srr_{K}\otimes L$-free direct
summands such that $\mathrm{Fil}^{i}/\mathrm{Fil}^{i-1}\simeq\srr_{K}(\delta_{i})$
for all $1\leq i\leq n$. Let $\ppar(\rho)$ denote the set of parameters
of $\rho$.%
\footnote{In the global setup of §1.2, we have $\ppar(\rho)=\prod_{v|p}\ppar(\rho|G_{F_{v}})$. %
} Note that a given representation may not admit any parameters at
all.

\textbf{Definition 6.1.4}. \cite{ColTri} \emph{A representation
$\rho$ is }trianguline \emph{if $\ppar(\rho)$ is nonempty.}

The most well-studied trianguline representations are the \emph{nearly
ordinary representations, }in which case the representation space
of $\rho$ itself admits a $G_{K}$-stable full flag $0=V^{(0)}\subset V^{(1)}\subset\cdots\subset V^{(n)}=V$
such that $\mathrm{Fil}^{i}=\mathbf{D}_{\mathrm{rig}}^{\dagger}(\rho^{(i)})$;
a parameter $\delta$ corresponds to a nearly ordinary structure on
$\rho$ if and only if $\delta_{i}(\varpi_{K})\in\mathcal{O}_{L}^{\times}$
for $1\leq i\leq n$. However, most trianguline representations are
irreducible in the category of Galois representations.

\subsubsection*{Example 6.1.1: de Rham representations}

Suppose $\rho:G_{K}\to\mathrm{GL}_{n}(\overline{\mathbf{Q}_{p}})$
is de Rham, and let $\mathrm{WD}(\rho)=(r,N)$ denote the associated
Frobenius-semisimple Weil-Deligne representation.

\textbf{Proposition 6.1.5. }\emph{The following are equivalent:}
\begin{description}
\item [{i.}] \emph{$\rho$ is trianguline.}
\item [{ii.}] \emph{$\rho$ becomes semistable over an abelian extension
of $K$.}
\item [{iii.}] \emph{$r$ is a direct sum of characters.}
\item [{iv.}] \emph{The local Langlands correspondent $\pi$ of $\mathrm{WD}(\rho)$
is a subquotient of a representation induced from a Borel subgroup
of $\mathrm{GL}_{n}(K)$.}
\end{description}
\emph{Proof (sketch).} The equivalence of i. and ii. follows from
Berger's dictionary \cite{BerphiNfilt} relating filtered $(\varphi,N,G_{K})$-modules
and $(\varphi,\Gamma_{K})$-modules. The equivalence of ii. and iii.
is an easy consequence of Fontaine's construction of $\mathrm{WD}(\rho)$
in terms of $\mathbf{D}_{\mathrm{pst}}(\rho)$. The equivalence of
iii. and iv. follows from Bernstein and Zelevinsky's work \cite{BZ,Zelevinsky},
and in particular from the fact that $r$ determines the cuspidal
support of $\pi$. $\square$\emph{ }

\subsubsection*{Example 6.1.2: Overconvergent modular forms}

Throughout this example, and the remainder of §6, we write $f$ for
an overconvergent cuspidal eigenform of finite slope and some tame
level $N_{f}$, with associated Galois representation $\rho_{f}:G_{\mathbf{Q}}\to\mathrm{GL}_{2}(L)$.
We define the \emph{weight }of $f$ as the unique continuous character
$w:\mathbf{Z}_{p}^{\times}\to L^{\times}$ such that $\det\rho_{f}|I_{\mathbf{Q}_{p}}\simeq w(\chi_{\mathrm{cyc}}^{-1})\chi_{\mathrm{cyc}}^{-1}$.
Set $k=2+\frac{\log w(1+p)}{\log(1+p)}$, so the Sen weights of $\rho_{f}$
are exactly $0$ and $k-1$. Let $\alpha_{f}$ be the $U_{p}$-eigenvalue
of $f$. We recall a brilliant result of Kisin \cite{KisinOCFM}:

\textbf{Theorem (Kisin). }\emph{The space of crystalline periods $\mathbf{D}_{\mathrm{crys}}^{+}(\rho_{f}|G_{\mathbf{Q}_{p}})^{\varphi=\alpha_{f}}$
is nonzero.}

We now associate a canonical parameter $\delta_{f}$ with $f$. Before
doing so, we partition the set of finite-slope overconvergent eigenforms
$f$ into five types:

1a. $k\in\mathbf{Z}_{\geq2}$ and $0\leq v(\alpha_{f})<k-1$.

1b. $k\in\mathbf{Z}_{\geq1}$, $v_{p}(\alpha_{f})=k-1$, and $\rho_{f}|G_{\mathbf{Q}_{p}}$
is de Rham and indecomposable.

2. $k\in\mathbf{Z}_{\geq1}$, $v_{p}(\alpha_{f})=k-1$, and $\rho_{f}|G_{\mathbf{Q}_{p}}$
is a direct sum of characters.

3a. $k\in\mathbf{Z}_{\geq1}$, $v_{p}(\alpha_{f})\geq k-1$, and $\rho_{f}|G_{\mathbf{Q}_{p}}$
is not de Rham.

3b. $k\notin\mathbf{Z}_{\geq1}$.

Forms of type 1 are always classical, while forms of type 3 are never
classical. 

If $f$ is of type 1 or 3b, then we define\[
\delta_{f}=(\mu_{\alpha},\mu_{\alpha^{-1}\eta}w(x_{0})^{-1}x_{0}^{-1})\]
where $\eta=\det\rho_{f}(\mathrm{Frob}_{p})$; here for $x\in\mathbf{Q}_{p}^{\times}$
we set $x_{0}=x|x|$ and $\mu_{\alpha}(x)=\alpha^{v_{p}(x)}$. If
$f$ is of type 2 or 3a, we define\[
\delta_{f}=(\mu_{p^{1-k}\alpha}x_{0}^{1-k},\mu_{\alpha^{-1}\eta}\varepsilon(x_{0})^{-1})\]
where $w(x)=x^{k-2}\varepsilon(x)$ with $\varepsilon$ of finite
order. (For a proof that these really are parameters, see e.g. Proposition
5.2 of \cite{ChPlongent}.) These are the unique parameters compatible
with Kisin's result, in the sense that the associated $\mathrm{Fil}^{1}\subset\mathbf{D}_{\mathrm{rig}}^{\dagger}(\rho_{f}|G_{\mathbf{Q}_{p}})$
satisfies\[
0\neq\mathbf{D}_{\mathrm{crys}}(\mathrm{Fil}^{1})\subseteq\mathbf{D}_{\mathrm{crys}}(\rho_{f}|G_{\mathbf{Q}_{p}})^{\varphi=\alpha_{f}}.\]

\subsection{The general conjecture}

Let $A$ be a finite étale $\mathbf{Q}_{p}$-algebra, and fix an integer
$n\geq1$. Let $\st=\st_{n,A}$ be the rigid analytic space parametrizing
continuous characters of $(A^{\times})^{n}$. $ $Let $\delta\in\st(\overline{\mathbf{Q}_{p}})$
be any character, so we may regard $\delta$ as an ordered $n$-tuple
of characters $\delta_{i}:A^{\times}\to\overline{\mathbf{Q}_{p}}^{\times}$
in the obvious way. Set $\Sigma=\mathrm{Hom}_{\mathbf{Q}_{p}-\mathrm{alg}}(A,\overline{\mathbf{Q}_{p}})$.
Given $\delta=(\delta_{1},\dots,\delta_{n})$ as above, there is some
sufficiently small open subgroup $U\subset A$ such that the linear
form $(\partial\delta_{i})(a)=$$\log(\delta_{i}(\mathrm{exp}(a))):U\to\overline{\mathbf{Q}_{p}}$
is well-defined, and we have $(\partial\delta_{i})(a)=\sum_{\sigma\in\Sigma}k_{\sigma,i}(\delta)\sigma(a)$
for some uniquely determined constants $k_{\sigma,i}(\delta)\in\overline{\mathbf{Q}_{p}}$.
When $A=K$ is a finite extension of $\mathbf{Q}_{p}$ and $\delta$
is a parameter of a trianguline representation $\rho:G_{K}\to\mathrm{GL}_{n}(L)$,
the multiset $\{-k_{\sigma,i}(\delta)\}_{1\leq i\leq n}$ is exactly
the set of $\sigma$-Sen weights of $\rho$.

\textbf{Definition 6.2.1. }\emph{$W(\delta)_{\sigma}<S_{n}$ is the
group generated by the transpositions $(ij)$ such that $k_{\sigma,i}(\delta)-k_{\sigma,j}(\delta)\in\mathbf{Z}$,
and \[
W(\delta)=\prod_{\sigma}W(\delta)_{\sigma}<S_{n}^{\Sigma}.\]
The reflections in $W(\delta)$ are those elements which are a transposition
in exactly one $W(\delta)_{\sigma}$ and trivial in $W(\delta)_{\sigma'}$
for all $\sigma'\neq\sigma$. }

For $g\in W(\delta)$, we define $g\cdot\delta$ by \[
(g\cdot\delta)_{i}(x)=\delta_{i}(x)\prod_{\sigma}\sigma(x)^{k_{\sigma,g^{-1}(i)}(\delta)-k_{\sigma,i}(\delta)}\qquad\forall1\leq i\leq n.\]
It's easy to check that this is a left action. Let $W(\delta)\cdot\delta$
be the orbit of $\delta$.

\textbf{Definition 6.2.2. }\emph{Given a character $\eta\in W(\delta)\cdot\delta$
and a reflection $g=(ij)_{\sigma}\in W(\delta)$ with $i<j$, $\eta$
}precedes $g\cdot\eta$ if $k_{\sigma,i}(\eta)-k_{\sigma,j}(\eta)\in\mathbf{Z}_{<0}$.
\emph{We define a partial ordering on $W(\delta)\cdot\delta$ by $\eta\preceq\eta'$
if either $\eta=\eta'$ or if there is a chain of characters $\eta_{0},\eta_{1},\dots,\eta_{N}\in W(\delta)\cdot\delta$
with $\eta_{0}=\eta$ and $\eta_{N}=\eta'$ such that $\eta_{i-1}$
precedes $\eta_{i}$ for all $1\leq i\leq N$.}

For any $\delta\in\st(\overline{\mathbf{Q}_{p}})$ we set \[
\st[\delta]=\left\{ \eta\in W(\delta)\cdot\delta\mid\delta\preceq\eta\right\} .\]
Now in the setting of Conjecture 1.2.5 we simply specialize all these
notions to $A=F\otimes_{\mathbf{Q}}\mathbf{Q}_{p}$, in which case\begin{eqnarray*}
\st_{n,F}(\overline{\mathbf{Q}_{p}}): & = & \st_{n,F\otimes_{\mathbf{Q}}\mathbf{Q}_{p}}(\overline{\mathbf{Q}_{p}})\\
 & = & \prod_{v|p}\st_{n,F_{v}}(\overline{\mathbf{Q}_{p}}).\end{eqnarray*}

Let us give a slight reformulation of Conjecture 1.2.5. For $K/\mathbf{Q}_{p}$
finite and $\rho:G_{K}\to\mathrm{GL}_{n}(\overline{\mathbf{Q}_{p}})$
a Galois representation, $\ppar(\rho)$ is naturally a subset of $\st_{n,K}$,
and we set (with $A=K$) \begin{eqnarray*}
\st[\rho] & = & \bigcup_{\delta\in\ppar(\rho)}\st[\delta]\\
 & \subset & \st_{n,K}(\overline{\mathbf{Q}_{p}}).\end{eqnarray*}
Given a global Galois representation $\rho:G_{F}\to\mathrm{GL}_{n}(\overline{\mathbf{Q}_{p}})$
as in §1.2, we define\begin{eqnarray*}
\st[\rho]^{\mathrm{loc}} & = & \prod_{v|p}\st[\rho|G_{F_{v}}]\\
 & \subset & \prod_{v|p}\st_{n,F_{v}}(\overline{\mathbf{Q}_{p}})\\
 & = & \st_{n,F}(\overline{\mathbf{Q}_{p}}).\end{eqnarray*}
Recall, on the other hand, the automorphically defined set $\st[\rho]\subset\st_{n,F}(\overline{\mathbf{Q}_{p}})$.
The following conjecture is easily seen to be equivalent to Conjecture
1.2.5.

\textbf{Conjecture 6.2.3. }\emph{Notation and assumptions as in }§\emph{1.2,
we have\[
\st[\rho]=\st[\rho]^{\mathrm{loc}}.\]
}

\textbf{Remarks. }

1. As we've already remarked, this conjecture is strongly analogous
with Serre's modularity conjecture and its generalizations, and the
formulation just given highlights this analogy. Furthering the analogy,
all of our partial results are of the form {}``for certain representations
$\rho$ and certain elements $\eta\in\st[\rho]^{\mathrm{loc}}$ we
have $\eta\in\st[\rho]$''. On the other hand, there doesn't seem
to be any obvious mod $p$ structure analogous with the set of parameters
of $\rho$ and the role they play in describing the total set $\st[\rho]^{\mathrm{loc}}$
.

2. This conjecture is rather nontrivial even when $\rho$ is classically
automorphic. For concreteness, suppose $\rho:G_{\mathbf{Q}}\to\mathrm{GL}_{n}(L)$
is crystalline with $n$ distinct Hodge-Tate weights. A \emph{refinement
}of $\rho$, in the terminology of \cite{BCast}, is a choice of an
ordering $\alpha_{\bullet}=\{\alpha_{1},\dots,\alpha_{n}\}$ on the
eigenvalues of the crystalline Frobenius $\varphi$ on $\mathbf{D}_{\mathrm{crys}}(\rho|G_{\mathbf{Q}_{p}})$.
Refinements are in bijection with complete $\varphi$-stable flags
$0=\mathcal{F}_{0}\subset\mathcal{F}_{1}\subset\cdots\subset\mathcal{F}_{n}=\mathbf{D}_{\mathrm{crys}}(\rho|G_{\mathbf{Q}_{p}})$
of $L$-vector spaces by the association $\mathcal{F}_{i}=\ker\prod_{1\leq j\leq i}(\varphi-\alpha_{j})$.
Any such flag, in turn, determines a triangulation of $\mathbf{D}_{\mathrm{rig}}^{\dagger}(\rho|G_{\mathbf{Q}_{p}})$
by setting $\mathrm{Fil}^{i}=\mathcal{F}_{i}\otimes_{\mathbf{Q}_{p}}\srr_{\mathbf{Q}_{p}}[\tfrac{1}{t}]\cap\mathbf{D}_{\mathrm{rig}}^{\dagger}(\rho|G_{\mathbf{Q}_{p}})$,
the intersection taking place in $\mathbf{D}_{\mathrm{crys}}(\rho|G_{\mathbf{Q}_{p}})\otimes_{\mathbf{Q}_{p}}\srr_{\mathbf{Q}_{p}}[\tfrac{1}{t}]\cong\mathbf{D}_{\mathrm{rig}}^{\dagger}(\rho|G_{\mathbf{Q}_{p}})[\tfrac{1}{t}]$.
It's easy to check that these associations determine bijections between
these three sets (orderings on $\varphi$-eigenvalues, complete $\varphi$-stable
flags, triangulations). Given a refinement $\alpha_{\bullet}$ with
$\delta=(\delta_{1},\dots,\delta_{n})$ the parameter of the associated
triangulation, it transpires that $\delta_{i}(x)=\mu_{\alpha_{i}}(x)x^{-k_{i}}$
where $k_{1},\dots,k_{n}$ is some ordering on the Hodge-Tate weights
of $\rho$. In particular, $W(\delta)=S_{n}$, and there's a unique
element $\delta^{\mathrm{cl}}\in W(\delta)\cdot\delta$ which is maximal
for the partial ordering $\preceq$. Again following \cite{BCast},
we say the refinement (or triangulation) is \emph{noncritical }if
$k_{1}<k_{2}<\cdots<k_{n}$, and \emph{critical }otherwise. It's easy
to see that the following conditions are equivalent:

i. The refinement $\alpha_{\bullet}$ is noncritical.

ii. The weight $\lambda$ of the putative point $x(\rho,\delta)$
is $B$-dominant.

iii. $\delta=\delta^{\mathrm{cl}}$, or equivalently, $\st[\delta]=\{\delta\}$.

Now the weight $\lambda$ of the putative point $x^{\mathrm{cl}}$
characterized by $\rho_{x}\simeq\rho$ and $\delta_{x}=\delta^{\mathrm{cl}}$
is $B$-dominant, and the subspace of $H^{\ast}(K_{1}(N),\scd_{\lambda})$
predicted by Conjecture 1.2.5 is compatible with the Fontaine-Mazur-Langlands
conjecture under the map\[
H^{\ast}(K_{1}(N),\scd_{\lambda})\to H^{\ast}(K_{1}(N),\scl_{\lambda})\]
(this calculation is sketched in \cite{Hthesis}). Qualitatively,
Conjecture 1.2.5 says the classical point $x^{\mathrm{cl}}$ has {}``companion
points'' if and only if $\delta$ is critical, and $x^{\mathrm{cl}}$
has more companion points the further $\delta$ is from being noncritical.

3. There is a purely local analogue of Conjecture 6.2.3, which seems
quite interesting in its own right. To formulate this, fix a finite
field $\mathbf{F}$ of characteristic $p$ and an absolutely irreducible
residual representation $\overline{\rho}:G_{K}\to\mathrm{GL}_{n}(\mathbf{F})$.
Let $\scx_{\overline{\rho}}$ be the rigid generic fiber of the universal
pseudodeformation space of $\overline{\rho}$, and let $X(\overline{\rho})\subset\scx_{\overline{\rho}}\times\st_{n,K}$
be the \emph{finite slope space }(Definition 2.10 of \cite{HeSch}).

\textbf{Conjecture 6.2.4. }\emph{For any point $x\in\scx_{\overline{\rho}}(\overline{\mathbf{Q}_{p}})$
with associated representation $\rho_{x}:G_{K}\to\mathrm{GL}_{n}(\overline{\mathbf{Q}_{p}})$,
we have\[
\mathrm{pr}_{2}\left(\mathrm{pr}_{1}^{-1}(x)\cap X(\overline{\rho})\right)=\st[\rho_{x}]\]
as subsets of $\st_{n,K}(\overline{\mathbf{Q}_{p}})$.}

We also strongly believe that if there is a functor\begin{eqnarray*}
\mathrm{Irr}_{n}(G_{K})_{/L} & \to & \mathrm{Ban}_{\mathrm{adm}}(\mathrm{GL}_{n}(K))_{/L}\\
\rho & \mapsto & \Pi(\rho)\end{eqnarray*}
which deserves to be called the \emph{p-}adic local Langlands correspondence,
then (modulo details of normalization) the set $\st[\rho]$ is exactly
the set of characters appearing in the locally analytic Jacquet module
of $\Pi(\rho)$.

\subsection{Evidence for two-dimensional Galois representations}

In this section we prove Theorem 1.2.6. Given $f$ an overconvergent
finite-slope eigenform, we write $N(\rho_{f})$ for the prime-to-$p$
Artin conductor of $\rho_{f}$. Note that \emph{a priori} the only
relation between the tame level $N_{f}$ and the Artin conductor is
the divisibility $N(\rho_{f})|N_{f}$, since we don't require that
$f$ be a newform. However, we have the following result.

\textbf{Proposition 6.3.1 (Level-lowering). }\emph{If $\overline{\rho}_{f}$
satisfies the hypotheses of Theorem 1.2.6, then there exists a finite-slope
eigenform $f_{0}$ of tame level $ $$N_{f_{0}}=N(\rho_{f})$ such
that $\rho_{f_{0}}\simeq\rho_{f}$.}

\emph{Proof. }This follows from Emerton's \emph{p-}adic version of
Mazur's principle (Theorem 6.2.4 of \cite{EmertonLG}) together
with the results in \cite{PaulinLG}. $\square$

Now, let $\rho$ be as in Theorem 1.2.6. By Corollary 1.2.2 of \cite{EmertonLG},
there is some $f$ and some continuous character $\nu:G_{\mathbf{Q}}\to\overline{\mathbf{Q}_{p}}^{\times}$
such that $\rho\simeq\rho_{f}\otimes\nu$. We may choose $\nu$ and
$f$ such that $\nu$ is unramified outside $p$ and $\infty$, and
then by the previous proposition we may choose $f$ such that $N_{f}=N(\rho)$.
Since Conjecture 1.2.5 is compatible with twisting by characters unramified
outside $p\infty$, we may assume $\nu=1$ and $\rho=\rho_{f}$. Let
$w:\mathbf{Z}_{p}^{\times}\to L^{\times}$ be the weight of $f$ as
in §6.1 above, and let $S_{w}^{\dagger}(\Gamma_{1}(N))$ be the linear
span of finite-slope overconvergent cusp forms of weight $w$ and
level $N$. For any character $\lambda=(\lambda_{1},\lambda_{2})$
in the $\mathrm{GL}_{2}/\mathbf{Q}$ weight space, work of Stevens \cite{Stfamilies}
and Bellaïche \cite{BeCrit} yields a noncanonical injection of semisimplified
Hecke modules\[
\beta_{\lambda}:S_{\lambda_{1}\lambda_{2}^{-1}}^{\dagger}(\Gamma_{1}(N))^{\mathrm{ss}}\otimes\lambda_{2}\circ(\mathrm{det}|\det|)\hookrightarrow H^{1}(K_{1}(N),\scd_{\lambda})^{\mathrm{ss}}\]
where of course the right-hand side denotes overconvergent cohomology
for $\mathrm{GL}_{2}$.

In general we have $|\ppar(\rho_{f})|\leq2$. Suppose first that $\delta_{f}$
is the unique parameter of $\rho_{f}$. This is true if and only if
either $f$ is type 1 and Steinberg at $p$, or type 3. If $f$ is
type 1 or type 3b, the point $x=x(\rho_{f},\delta_{f})$ we seek is
the point associated with the eigenspace of $f$ under the map $\beta_{(w,0)}$.
If $f$ is type 3a, with weight $w(x)=x^{k-2}\varepsilon$, then by
deep work of Coleman, $f$ has a \emph{companion form $g$} \cite{Colemanclassicaloc},
namely a form $g$ of type 3b such that $\rho_{f}\simeq\rho_{g}\otimes\chi_{\mathrm{cyc}}^{1-k}$
and $\alpha_{f}=p^{k-1}\alpha_{g}$. It's easy to see that $g$ has
weight $x^{-k}\varepsilon$, and the eigenspace of $g\otimes1$ under
the map $\beta_{(x^{-1}\varepsilon,x^{k-1})}$ corresponds to the
point $x(\rho_{f},\delta_{f})$ predicted by Conjecture 1.2.4.

Now suppose that $|\ppar(\rho_{f})|=2$, so $f$ is type 1 or 2. The
point $x(\rho_{f},\delta_{f})$ is given exactly as in the type 1
or type 3b subcases of the previous case, respectively. There is a
classical newform $\mathbf{f}\in S_{k}(\Gamma_{1}(Np^{n}))$ for some
$n\geq0$ which is {}``partially unramified principal series'' at
$p$, such that $f$ is in the same generalized eigenspace as the
refinement $\mathbf{f}_{\alpha}$ of $\mathbf{f}$. Let $\varepsilon_{\mathbf{f}}$
be the $p$-part of the nebentypus of $\mathbf{f}$. The second parameter
$\delta'$ is characterized by the fact that $\mathbf{f}\otimes\varepsilon_{\mathbf{f}}^{-1}$
admits a type 1 refinement $h$ with parameter $\delta'\otimes\varepsilon_{\mathbf{f}}$,
and $x(\rho_{f},\delta')$ corresponds to the eigenspace of $h$ under
the map $\beta_{(w_{h}\varepsilon_{\mathbf{f}},\varepsilon_{\mathbf{f}})}$.
$\square$

\subsection{Evidence for three- and four-dimensional Galois representations}

In this section we prove Theorems 1.2.7 and 1.2.8. 

\emph{Proof of Theorem 1.2.7. }Notation as in the theorem, let us
show that the character $\mathrm{sym}^{2}\delta_{f}$ really is a
parameter of $\mathrm{sym}^{2}\rho_{f}$. For brevity set $\mathbf{D}_{f}=\mathbf{D}_{\mathrm{rig}}^{\dagger}(\rho_{f}|G_{\mathbf{Q}_{p}})$
and $\mathbf{D}_{\mathrm{sym}^{2}f}=\mathbf{D}_{\mathrm{rig}}^{\dagger}(\mathrm{sym}^{2}\rho_{f}|G_{\mathbf{Q}_{p}})$.
We may realize $\mathbf{D}_{\mathrm{sym}^{2}f}$ as the $\srr$-span
of symmetric tensors in $\mathbf{D}_{f}\otimes_{\srr}\mathbf{D}_{f}$.
Suppose $\mathbf{D}_{f}$ has a triangulation\[
0\to\srr(\delta_{1})\to\mathbf{D}_{f}\to\srr(\delta_{2})\to0.\]
Let $v_{1},v_{2}$ be a basis for $\mathbf{D}_{f}$ with $v_{1}$
spanning $\srr(\delta_{1})$. The filtration\begin{multline*}
0\subset\srr(\delta_{1}^{2})\simeq\mathrm{Span}_{\srr}(v_{1}\otimes v_{1})\subset\mathrm{Span}_{\srr}(v_{1}\otimes v_{1},v_{1}\otimes v_{2}+v_{2}\otimes v_{1})\\
\subset\mathrm{Span}_{\srr}(v_{1}\otimes v_{1},v_{1}\otimes v_{2}+v_{2}\otimes v_{1},v_{2}\otimes v_{2})\simeq\mathbf{D}_{\mathrm{sym}^{2}f}\end{multline*}
then exhibits $(\delta_{1}^{2},\delta_{1}\delta_{2},\delta_{2}^{2})$
as an element of $\ppar(\mathrm{sym}^{2}\rho_{f})$.

Since $\mathrm{sym}^{2}\rho_{f}$ is assumed irreducible, $\rho_{f}$
is neither reducible or dihedral up to twist. If $f$ is of type 1
or 3b, $f$ defines a point $x_{f}=x(\rho_{f},\delta_{f})\in\scc_{0}^{\mathrm{ncm}}(\tau)$
(with notation as in §5.4 and §6.3), and the point $\mathbf{s}(x_{f})\in\sx_{\mathrm{GL}_{3}/\mathbf{Q},N(\mathrm{sym}^{2}\tau)}$
is the point predicted by Conjecture 1.2.4. If $f$ is of type 2 or
3a with companion form $g$, we take a suitable twist of $\mathbf{s}(x_{g})$.
$\square$

\emph{Proof of Theorem 1.2.8. }Notation as in the theorem, we first
demonstrate that the claimed characters are parameters. Suppose $\mathbf{D}_{f}=\mathbf{D}(\rho_{f}|G_{\mathbf{Q}_{p}})$
and $\mathbf{D}_{g}=\mathbf{D}(\rho_{g}|G_{\mathbf{Q}_{p}})$ admit
triangulations\[
0\to\srr(\delta_{f,1})\to\mathbf{D}_{f}\to\srr(\delta_{f,2})\to0\]
and likewise for $\mathbf{D}_{g}$. Let $v_{1},v_{2}$ be a basis
for $\mathbf{D}_{f}$ with $v_{1}$ generating $\srr(\delta_{f,1})$,
and let $w_{1},w_{2}$ be an analogous basis for $\mathbf{D}_{g}$.
The character $\delta_{f}\boxtimes\delta_{g}$ is then the parameter
associated with the triangulation\[
0\subset\srr(\delta_{f,1})\otimes\srr(\delta_{g,1})\subset\srr(\delta_{f,1})\otimes\mathbf{D}_{g}\subset\mathrm{Span}_{\srr}\left(v_{1}\otimes w_{1},v_{1}\otimes w_{2},v_{2}\otimes w_{1}\right)\subset\mathbf{D}_{f}\otimes\mathbf{D}_{g},\]
and $\delta_{g}\boxtimes\delta_{f}$ corresponds to the triangulation
\[
0\subset\srr(\delta_{f,1})\otimes\srr(\delta_{g,1})\subset\srr(\delta_{g,1})\otimes\mathbf{D}_{f}\subset\mathrm{Span}_{\srr}\left(v_{1}\otimes w_{1},v_{2}\otimes w_{1},v_{1}\otimes w_{2}\right)\subset\mathbf{D}_{f}\otimes\mathbf{D}_{g},\]

Notation as in §5.5, if $f,g$ are both of type 1 and/or type 3b then
$\mathbf{t}(x_{f},x_{g})$ and $\mathbf{t}(x_{g},x_{f})$ (with $x_{f},x_{g}$
as in the previous proof) are exactly the points $x(\rho_{f}\otimes\rho_{g},\delta_{f}\boxtimes\delta_{g})$
and $x(\rho_{f}\otimes\rho_{g},\delta_{g}\boxtimes\delta_{f})$. The
general case is similar: the points we seek are given by twisting
the points obtained by applying $\mathbf{t}$ to the companion points
of $f$ and/or $g$. $\square$

\appendix

\section{Some commutative algebra}

In this appendix we collect some results relating the projective dimension
of a module $M$ and its localizations, the nonvanishing of certain
Tor and Ext groups, and the heights of the associated primes of $M$.
We also briefly recall the definition of a perfect module, and explain
their basic properties. These results are presumably well-known to
experts, but they are not given in our basic reference \cite{Matcrt}.

Throughout this subsection, $R$ is a commutative Noetherian ring
and $M$ is a finite $R$-module. Our notations follow \cite{Matcrt},
with one addition: we write $\mathrm{mSupp}(M)$ for the set of maximal
ideals in $\mathrm{Supp}(M)$.

\textbf{Proposition A.1. }\emph{There is an equivalence\[
\mathrm{projdim}_{R}(M)\geq n\Leftrightarrow\mathrm{Ext}_{R}^{n}(M,N)\neq0\,\mathrm{for\, some}\, N\in\mathrm{Mod}_{R}.\]
}See e.g. p. 280 of \cite{Matcrt} for a proof.

\textbf{Proposition A.2. }\emph{The equality\[
\mathrm{projdim}_{R}(M)=\mathrm{sup}_{\mathfrak{m}\in\mathrm{mSupp}(M)}\mathrm{projdim}_{R_{\mathfrak{m}}}(M_{\mathfrak{m}})\]
holds.}

\emph{Proof.} Any projective resolution of $M$ localizes to a projective
resolution of $M_{\mathfrak{m}}$, so $\mathrm{projdim}_{R_{\mathfrak{m}}}(M_{\mathfrak{m}})\leq\mathrm{projdim}_{R}(M)$
for all $\mathfrak{m}$. On the other hand, if $\mathrm{projdim}_{R}(M)\geq n$,
then $\mathrm{Ext}_{R}^{n}(M,N)\neq0$ for some $N$, so $\mathrm{Ext}_{R}^{n}(M,N)_{\mathfrak{m}}\neq0$
for some $\mathfrak{m}$; but $\mathrm{Ext}_{R}^{n}(M,N)_{\mathfrak{m}}\simeq\mathrm{Ext}_{R_{\mathfrak{m}}}^{n}(M_{\mathfrak{m}},N_{\mathfrak{m}})$,
so $\mathrm{projdim}_{R_{\mathfrak{m}}}(M_{\mathfrak{m}})\geq n$
for some $\mathfrak{m}$ by Proposition A.1. $\square$

\textbf{Proposition A.3. }\emph{For $M$ any finite $R$-module, the
equality\[
\mathrm{projdim}_{R}(M)=\mathrm{sup}_{\mathfrak{m}\in\mathrm{mSupp}(M)}\mathrm{sup}\left\{ i|\mathrm{Tor}_{i}^{R}(M,R/\mathfrak{m})\neq0\right\} \]
holds. If furthermore $\mathrm{projdim}_{R}(M)<\infty$ then the equality
\[
\mathrm{projdim}_{R}(M)=\mathrm{sup}\left\{ i|\mathrm{Ext}_{R}^{i}(M,R)\neq0\right\} \]
holds as well.}

\emph{Proof. }The module $\mathrm{Tor}_{i}^{R}(M,R/\mathfrak{m})$
is a finite-dimensional $R/\mathfrak{m}$-vector space, so localization
at $\mathfrak{m}$ leaves it unchanged, yielding\begin{eqnarray*}
\mathrm{Tor}_{i}^{R}(M,R/\mathfrak{m}) & \simeq & \mathrm{Tor}_{i}^{R}(M,R/\mathfrak{m})_{\mathfrak{m}}\\
 & \simeq & \mathrm{Tor}_{i}^{R_{\mathfrak{m}}}(M_{\mathfrak{m}},R_{\mathfrak{m}}/\mathfrak{m}).\end{eqnarray*}
Since the equality $\mathrm{projdim}_{S}(N)=\mathrm{sup}\left\{ i|\mathrm{Tor}_{i}^{S}(N,S/\mathfrak{m}_{S})\neq0\right\} $
holds for any local ring $S$ and any finite $S$-module $N$ (see
e.g. Lemma 19.1.ii of \cite{Matcrt}), the first claim now follows
from Proposition A.2.

For the second claim, we first note that if $S$ is a local ring and
$N$ is a finite $S$-module with $\mathrm{projdim}_{S}(N)<\infty$,
then $\mathrm{projdim}_{S}(N)=\mathrm{sup}\{i|\mathrm{Ext}_{S}^{i}(N,S)\ne0\}$
by Lemma 19.1.iii of \cite{Matcrt}. Hence by Proposition A.2 we have\begin{eqnarray*}
\mathrm{projdim}_{R}(M) & = & \mathrm{sup}_{\mathfrak{m}\in\mathrm{mSupp}(M)}\mathrm{sup}\left\{ i|\mathrm{Ext}_{R_{\mathfrak{m}}}^{i}(M_{\mathfrak{m}},R_{\mathfrak{m}})\neq0\right\} \\
 & = & \mathrm{sup}\{i|\mathrm{Ext}_{R}^{i}(M,R)_{\mathfrak{m}}\neq0\,\mathrm{for\, some}\,\mathfrak{m}\}\\
 & = & \mathrm{sup}\{i|\mathrm{Ext}_{R}^{i}(M,R)\neq0\},\end{eqnarray*}
as desired. $\square$

\textbf{Proposition A.4. }\emph{If $R$ is a Cohen-Macaulay ring,
$M$ is a finite $R$-module of finite projective dimension, and $\mathfrak{p}$
is an associated prime of $M$, then $\mathrm{ht}\mathfrak{p}=\mathrm{projdim}_{R_{\mathfrak{p}}}(M_{\mathfrak{p}})$.
In particular, $\mathrm{ht}\mathfrak{p}\leq\mathrm{projdim}_{R}(M)$.}

\emph{Proof. }Supposing $\mathfrak{p}$ is an associated prime of
$M$, there is an injection $R/\mathfrak{p}\hookrightarrow M$; this
localizes to an injection $R_{\mathfrak{p}}/\mathfrak{p}\hookrightarrow M_{\mathfrak{p}}$,
so $\mathrm{depth}_{R_{\mathfrak{p}}}(M_{\mathfrak{p}})=0$. Now we
compute\begin{eqnarray*}
\mathrm{ht}\mathfrak{p} & = & \mathrm{dim}(R_{\mathfrak{p}})\\
 & = & \mathrm{depth}_{R_{\mathfrak{p}}}(R_{\mathfrak{p}})\;(\mathrm{by\, the\, CM\, assumption})\\
 & = & \mathrm{depth}_{R_{\mathfrak{p}}}(M_{\mathfrak{p}})+\mathrm{projdim}_{R_{\mathfrak{p}}}(M_{\mathfrak{p}})\;(\mathrm{by\, the\, Auslander}-\mathrm{Buchsbaum\, formula)}\\
 & = & \mathrm{projdim}_{R_{\mathfrak{p}}}(M_{\mathfrak{p}}),\end{eqnarray*}
whence the result. $\square$

Now we single out an especially nice class of modules, which are equidimensional
in essentially every sense of the word. Recall the \emph{grade }of
a module $M$, written $\mathrm{grade}_{R}(M)$, is the $\mathrm{ann}_{R}(M)$-depth
of $R$; by Theorems 16.6 and 16.7 of \cite{Matcrt}, \[
\mathrm{grade}_{R}(M)=\mathrm{inf}\{i|\mathrm{Ext}_{R}^{i}(M,R)\neq0\},\]
so quite generally $\mathrm{grade}_{R}(M)\leq\mathrm{projdim}_{R}(M)$.

\textbf{Definition A.5. }\emph{A finite $R$-module $M$ is }perfect
\emph{if $\mathrm{grade}_{R}(M)=\mathrm{projdim}_{R}(M)<\infty$.}

\textbf{Proposition A.6. }\emph{Let $R$ be a Noetherian ring, and
let $M$ be a perfect $R$-module, with $\mathrm{grade}_{R}(M)=\mathrm{projdim}_{R}(M)=d$.
Then for any $\mathfrak{p}\in\mathrm{Supp}(M)$ we have $\mathrm{grade}_{R_{\mathfrak{p}}}(M_{\mathfrak{p}})=\mathrm{projdim}_{R_{\mathfrak{p}}}(M_{\mathfrak{p}})=d$.
If furthermore $R$ is Cohen-Macaulay, then $M$ is Cohen-Macaulay
as well, and every associated prime of $M$ has height $d$.}

\emph{Proof. }The grade of a module can only increase under localization
(as evidenced by the Ext definition above), while the projective dimension
can only decrease; on the other hand, $\mathrm{grade}_{R}(M)\leq\mathrm{projdim}_{R}(M)$
for any finite module over any Noetherian ring. This proves the first
claim.

For the second claim, Theorems 16.6 and 17.4.i of \cite{Matcrt} combine
to yield the formula\[
\mathrm{dim}(M_{\mathfrak{p}})+\mathrm{grade}_{R_{\mathfrak{p}}}(M_{\mathfrak{p}})=\mathrm{dim}(R_{\mathfrak{p}})\]
for any $\mathfrak{p}\in\mathrm{Supp(}M)$. The Auslander-Buchsbaum
formula reads\[
\mathrm{depth}_{R_{\mathfrak{p}}}(M_{\mathfrak{p}})+\mathrm{projdim}_{R_{\mathfrak{p}}}(M_{\mathfrak{p}})=\mathrm{depth}_{R_{\mathfrak{p}}}(R_{\mathfrak{p}}).\]
But $\mathrm{dim}(R_{\mathfrak{p}})=\mathrm{depth}_{R_{\mathfrak{p}}}(R_{\mathfrak{p}})$
by the Cohen-Macaulay assumption, and $\mathrm{grade}_{R_{\mathfrak{p}}}(M_{\mathfrak{p}})=\mathrm{projdim}_{R_{\mathfrak{p}}}(M_{\mathfrak{p}})$
by the first claim. Hence $\mathrm{depth}_{R_{\mathfrak{p}}}(M_{\mathfrak{p}})=\mathrm{dim}(M_{\mathfrak{p}})$
as desired. The assertion regarding associated primes is immediate
from the first claim and Proposition A.4. $\square$

\section{The dimension of irreducible components}

\begin{center}
{\large by James Newton}
\par\end{center}{\large \par}

In this appendix we use the results of the above article to give some
additional evidence for Conjecture 1.1.4. In the notation and terminology
of Section 1 above, we prove

\textbf{Proposition B.1. }\emph{Any irreducible component of $\mathscr{X}_{\mathbf{G},K^{p}}$
containing a given point $x$ has dimension at least $\mathrm{dim}(\scw_{K^{p}})-l(x)$.}

Note that Proposition 5.7.4 of \cite{UrEigen} implies that at least
one of these components has dimension at least $\mathrm{dim}(\scw_{K^{p}})-l(x)$.
This is stated without proof in that reference, and is due to G. Stevens
and E. Urban. We learned the idea of the proof of this result from
E. Urban --- in this appendix we adapt that idea and make essential
use of Theorem 3.3.1 (in particular the \textquoteleft{}Tor spectral
sequence\textquoteright{}) to provide a fairly simple proof of Proposition
B.1.

We place ourselves in the setting of Sections 1 and 4.3. In particular,
$\mathbf{G}$ is a reductive group over $\mathbf{Q}$. Fix an open
compact subgroup $K^{p}\subset\mathbf{G}(\mathbf{A}_{f}^{p})$ and
a slope datum $(U_{t},\Omega,h)$. Suppose that $\mathfrak{M}$ is
a maximal ideal of $\mathbf{T}_{\Omega,h}(K^{p})$ corresponding to
a point $x\in\scx_{\mathbf{G},K^{p}}(\overline{\mathbf{Q}_{p}})$.
Denote by $\mathfrak{m}$ the contraction of $\mathfrak{M}$ to $\sco(\Omega)$.
Let $\mathscr{P}$ be a minimal prime of $\mathbf{T}_{\Omega,h}(K^{p})$
contained in $\mathfrak{M}$. Since $H^{\ast}(K^{p},\scd_{\Omega})_{\leq h}$
is a finite faithful $\mathbf{T}_{\Omega,h}(K^{p})$-module, minimal
primes of $\mathbf{T}_{\Omega,h}(K^{p})$ are in bijection with minimal
elements of \[
\mathrm{Supp}_{\mathbf{T}_{\Omega,h}(K^{p})}(H^{\ast}(K^{p},\scd_{\Omega})_{\leq h});\]
by Theorem 6.5 of \cite{Matcrt}, minimal elements of the latter set
are in bijection with minimal elements of \[
\mathrm{Ass}_{\mathbf{T}_{\Omega,h}(K^{p})}(H^{*}(K^{p},\scd_{\Omega})_{\le h}).\]

\textbf{Definition B.2.} \emph{Denote by $r$ the minimal index $i$
such that $\mathscr{P}$ is in the support of $H^{i}(K^{p},\scd_{\Omega})_{\le h,\mathfrak{M}}$,
and by $q$ the minimal index $i$ such that $H^{i}(K^{p},\scd_{\lambda_{x}})_{(\ker\phi_{x})}\neq0$.}

Let $\wp$ denote the contraction of $\mathscr{P}$ to a prime of
$\sco(\Omega)_{\mathfrak{m}}$; in particular, $\wp$ is an associated
prime of $H^{r}(K^{p},\scd_{\Omega})_{\le h,\mathfrak{M}}$. The ring
$\sco(\Omega)_{\mathfrak{m}}$ is a regular local ring. The localisation
$\sco(\Omega)_{\wp}$ is therefore a regular local ring, with maximal
ideal $\wp\sco(\Omega)_{\wp}$. We let $(x_{1},...,x_{d})$ denote
a regular sequence generating $\wp\sco(\Omega)_{\wp}$. After multiplying
the $x_{i}$ by units in $\sco(\Omega)_{\wp}$, we may assume that
the $x_{i}$ are in $\sco(\Omega)$. Note that $(x_{1},...,x_{d})\sco(\Omega)_{\mathfrak{m}}$
may be a proper submodule of $\wp$. Nevertheless, we have \[
d=\mathrm{dim}(\sco(\Omega)_{\wp})=\mathrm{ht}(\wp).\]
 We will show that $d\le l(x)$.

Denote by $A_{i}$ the quotient $\sco(\Omega)_{\wp}/(x_{1},...,x_{i})\sco(\Omega)_{\wp}$
and denote by $\Sigma_{i}$ the Zariski closed subspace of $\Omega$
defined by the ideal $(x_{1},...,x_{i})\sco(\Omega)$. The affinoids
$\Sigma_{i}$ may be non-reduced. Note that $A_{i}=\sco(\Sigma_{i})_{\wp}$
and $\sco(\Sigma_{i+1})=\sco(\Sigma_{i})/x_{i+1}\sco(\Sigma_{i})$.

\textbf{Lemma B.3. }\emph{The space \[
H^{r-d}(K^{p},\scd_{\Sigma_{d}})_{\le h,\mathscr{P}}\]
is non-zero.} 

\emph{Proof. }By induction, it suffices to prove the following: let
$i$ be an integer satisfying $0\le i\le d-1$. Suppose \[
H^{r-i}(K^{p},\scd_{\Sigma_{i}})_{\le h,\mathscr{P}}\]
 is a non-zero $A_{i}$-module, with $\wp A_{i}$ an associated prime,
and \[
H^{t}(K^{p},\scd_{\Sigma_{i}})_{\le h,\mathscr{P}}=0\]
 for every $t<r-i$. Then \[
H^{r-i-1}(K^{p},\scd_{\Sigma_{i+1}})_{\le h,\mathscr{P}}\]
 is a non-zero $A_{i+1}$-module, with $\wp A_{i+1}$ an associated
prime, and \[
H^{t}(K^{p},\scd_{\Sigma_{i+1}})_{\le h,\mathscr{P}}=0\]
 for every $t<r-i-1$.

Note that the hypothesis of this claim holds for $i=0$, by the minimality
of $r$. Suppose the hypothesis is satisfied for $i$. It will suffice
to show that 
\begin{itemize}
\item $H^{t}(K^{p},\scd_{\Sigma_{i+1}})_{\le h,\mathscr{P}}=0$ for every
$t<r-i-1$ 
\item there is an isomorphism of non-zero $A_{i}$-modules \[
\iota:\mathrm{Tor}_{1}^{A_{i}}(H^{r-i}(K^{p},\scd_{\Sigma_{i}})_{\le h,\mathscr{P}},A_{i}/x_{i+1}A_{i})\cong H^{r-i-1}(K^{p},\scd_{\Sigma_{i+1}})_{\le h,\mathscr{P}}.\]
 
\end{itemize}
Indeed, the left hand side (which we denote by $T$) of the isomorphism
$\iota$ is given by the $x_{i+1}$-torsion in $H^{r-i}(K^{p},\scd_{\Sigma_{i}})_{\le h,\mathscr{P}}$,
so a non-zero $A_{i}$-submodule of $H^{r-i}(K^{p},\scd_{\Sigma_{i}})_{\le h,\mathscr{P}}$
with annihilator $\wp A_{i}$ immediately gives a non-zero $A_{i+1}$-submodule
of $T$ with annihilator $\wp A_{i+1}$.

Both the claimed facts are shown by studying the localisation at $\mathscr{P}$
of the spectral sequence \[
E_{2}^{s,t}:\mathrm{Tor}_{-s}^{A(\Sigma_{i})}(H^{t}(K^{p},\scd_{\Sigma_{i}})_{\le h},A(\Sigma_{i+1}))\Rightarrow H^{s+t}(K^{p},\scd_{\Sigma_{i+1}})_{\le h}\]
(cf. Remark 3.3.2). After localisation at $\mathscr{P}$, the spectral
sequence degenerates at $E_{2}$. This is because we have a free resolution
\[
0\rightarrow A_{i}\overset{\times x_{i+1}}{\rightarrow}A_{i}\rightarrow A_{i+1}\rightarrow0\]
of $A_{i+1}$ as an $A_{i}$-module (we use the fact that $x_{i+1}$
is not a zero-divisor in $A_{i}$), so $(E_{2}^{s,t})_{\mathscr{P}}$
vanishes whenever $s\notin\{-1,0\}$. Moreover, since \[
H^{t}(K^{p},\scd_{\Sigma_{i}})_{\le h,\mathscr{P}}=0\]
for every $t<r-i$, we know that $(E_{2}^{s,t})_{\mathscr{P}}$ vanishes
for $t<r-i$. The existence of the isomorphism $\iota$ and the desired
vanishing of $H^{t}(K^{p},\scd_{\Sigma_{i+1}})_{\le h,\mathscr{P}}$
are therefore demonstrated by the spectral sequence, since the only
non-zero term $(E_{2}^{s,t})_{\mathscr{P}}$ contributing to $(E_{\infty}^{r-i-1})_{\mathscr{P}}$
is given by $s=-1,t=r-i$, whilst $(E_{2}^{s,t})_{\mathscr{P}}=0$
for all $(s,t)$ with $s+t<r-i-1$. $\square$

\textbf{Corollary B.4.}\emph{ We have an inequality $r-d\ge q$. Since
$r\le q+l$ we obtain $d\le l$. In particular $\wp$ has height $\le l$,
so the irreducible component of $\mathbf{T}_{\Omega,h}(K^{p})$ corresponding
to $\mathscr{P}$ has dimension $\ge$ $\mathrm{dim}(\Omega)-l$.} 

\emph{Proof. }It follows from Proposition 4.5.2 (with $\Omega$ replaced
by $\Sigma_{d}$) that \[
H^{i}(K^{p},\scd_{\Sigma_{d}})_{\le h,\mathfrak{M}}\]
is zero for $i<q$. Our Lemma therefore implies that $r-d\ge q$.
The conclusion on dimensions follows from the observation made in
Section 4.5 that $\mathbf{T}_{\Omega,h}(K^{p})/\mathscr{P}$ has the
same dimension as $\sco(\Omega)/\wp$. $\square$ 

Proposition B.1 follows immediately from the Corollary. We have also
shown that if $d=l$, then $r=q+l$.

\bibliographystyle{amsalpha}
\nocite{*}
\def\cprime{$'$}
\providecommand{\bysame}{\leavevmode\hbox to3em{\hrulefill}\thinspace}
\providecommand{\MR}{\relax\ifhmode\unskip\space\fi MR }
\providecommand{\MRhref}[2]{%
  \href{http://www.ams.org/mathscinet-getitem?mr=#1}{#2}
}
\providecommand{\href}[2]{#2}

\end{document}